\documentclass[finalsize,finaltitle,draft]{zaamodified}

\DeclareMathOperator{\const}{const}

\newcommand{\ef}{\end{equation}}
\chardef\bslash=`\\ 

\theoremstyle{theorem}\newtheorem{corl}{Corollary}[theorem]
\theoremstyle{theorem}

\theoremstyle{definition}

\newcommand{\thmref}[1]{Theorem~\ref{#1}}
\newcommand{\secref}[1]{Section~\ref{#1}}

\newcommand{\lemref}[1]{Lemma~\ref{#1}}
\newcommand{\corref}[1]{Corollary~\ref{#1}}
\newcommand{\exampref}[1]{Example~\ref{#1}}

\renewcommand{\k}{\varkappa}

 \renewcommand{\sectionmark}[1]{}

\newcommand{\ve}{\varepsilon}

\newcommand{\iy}{\infty}

 \newcommand{\cupl}{\operatornamewithlimits{\bigcup}\limits}

\newcommand{\fc}{\frac}
\newcommand{\g}{\gamma}

 \newcommand{\dl}{\delta}
\newcommand{\Dl}{\Delta}

\newcommand{\om}{\omega}

\newcommand{\vp}{\varphi}

\newcommand{\doe}{\overset{\text{def}}{=}}
 
\newcommand{\loc} {\operatorname{loc}}

\newcommand{\D} {\Delta}

\renewcommand{\a}{\alpha}

\newcommand{\n}{\noindent}
\newcommand{\s}{\smallskip}
\newcommand{\m}{\medskip}

\begin{document}

\title{Conditions for Correct Solvability \\ of a Simplest
Singular Boundary Value Problem \\ of General Form}
\runtitle{Conditions for Correct Solvability}

\author{N.~A. Chernyavskaya and L.~A. Shuster}
\runauthor{N.~A. Chernyavskaya and L.~A. Shuster}

\address{A. Chernyavskaya:  Department of Mathematics and Computer Science,
Ben-Gurion University of the Negev,
P.O.B. 653, Beer-Sheva, 84105, Israel; nina@math.bgu.ac.il}
\address{A. Shuster:
Department of Mathematics and Computer Science,
Bar-Ilan University, 52900 Ramat Gan, Israel;}

\abstract{We consider the singular boundary value problem
\begin{align*} 
-r(x)y'(x)+q(x)y(x)&=f(x),\quad x\in \mathbb R \\[0.2cm]
\lim_{|x|\to\iy}y(x)&=0,
\end{align*}
where $f \in L_p(\mathbb R),$\ $p\in[1,\iy]$ $(L_\iy(\mathbb R):=C(\mathbb R)),
$\ $r$ is a continuous positive function
for $x \in \mathbb R$, \ $q \in L_1^{\loc}(\mathbb R), q \ge 0.$
A solution of this problem is, by definition, any absolutely continuous
function $y$ satisfying the limit condition and almost everywhere the
differential equation.
This problem  is called correctly solvable in a given space
$L_p(\mathbb R)$ if for
any function $f \in L_p(\mathbb R)$ it has a unique solution
$y \in L_p(\mathbb R)$ and
if the following inequality holds with an absolute constant $c_p\in (0,\iy):$
$$
\|y\|_{L_p(\mathbb R)}\le c_p\|f\|_{L_p(\mathbb R)},\quad \ f \in L_p(\mathbb R).
$$
We find minimal requirements for $r$ and $q$ under which the
above problem is correctly solvable in $L_p(\mathbb R).$}




\maketitle


\section{Introduction}\label{introduction}
\setcounter{equation}{0} \numberwithin{equation}{section}

We consider the singular boundary value problem
\begin{align}\label{1.1}
-r(x)y'(x)+q(x)y(x)&=f(x),\quad x\in \mathbb R  \\[0.2cm]
\label{1.2}
\lim_{|x|\to\iy}y(x)&=0.
\end{align}
Here and throughout the sequel $f \in L_p(\mathbb R),$\ $p\in[1,\iy]$\
$(L_\iy(\mathbb R):=C(\mathbb R))$
and
\begin{align}\label{1.3}
0<r &\in C^{\loc}(R) \\[0.1cm]
\label{1.4}
0\le q &\in L_1^{\loc}(R).
\end{align}
(In \eqref{1.3} we use the symbol $C^{\loc}(\mathbb R)$ to denote the set
of functions defined and
continuous on $\mathbb R.)$
Throughout the paper, we assume that the above conventions are satisfied.
We also define a solution of \eqref{1.1}--\eqref{1.2} as any absolutely
continuous function $y(\cdot)$ satisfying \eqref{1.2} and
satisfying \eqref{1.1} almost everywhere on $\mathbb R.$
Our goal is to establish a criterion for correct solvability of
\eqref{1.1}--\eqref{1.2} in $L_p(\mathbb R).$
We call problem \eqref{1.1}--\eqref{1.2} correctly solvable in a given space
$L_p(\mathbb R)$ (see
\cite[Ch.III, \S6, no.2]{Co}),  if both the following conditions hold:

\begin{tightitemize}
\item[\bf I)] For every function $f \in L_p(\mathbb R)$ there exists a
unique solution $y \in L_p(\mathbb R)$ of
\eqref{1.1}--\eqref{1.2}.
\item[\bf II)] The solution $y \in L_p(\mathbb R)$ of \eqref{1.1}--\eqref{1.2}
satisfies inequality
\eqref{1.5} with an absolute constant $c_p\in (0,\iy):$
\begin{equation}\label{1.5}
\|y\|_p\le c_p\|f\|_p,\quad \forall\ f\in L_p(\mathbb R).
\end{equation}
\end{tightitemize}
This problem is attractive for the following reasons.
First, \eqref{1.1}--\eqref{1.2} is a simplest singular boundary value problem
of general form
since this problem is one-dimensional, the order of the corresponding
differential
operation\linebreak
$\ell=-r(x)\fc{d}{dx}+q(x)$ is minimal possible, and, finally, the coefficient
at
the highest derivative in the operation $\ell$ is, in general, variable.
Hence our main argument for studying \eqref{1.1}--\eqref{1.2} is as follows:
Since singular
boundary value problems have been systematically studied for many years, it is
not natural to
ignore the simplest one.
Our second argument is closely related to the first one:  The main properties
of
\eqref{1.1}--\eqref{1.2} could serve as a useful illustration in exposition of
the general
theory of singular boundary value problems.
(Note that equation \eqref{1.1} was already used for that purpose (see
\cite[Ch.III]{D}, where
the systematic exposition of the theory of boundary value problems begins with
an exhaustive
description of regular boundary  value problems for equations \eqref{1.1} with
$r(x)\equiv1).$
To explain the third more specific reason for studying
\eqref{1.1}--\eqref{1.2}, we need some
additional details.
This paper is a continuation of the papers \cite{C}, \cite{CS1}, \cite{CS2}.
The main result of this study started in \cite{CS1}, \cite{CS2} was obtained in
\cite{C} and is stated in the following theorem.
\begin{theorem}\label{theorem1.1}
Let $r(x)\equiv 1,$\ $x\in \mathbb R.$
For $p\in[1,\iy)$ problem \eqref{1.1}--\eqref{1.2} is correctly solvable in
$L_p(\mathbb R)$ if and
only if \eqref{1.6} holds:
\begin{equation}\label{1.6}
\exists\ a\in (0,\iy):\quad \inf_{x\in \mathbb R}\int_{x-a}^{x+a} q(t)dt>0.
\end{equation}
Moreover, problem \eqref{1.1}--\eqref{1.2} is correctly solvable in $C(\mathbb R)$ if
and only if
\eqref{1.7} holds:
\begin{equation}\label{1.7}
\lim_{|x|\to\iy}\int_{x-a}^{x+a}q(t)dt=\iy,\quad\forall\ a\in (0,\iy).
\end{equation}
\end{theorem}

(Note that condition \eqref{1.7} can be easily rewritten in a slightly
different equivalent form (see
\eqref{6.2}). Let us briefly analyze \eqref{1.6} (the difference between
criteria \eqref{1.6} and
\eqref{1.7} can be explained only by the difference between the integral and
the uniform
metric in $L_p(\mathbb R)$ and $C(\mathbb R),$ respectively).
First note that condition \eqref{1.6} does not depend on $p\in [1,\iy),$ and
therefore  problem
\eqref{1.1}--\eqref{1.2} with $r(x)\equiv 1, x \in \mathbb R$ is correctly
solvable (or unsolvable) in $L_p(\mathbb R)$
for all  $p\in[1,\iy).$
Let us now look for properties of $q(\cdot)$ guaranteeing that condition
\eqref{1.6} holds (or
does not hold).

\s
Clearly, under the following condition \eqref{1.8} (resp., \eqref{1.9})
\begin{equation}\label{1.8}\inf_{x\in \mathbb R}q(x)>0
\end{equation}
\begin{equation}\label{1.9}
\lim_{x\to-\iy} q(x)=0\quad\text{or}\quad \lim_{x\to\iy}q(x)=0
\end{equation}
inequality \eqref{1.6} holds (resp., does not hold).
Moreover, in the cases where both \eqref{1.8} and \eqref{1.9} do not hold,
condition
\eqref{1.6} is in particular a requirement for the type of oscillation of
$q(\cdot).$
(We say that a non-negative function $q(\cdot)$ oscillates on
$\mathbb R$ if it vanishes on
a non-bounded subset of $\mathbb R$).
For example, problem \eqref{1.1}-\eqref{1.2} with
\begin{equation}\label{1.10}
r(x)\equiv 1,\quad q(x)=1+\sin (|x|^\theta),\quad x\in \mathbb R,\; \theta\in(0,\iy)
\end{equation}
is correctly solvable in $L_p(\mathbb R) $ if and only if $\theta\ge1$
(see \cite{C}, \cite{CS2}).

\s
Here the difference in oscillation properties of $q(\cdot)$ becomes especially
obvious if one
compares the graphs of $q(\cdot)$ for $\theta\ge 1$ and $\theta<1.$
Thus, if $r(x)\equiv 1,$ $x\in \mathbb R$, one restricts himself to the cases
\eqref{1.8} and \eqref{1.9} and
to the oscillatory function $q(\cdot)$, then the situation of
\eqref{1.9} is the only one where
\eqref{1.1}--\eqref{1.2} has no ``a priori chance" to be correctly solvable in
$L_p(\mathbb R).$
Therefore our third argument for studying \eqref{1.1}--\eqref{1.2} with
$r(x)\not\equiv\const,$ $x\in \mathbb R$, is the following.
We want to answer the following questions: For
what equations \eqref{1.1} with $r(x)\not\equiv\const,$\ $x\in \mathbb R$,
does the validity of \eqref{1.9} imply that
\eqref{1.1}--\eqref{1.2} is not correctly solvable in $L_p(\mathbb R)?$
Does the class of correctly solvable problems \eqref{1.1}--\eqref{1.2} include
equations
\eqref{1.1} with the property
\begin{equation}\label{1.11}
\lim_{x\to-\iy} r(x)=0\quad\text{or}\quad \lim_{x\to\iy}r(x)=0 \; ?
\end{equation}
(Concerning \eqref{1.9}, see \thmref{Theorem2.7} in \secref{MResults},
\thmref{theorem9.1} in Section 9 and \thmref{theorem1.2}
below. As to \eqref{1.11}, see the example in \secref{Example} and
\thmref{theorem1.2} below.)
We now note the main feature of problem \eqref{1.1}--\eqref{1.2}.
Its solution $y(\cdot),$ if it exists, can be explicitly expressed in terms of
initial data, and it has
the following form (see \lemref{lemma4.1} in \secref{PMR}):
\begin{equation}\label{1.12}
y(x)=(Gf)(x)\doe
\int_x^\iy\fc{1}{r(t)}\exp\bigg(-\int_x^t\fc{q(s)}{r(s)}ds\bigg)
f(t)dt,\quad x\in \mathbb R.
\end{equation}

Boundary value problems for differential equations of order $n\ge2$ do not have
such a property.
Thus, in order to study correct solvability of \eqref{1.1}--\eqref{1.2} in
$L_p(\mathbb R),$  one has to
find minimal requirements to $r(\cdot)$ and $q(\cdot)$ under which the integral
operator $G: L_p(\mathbb R)\to L_p(\mathbb R)$ is bounded, and,
in addition, all the functions from its image vanish on $\pm\iy.$
We emphasize that conditions for an integral operator of the form \eqref{1.12}
to be bounded in
$L_p(\mathbb R)$ may be easily found by means of Hardy-type inequalities (see \cite{OK}
 and Theorems \ref{theorem3.1} and \ref{theorem3.2} in
\secref{Auxiliary}). Therefore our goal reduces to finding a criterion for
the validity
of \eqref{1.2} which would serve as a completion to these conditions.
We specially note that such a restriction of goals of our investigation does
not make the
problem of studying \eqref{1.1}--\eqref{1.2} less important.
This a priori statement is based on our papers \cite{CS3}, \cite{CS4}, where it
is
shown (in the example of linear differential equations of order 2) that the
study of
boundary properties of solutions of singular differential equations of order 2
is a
problem separate from the problem of correct solvability of a differential
equation
in a given space, and thus it requires additional study.

\s
Let us now briefly describe our results (see \secref{MResults}).
Correct solvability of \eqref{1.1}--\eqref{1.2} in $L_p(\mathbb R)$ for
$p\in(1,\iy),$\ $p=1$ and
$p=\iy$, is studied separately since the known precise statements  on estimating
the norm of the
operator $G: L_p(\mathbb R)\to L_p(\mathbb R)$ (see \eqref{1.12}) can be divided into exactly
three cases (see
\secref{Auxiliary}).
Thus in \secref{MResults} we give three groups of theorems, each of which
contains a general
(unconditional) criterion for correct solvability  of \eqref{1.1}--\eqref{1.2}.
Theorems \ref{Theorem2.1}, \ref{Theorem2.3} and \ref{Theorem2.5} and a
particular criterion
of Theorems \ref{Theorem2.2}, \ref{Theorem2.4} and \ref{Theorem2.6}, which can
be applied
under a certain a priori assumption on $r(\cdot)$ and $q(\cdot)$ (see
\secref{MResults}, \eqref{2.11}).
Theorems \ref{Theorem2.1}, \ref{Theorem2.3} and \ref{Theorem2.5} contain
conditions expressed in
a non-local form, and it might be hard to check them in concrete cases.
Statements of Theorems \ref{Theorem2.2}, \ref{Theorem2.4} and \ref{Theorem2.6}
obtained with the
help of general criteria contain only conditions expressed in a local form.
It is much easier to check them in concrete cases (compare with the conditions
of
Theorems \ref{Theorem2.1}, \ref{Theorem2.3} and \ref{Theorem2.5}) since one can
use standard
tools of local analysis (see the example in \secref{Example}).

\s
In addition, following a suggestion of S. Luckhaus,   we complement
Theorems
\ref{Theorem2.1}--\ref{Theorem2.6} by  \thmref{theorem1.2}  where
the conditions for correct solvability of problem \eqref{1.1}--\eqref{1.2} in
the space $L_p(\mathbb R),$\ $p\in[1,\iy]$, are expressed only in terms
of the coefficients $r(\cdot)$ and $q(\cdot)$ of equation
\eqref{1.1}. This theorem can be viewed as an explicit example of the
relationship between $r(\cdot)$,
$q(\cdot)$ and the parameter $p\in [1,\iy],$ fixing the choice of the space
$L_p(\mathbb R)$, which guarantees
the correct solvability of the problem under investigation.

\begin{theorem}\label{theorem1.2} {\rm (\S$9$)}
Suppose that the following conditions hold:
\begin{tightitemize}
\item[\bf 1)] The functions $r(x)$ and $q(x)$ are positive and continuous for
$x\in \mathbb R$.
\item[\bf 2)] There exists $a\ge 1$ and $b>0$ and an interval $(\alpha,\beta)$ such
that for all $x\notin (\alpha,\beta)$
\begin{equation}\label{1.13}
\frac{1}{a}\le \frac{r(t)}{r(x)},\ \frac{q(t)}{q(x)}\le a
\qquad\text{\rm for}\; |t-x|\le b\,\frac{r(x)}{q(x)} .
\end{equation}
\item[\bf 3)]  The following inequality is valid:
\begin{equation}\label{1.14}
\gamma\doe a^2\exp\Big(-\frac{b}{a^2}\Big)\le \frac{1}{3}.
\end{equation}
\end{tightitemize}
Then problem \eqref{1.1}--\eqref{1.2} is correctly solvable in
$L_p(\mathbb R),$\ $p\in[1,\iy]$, if and only if
the conditions from the following table are satisfied:

\bigskip

\centerline{\vbox {\offinterlineskip
 \hrule
\halign{&\vrule
#&
\strut\hfil
#\tabskip= .02em
\cr
height4pt
 &\omit&&\omit&&\omit&& \omit&&\omit\hfill && \cr
&& \quad Space&& \quad$L_1(\mathbb R)$    \quad
&&\ $L_p(\mathbb R),\ 1<p<\infty\ $&& \quad $C(\mathbb R)$
  &&\cr
 \noalign{\hrule}
&\omit&&\omit&&\omit&& \omit&&\omit\hfill && \cr
 && \quad Conditions for\quad&&\quad$r_0>0$  && \qquad $\sigma_{p'}>0$\qquad
&&\quad\ $q(x)\to\infty$ \hfil &&\cr
&& \quad correct solvability\quad&&\quad$q_0>0$  && \qquad$q_0>0$\qquad
&&\quad as $|x|\to\infty$ &&
  \hfil &&\cr
&\omit&&\omit&&\omit&& \omit&&\omit\hfill && \cr }
  \hrule} \ }

\m\n
Here $p'=\frac{p}{p-1}$ for $ p\in(1,\iy)$ and
\begin{equation}\label{1.16}
r_0=\inf_{x\in \mathbb R} r(x),\quad  q_0=\inf_{x\in \mathbb R}q(x),
\quad \sigma_{p'}=\inf_{x\in R}r(x)^{\frac{1}{p}}q(x)^{\frac{1}{p'}},\
\end{equation}
\end{theorem}

\begin{remark}\label{rem1.3}
The simplicity of the conditions in the above table can be explained by the
fact that the function
$q(\cdot)$ in \thmref{theorem1.2} is assumed to be non-oscillating.
\end{remark}

It should be emphasized that in this work we use tools
developed by the authors when they studied boundary properties of solutions of
the
Sturm-Liouville equation (see \cite{CS3}, \cite{CS4}).
These tools originated under the direct influence of ideas and techniques known
as
Otelbaev's localization methods (which is the terminology of M.~O. Otelbaev).
We refer to \cite{MO}, \cite{O} for a detailed exposition of these tools
and effective applications
to various problems of analysis and the theory of differential and difference
operators.
We also note that one can find more detailed and concrete comments to all the
statements,
including technical ones, in the course of the proofs.

\subsection*{Acknowledgment}  The authors   thank  Prof. Z.~S.
Grinshpun and Prof. S. Luckhaus for useful remarks and
suggestions which significantly improved the paper.


\section{Main Results}\label{MResults}

In this section we give the main results of the paper concerning correct
solvability of
\eqref{1.1}--\eqref{1.2} in $L_p(\mathbb R),$\ $p\in[1,\iy]$, including some comments.
Here and throughout the sequel, the symbols $c,c(\cdot),c_1,c_2,\dots$ denote
absolute positive constants which are not
essential for exposition and may differ even within a single chain of
calculations.
Let us introduce an auxiliary function $d(\cdot)$ which will be systematically used
below.
We temporarily assume that in addition to \eqref{1.3}--\eqref{1.4}, we have
also \eqref{2.1}:
\begin{equation} \label{2.1}
S_1\doe\int_{-\iy}^0\fc{q(t)}{r(t)}dt\,=\,\iy.
\end{equation}
Then for every $x\in \mathbb R$, we have a uniquely determined function
$d(\cdot)$ where
\begin{equation}\label{2.2}
d(x)=\inf_{d>0}\left\{d:\int_{x-d}^{x+d}\fc{q(t)}{r(t)}dt=2\right\}
\end{equation}
(see \S\S$3$-$7$ for various properties of $d(\cdot)$).
Note that the function $d(\cdot)$ is mainly applied as ``counter"   of the exponent
in
the inner integral of formula \eqref{1.12}.

\begin{theorem}\label{Theorem2.1} {\rm (\S$4$)}
Let $p\in(1,\iy),$ \ $p'=p\,(p-1)^{-1}$.
problem \eqref{1.1}--\eqref{1.2} is correctly solvable in $L_p(\mathbb R)$
if and only if the following
conditions hold together:
\begin{align}
&{\bf 1)}\hspace{0.8cm}
M_p \doe \sup_{x\in \mathbb R}M_p(x)\,<\,\iy ,
\; \text{ where }  \label{2.3}\\
\begin{split}
&\phantom{{\bf 1)}\hspace{0.25cm}}
M_p(x)\doe
\bigg[\int_{-\iy}^x\exp\bigg(-p\int_t^x\fc{q(\xi)}{r(\xi)}
d\xi\bigg)dt\bigg]^{\frac{1}{p}} \hspace{4.cm}\\
&\hspace{2.1cm}\cdot
\bigg[\int_x^\iy\fc{1}{r(t)^{p'}}\exp
\bigg(-p'\int_x^t\fc{q(\xi)}{r(\xi)}d\xi\bigg)dt\bigg]^{\frac{1}{p'}}
\label{2.4}
\end{split} \\[0.1cm]
&{\bf 2)}\hspace{0.94cm} S_1=\iy \quad \text{\rm (see \eqref{2.1})}\label{2.5}\\
&{\bf 3)}\hspace{0.73cm}
A_{p'} \doe \sup_{x\in \mathbb R}A_{p'}(x)\,<\,\iy,
\; \text{ where } \label{2.6}\\
&\phantom{{\bf 1)}\  }A_{p'}(x) \doe \int_{x-d(x)}^{x+d(x)}\fc{dt}{r(t)^{p'}}\label{2.7}
\end{align}

\n
Condition \eqref{2.3} is a criterion for the operator
$G: L_p(\mathbb R)\to L_p(\mathbb R)$ {\rm (}see \eqref{1.12}{\rm )} to
be bounded, and the following inequalities hold:
\begin{equation}\label{2.8}
M_p\le\|G\|_{p\to p}\le (p)^{\frac{1}{p}}(p')^{\frac{1}{p'}}M_p,\quad p\in (1,\iy).
\end{equation}
\end{theorem}

\begin{corl}\label{Corollary2.1.1} {\rm (\S$4$)}
Let $p\in (1,\iy),$ and suppose that \eqref{1.1}--\eqref{1.2} is correctly
solvable in $L_p(\mathbb R).$
Then its solution $y(\cdot)$ satisfies the inequality
\begin{equation}\label{2.9}
\Big\|\Big(\fc{q}{r}\Big)^{\frac{1}{p}}y\Big\|_p\le
c(p)\|f\|_p \quad\forall \ f \in L_p(\mathbb R).
\end{equation}
\end{corl}

\begin{remark}\label{Remark2.1}
\thmref{Theorem2.1} is ``not convenient" for application to concrete problems
\eqref{1.1}--\eqref{1.2}.
Its main disadvantage is as follows: condition \eqref{2.3} is non-local, and
therefore may be
hard for checking.
In contrast with \eqref{2.3}, condition\eqref{2.6} is significantly easier to
check since
usually one can find  two-sided estimates sharp by order for the function
$d(\cdot)$ (see
\thmref{Theorem2.8} below and \secref{Example}).
Moreover, under an additional -- to \eqref{1.3}--\eqref{1.4} -- assumption
\begin{equation}\label{2.10}
d_0\doe\sup_{x\in \mathbb R}d(x)\,<\,\iy ,
\end{equation}
conditions \eqref{2.3} and \eqref{2.5} hold (see \thmref{Theorem2.2} below),
and therefore checking 1)~--~3) reduces to checking 3) only, and can be made by
local tools (see
\secref{Example}).
Thus application of \eqref{2.10} simplifies \thmref{Theorem2.1}.
Checking \eqref{2.10}, in turn, may be simplified because of the following
lemma.
\end{remark}

\begin{lemma}\label{lemmama2.1} {\rm (\S \ref{Auxiliary})}
The conditions $S_1=\iy$ and $d_0<\iy$ {\rm (}see \eqref{2.1}
and \eqref{2.10}{\rm )} hold
together if and only if
\begin{equation}\label{2.11}
\exists\ a\in(0,\iy): \quad  B(a)\doe\inf_{x\in \mathbb R}
\int_{x-a}^{x+a}\fc{q(t)}{r(t)}\,dt\,>\,0.
\end{equation}
\end{lemma}

In the next statement, conditions for correct solvability of
\eqref{1.1}--\eqref{1.2} in
$L_p(\mathbb R)$ contain only local requirements.

\begin{theorem}\label{Theorem2.2} {\rm (\S$4$)}
Let $p\in(1,\iy),$ and suppose that \eqref{2.11} holds.
Then \eqref{1.1}--\eqref{1.2} is correctly solvable in $L_p(\mathbb R)$
if and only if $A_{p'}<\iy$
{\rm (}see \eqref{2.6}{\rm )}.
It holds
$$
\|G\|_{p\to p}\le cd_0^{\frac{1}{p}}A_{p'}{}^{\frac{1}{p'}}.
$$
\end{theorem}

\begin{corl}\label{Corollary2.2.1} {\rm (\S$4$)}
Let $p\in (1,\iy),$ and suppose that
\begin{equation}\label{2.12}
\exists\ c\in [1,\iy):\quad c^{-1}\le r(x)\le c,\ x\in \mathbb R.
\end{equation}
Then \eqref{1.1}--\eqref{1.2} is correctly solvable in $L_p(\mathbb R)$ if and
only if \eqref{1.6} holds.
\end{corl}

In the Theorems \ref{Theorem2.3} and \ref{Theorem2.4} we give conditions for
correct solvability of
\eqref{1.1}--\eqref{1.2} in $L_1(\mathbb R).$

\begin{theorem}\label{Theorem2.3} {\rm (\S$5$)}
problem \eqref{1.1}--\eqref{1.2} is correctly solvable in $L_1(\mathbb R)$ if
and only if the
following conditions hold together:
\begin{align}
&{\bf 1)}\hspace{0.8cm}  S_1=\iy\quad\text{{\rm (}see \eqref{2.1}{\bf )}}
\hspace{7cm}
\label{2.13}\\
&{\bf 2)}\hspace{0.8cm}
r_0\doe\inf_{x\in \mathbb R}r(x)\,>\,0
\label{2.14}\\
\begin{split}
&{\bf 3)}\hspace{0.6cm} M_1\doe\sup_{x\in \mathbb R}M_1(x)\,<\,\iy,
\; \text{ where }\label{2.15} \\
&\hspace{0.42cm}M_1(x)\doe\fc{1}{r(x)}\int_{-\iy}^x\exp\bigg(-\int_t^x\fc{q(\xi)}{r(\xi)}d\xi\bigg)dt.
\end{split}
\end{align}
Condition \eqref{2.15} is a criterion for the operator
$G: L_1(\mathbb R)\to L_1(\mathbb R)$
to be bounded, and we have
\begin{equation}\label{2.16}
\|G\|_{1\to 1}=M_1.
\end{equation}
\end{theorem}

\begin{corl}\label{Corollary2.3.1} {\rm (\S$5$)}
Suppose that \eqref{1.1}--\eqref{1.2} is correctly solvable in $L_1(\mathbb R).$
Then its solution $y(\cdot)$ satisfies the following inequality
{\rm (}see \eqref{2.14}{\rm )}:
\begin{equation}\label{2.17}
\|y'\|_1+\Big\|\fc{q}{r}y\Big\|_1\le \fc{3}{r_0} \|f\|_1 \quad \forall\
f\in L_1(\mathbb R).
\end{equation}
\end{corl}
\begin{corl}\label{Corollary2.3.2} {\rm (\S$5$)}
Under assumption \eqref{2.12}, problem \eqref{1.1}--\eqref{1.2} is correctly
solvable in $L_1(\mathbb R)$ if
and only if \eqref{1.6} holds.
\end{corl}
\begin{remark}\label{Remark2.2}
Inequality \eqref{2.17} shows that on the set of solutions of
\eqref{1.1}--\eqref{1.2}, equation \eqref{1.1}
has the property of separability in the space $L_{\theta(\cdot)}^1(\mathbb R)$
where $\theta(x)=\fc{1}{r(x)},$\ $x\in \mathbb R.$
Note that  the problem of separability of singular differential operators was
first considered
in \cite{EG1}, \cite{EG2}.
\end{remark}

In the next theorem, conditions for correct solvability of
\eqref{1.1}--\eqref{1.2}
in $L_1(\mathbb R)$ contain only local requirements.

\begin{theorem}\label{Theorem2.4} {\rm (\S$5$)}
Suppose that condition \eqref{2.11} holds.
Then \eqref{1.1}--\eqref{1.2} is correctly solvable in $L_1(\mathbb R)$ if and only if
$r_0>0$ {\rm (}see \eqref{2.14}{\rm )}.
We have the inequality  $\|G\|_{1\to1}\le cr_0^{-1}d_0$
{\rm (}see \eqref{2.10}{\rm )}.
\end{theorem}
Below, in Theorems \ref{2.5} and \ref{2.6}, we give conditions for correct
solvability of
\eqref{1.1}--\eqref{1.2} in $C(\mathbb R).$

\begin{theorem}\label{Theorem2.5} {\rm (\S$6$)}
problem \eqref{1.1}--\eqref{1.2} is correctly solvable in $C(\mathbb R)$ if and
only if
\begin{equation}\label{2.18}
\lim_{|x|\to\iy}\int_x^\iy
\fc{1}{r(t)}\exp\bigg(-\int_x^t\fc{q(\xi)}{r(\xi)}d\xi\bigg)dt=0.
\end{equation}
Moreover, if \eqref{1.1}--\eqref{1.2} is correctly solvable in $C(\mathbb R),$
then $S_1=\iy$ {\rm (}see \eqref{2.1}{\rm )}.
Finally, it holds
\begin{equation}\label{2.19}
\|G\|_{C(\mathbb R)\to C(\mathbb R)}=\sup_{x\in \mathbb R}\int_x^\iy
\fc{1}{r(t)}\exp\bigg(-\int_x^t\fc{q(\xi)}{r(\xi)}d\xi\bigg)dt.
\end{equation}
\end{theorem}
\begin{corl}\label{Corollary2.5.1} {\rm (\S$6$)}
Under assumption \eqref{2.12}, problem \eqref{1.1}--\eqref{1.2} is correctly
solvable in $C(\mathbb R)$
if and only if \eqref{1.7} holds.
\end{corl}
A local form of the condition  for correct solvablility of
\eqref{1.1}--\eqref{1.2} is given in
the following theorem.
\begin{theorem}\label{Theorem2.6} {\rm (\S$6$)}
Under assumption \eqref{2.11}, problem \eqref{1.1}--\eqref{1.2} is correctly
solvable in $C(\mathbb R)$ if and only if $\tilde A_0=0,$ where $\tilde
A_0=\lim\limits_{|x|\to\iy}\tilde{A}(x).$ Here
\begin{equation}\label{2.20}
\tilde{A}(x)=\int_{x-d(x)}^{x+d(x)}\fc{dt}{r(t)}, \quad x \in R
\end{equation}
It holds
\begin{equation}\label{2.21}
\|G\|_{C(\mathbb R)\to C(\mathbb R)}\le c\sup_{x\in \mathbb R}
\int_{x-d(x)}^{x+d(x)}\fc{dt}{r(t)}.
\end{equation}
\end{theorem}

The next theorem answers the question from \secref{introduction} in connection
with condition
\eqref{1.9}.

\begin{theorem}\label{Theorem2.7} {\rm (\S$7$)}
Suppose that the conditions \eqref{2.11} and \eqref{1.9} hold.
Then for every $p\in[1,\iy]$, problem \eqref{1.1}--\eqref{1.2} cannot be
correctly solvable in
$L_p(\mathbb R).$
\end{theorem}

To conclude the section, we present one of the methods for obtaining two-sided
estimates for the
function $d(\cdot)$ (see \S\S$3$ and $8$ for additional details).

\begin{theorem}\label{Theorem2.8} {\rm (\S$3$)}
Suppose that there exist a continuous positive function $q_1(\cdot)$
on $\mathbb R$ and a function
$q_2 \in L_1^{\loc}(\mathbb R)$ such that
$q(x)=q_1(x)+q_2(x)$ for $x\in \mathbb R$ and, in
addition,
$\k_1(x)\to 0,$\ $\k_2(x)\to0$ as $|x|\to\iy,$
where
\begin{align}\label{2.22}
\k_1(x)&\doe\sup_{|z|\le 2
\fc{r(x)}{q_1(x)}}\left|\int_0^z
\left[\fc{q_1(x+t)}{r(x+t)}-2\fc{q_1(x)}{r(x)}+\fc{q_1(x-t)}
{r(x-t)}\right]dt\right| \\
\label{2.23}
\k_2(x)&\doe\sup_{|z|\le
2\fc{r(x)}{q_1(x)}}\left|\int_{x-z}^{x+z}\fc{q_2(t)}{r(t)} dt\right|.
\end{align}
Then it hold the relations
\begin{align}\label{2.24}
\fc{q_1(x)}{r(x)}d(x)=1+\ve(x),&\quad |\ve(x)|\le c_1[\k_1(x)+\k_2(x)],\quad
|x|\gg1 \\
\label{2.25}
c_2^{-1}&\le \fc{q_1(x)}{r(x)}d(x)\le c_2,\quad x\in \mathbb R.
\end{align}
\end{theorem}

\begin{remark}\label{rem2.9}
We refer to \S$9$ for some assertions related to \thmref{theorem1.2} (see~\S$1$) and
the complete formulations of the
Theorems~\ref{Theorem2.6} and \ref{Theorem2.7}.
\end{remark}


\section{Auxiliary Results and Technical Assertions}\label{Auxiliary}

In the first part of this section we present various assertions used in the
proofs of
Theorems~\ref{Theorem2.1}--\ref{Theorem2.7}. In the second part we prove some
of them.

\subsection{Statement of auxiliary results}\label{Statement}

\begin{theorem}\label{theorem3.1} {\rm \cite[Ch.I, \S$6$, no.6.5]{OK}}.
Let $AC(R)^{(+)}$ be the set of absolutely continuous functions $F(\cdot)$ on
$\mathbb R$ such that
\begin{equation}\label{3.1}
\lim_{x\to\iy} F(x)=0,
\end{equation}
and let $w(\cdot)$ and $v(\cdot)$ be measurable and positive functions on $R.$
Then the Hardy inequality
\begin{equation}\label{3.2}
\|w^{\frac{1}{p}}F\|_p\le C\|v^{\frac{1}{p}}F'\|_p
\end{equation}
holds for all $F\in AC(R)^{(+)}$ if and only if $H^{(p)}<\iy$,
where
\begin{align}\label{3.3}
\begin{split}
H^{(p)}&=\sup_{x\in \mathbb R}H^{(p)}(x)  \\
H^{(p)}(x)&=\left(\int_{-\iy}^xw(t)dt\right)^{\frac{1}{p}}
\left(\int_x^\iy
v(t)^{-\frac{p'}{p}}dt\right)^{\frac{1}{p'}},\quad x\in \mathbb R,
\end{split}
\end{align}
and $p\in (1,\iy),\ p'=\fc{p}{p-1}.$
Moreover, the following estimates hold for the smallest constant $C$ in
\eqref{3.2}:
\begin{equation}\label{3.4}
H^{(p)}\le C\le(p)^{\frac{1}{p}}(p')^{\frac{1}{p'}}H^{(p)}.
\end{equation}
\end{theorem}
\begin{theorem}\label{theorem3.2} {\rm \cite[Ch.2, \S7]{O}}
Let $p\in (1,\iy)$ and let $\mu( \cdot),\theta(\cdot)$ be continuous and
positive functions on $\mathbb R$. Denote by $K $ the integral operator
\begin{equation}\label{3.5}
(K  f)(t)=\mu(t)\int_t^\iy \theta(\xi)f(\xi)d\xi,\quad t\in \mathbb R.
\end{equation}
Then the operator $K : L_p(\mathbb R)\to L_p(\mathbb R)$ is bounded if and
only if $H_p<\iy$,
where
\begin{align}
\begin{split}
H_p&=\sup\limits_{x\in \mathbb R} H_p(x)  \\
\label{3.6}
H_p(x)&=\left[\int_{-\iy}^x \mu(t)^pdt\right]^{\frac{1}{p}}\cdot
\left[\int_x^\iy\theta(t)^{p'}dt\right]^{\frac{1}{p'}},\quad p'=\frac{p}{p-1}.
\end{split}
\end{align}
There hold the inequalities
\begin{equation}\label{3.7}
H_p\le\|K \|_{p\to p}\le (p)^{\frac{1}{p}}(p')^{\frac{1}{p'}}H_p.
\end{equation}
\end{theorem}

\begin{theorem} \label{theorem3.3} {\rm \cite[Ch.V, \S2, no.4-5]{KA}}
Let $K $ be the integral operator \eqref{3.5}.
Then
\begin{align}\label{3.11}
\|K \|_{1\to 1}&=\sup_{x\in \mathbb R}\theta(x)\int_{-\iy}^x\mu(t)dt \\
\label{3.12}
\|K \|_{C(R)\to C(R)}&=\sup_{x\in \mathbb R}\mu(x)\int_x^\iy\theta(t)dt.
\end{align}
\end{theorem}

\begin{lemma}\label{lemma3.1}
Let a function $\Phi(x,d)$ be defined for $x\in \mathbb R,$\ $d\in \mathbb R$
and satisfy the following properties:

\begin{tightitemize}[3mm]
\item[\bf 1)] $\Phi(x,d)$ is continuous in $d\in \mathbb R $ for every fixed
$x\in \mathbb R$
\item[\bf 2)] $\Phi(x,0)=0,$\ $\lim\limits_{d\to\iy}\Phi(x,d)=\iy,\
x\in \mathbb R$
\item[\bf 3)] $\Phi(x,d)$ does not decrease for $d\ge 0$ for every fixed
$x\in \mathbb R.$
\end{tightitemize}
Then for every $x\in \mathbb R$ the equation $\Phi(x,d)=1$ in $d\ge0$
has at least one solution. Let
\begin{equation}\label{3.13}
\hat d(x)=\inf_{d\ge0}\{d:\Phi(x,d)=1\}.
\end{equation}
If, in addition to {\rm 1)~--~3)}, the condition
\begin{equation}
\begin{aligned}\label{3.14}
&\text{{\bf 4)}} \quad
h\cdot[\Phi(x,d)-\Phi(x+h,d-h)]\ge0 \hspace{5cm}\\
&\phantom{{\bf 4)} \quad} h\cdot[\Phi(x+h,d+h)-\Phi(x,d)]\ge 0
\end{aligned}
\end{equation}
holds for every given $x\in \mathbb R$ and $|h|\le d\le \hat d(x)$,
then the following inequalities hold:
\begin{alignat}{2}\label{3.15}
&|\hat d(x+h)-\hat d(x)|\le|h| &\quad &\text{\rm for } |h|\le\hat d(x) \\
\label{3.16}
&1-\ve\le\fc{\hat d(t)}{\hat d(x)}\le 1+\ve &\quad &\text{\rm for }
|t-x|\le\ve\hat d(x),\; \ve\in[0,1].
\end{alignat}
Moreover, under the conditions {\rm 1)~--~4)}, the function $\hat d(\cdot)$
is continuous on $\mathbb R.$
\end{lemma}

\begin{example}\label{examp3.1} Let $S_1=\infty$ (see \eqref{2.1})
Set
$$\Phi(x,d)=\fc{1}{2}\int_{x-d}^{x+d}\fc{q(\xi)}{r(\xi)}d\xi,\quad x\in
R,\; d\in \mathbb R.$$
Then, clearly, the hypotheses 1)~--~4) of \lemref{lemma3.1} are satisfied, and
in this case we have
$\hat d(x)=d(x)$ for $x\in \mathbb R$ (see (2.2)).
In particular, $d(x)$ is continuous for $x\in \mathbb R$ and satisfies
\eqref{3.15}~--~\eqref{3.16}.
\end{example}
\begin{definition}\label{definition3.1}
Suppose that we are given:\ $x\in \mathbb R,$ a positive continuous function
$\varkappa(t),$ for $t\in \mathbb R$, a
sequence $\{x_n\}_{n\in N'},$\ $N'=\{\pm 1,\pm2,\dots\}.$
Consider the segments $\Dl_n=[\Dl_n^-,\Dl_n^+],$\ $\Dl_n^\pm=x_n\pm
\varkappa(x_n).$
We say that a sequence of segments $\{\Dl_n\}_{n=1}^\iy$\
$(\{\Dl_n\}_{n=-\iy}^{-1})$ forms an
$R(x,\varkappa(\cdot))$-{\it covering} of $[x,\iy)$ (resp., $(-\iy,x])$ if the
following conditions
hold:

\begin{tightitemize}
\item[1)] $\Dl_n^+=\Dl_{n+1}^-$ for $n\ge 1$ (resp., $\Dl_{n-1}^+=\Dl_n^-$
for $n\le -1)_{\phantom{\big)}}$
\item[2)] $\Dl_1^-=x$ (resp., $\Dl_{-1}^+=x)$, \
$\cupl\limits_{n\ge1}\Dl_n=[x,\iy)$
(resp.,
$\cupl\limits_{n=\le-1} \Dl_n=(-\iy,x])$.
\end{tightitemize}
\end{definition}

\begin{lemma}\label{lemma3.2}
Suppose that a positive continuous function $\varkappa(\cdot)$ on $\mathbb R$
satisfies the following relations:
\begin{equation}\label{3.17}
 \underset {t\to\iy}\varliminf (t-\varkappa(t))=\iy\quad
\text{(resp.,}\quad
 \underset{t\to-\iy}\varlimsup (t+\varkappa(t))=-\iy).
\end{equation}
Then for every $x\in \mathbb R$ there exists an $R(x,\varkappa(\cdot))$-covering of
$[x,\iy)$ (resp., an
$R(x,\varkappa(\cdot))$-covering of $(-\iy,x]).$
\end{lemma}

\begin{lemma}\label{lemma3.3}
Let $S_2=\iy$ (resp., $S_1=\iy$, see \eqref{2.1}), where
\begin{equation}\label{3.18}
S_2=\int_0^\iy\fc{q(t)}{r(t)}dt.
\end{equation}
Then for every $x\in \mathbb R$ there exists an $R(x,d(\cdot))$-covering of $[x,\iy)$
(resp., an
$R(x,d(\cdot))$-covering of $(-\iy,x]).$
\end{lemma}

\subsection{Proofs of auxiliary assertions}\label{Proofsaux}

\renewcommand{\qedsymbol}{}
\begin{proof}[Proof of Theorem 3.2] \ {\tt Necessity}.

\s
Let $p\in (1,\iy)$, and suppose that the operator
$K:L_p(\mathbb R)\to L_p(\mathbb R)$ is bounded.
Denote by $[t_1,t_2]$ the arbitrary finite interval and
$$
f_0(\xi)=\begin{cases}
\theta(\xi)^{p'-1},&\quad\text{if $\xi\in[t_1,t_2]$}\\
0,&\quad \text{if    $\xi\notin[t_1,t_2]$}.
\end{cases}
$$
From the continuity of $\theta(\cdot)$, it follows that $f_0\in L_p(\mathbb R)$
and
\begin{equation}\label{3.16a}
\|f_0\|_p^p=\int_{t_1}^{t_2}\theta(\xi)^{p(p'-1)}d\xi
=\int_{t_1}^{t_2}\theta(\xi)^{p'}d\xi.
\end{equation}
Moreover,
\begin{align} \label{3.17a}
\begin{split}
\|Kf_0\|_p^{p}&=
\int_{-\iy}^{t_1}\left[\mu(t)\int_t^\iy\theta(\xi)f_0(\xi)d\xi\right]^pdt \\
&\phantom{\ge }
+\int_{t_1}^\iy\left[\mu(t)\int_t^\iy\theta(\xi)f_0(\xi)d\xi\right]^pdt \\
&\ge\int_{-\iy}^{t_1}\left[\mu(t)\int_t^\iy\theta(\xi)f_0(\xi)d\xi\right]^pdt\\
&\ge\int_{-\iy}^{t_1}\mu(t)^pdt\left(\int_{t_1}^\iy\theta(\xi)f_0(\xi)d\xi\right)^p\\
&=\int_{-\iy}^{t_1}\mu(t)^pdt\left(\int_{t_1}^{t_2}\theta(\xi)^{p'}d\xi\right)^p.
\end{split}
\end{align}
Now according to \eqref{3.16a} and \eqref{3.17a} get:
\begin{align*}
\left(\int_{-\iy}^{t_1}\mu(t)^pdt\right)^{\frac{1}{p}}
\left(\int_{t_1}^{t_2}\theta(\xi)^{p'}d\xi\right) &\le\|Kf_0\|_p \\[0.2cm]
&\le\|K\|_{p\to p}\|f_0\|_p\\
&  =\|K\|_{p\to p}
\left(\int_{t_1}^{t_2}\theta(\xi)^{p'}d\xi\right)^{\frac{1}{p}},
\end{align*}
and hence
\begin{align}\label{3.18a}
\|K\|_{p\to p}\ge
\left(\int_{-\iy}^{t_1}\mu(t)^pdt\right)^{\frac{1}{p}}
\left(\int_{t_1}^{t_2}\theta(\xi)^{p'}d\xi\right)^{\frac{1}{p'}}.
\end{align}
Since $t_1,$ $t_2$ are arbitrary, from \eqref{3.18a} we obtain the lower
estimate of \eqref{3.7}.
\end{proof}
\begin{proof}[Proof of Theorem 3.2] \ {\tt Sufficiency}.

\s
In \eqref{3.2} set $w(x)=\mu(x)^p,$\ $v(x)=\theta(x)^{-p}.$
Then $H^{(p)}=H_p<\iy$ and \eqref{3.2} holds for every $F\in A C(R)^{(+)}$
by
\thmref{theorem3.1}. Denote
\begin{equation}\label{3.19a}
F(x)=\int_x^\iy\theta(s)f(s)ds,\quad \forall f\in L_p(\mathbb R),
\; x\in \mathbb R.
\end{equation}
Since $H_p<\iy,$ then $\theta\in L_{p'}(x,\iy)$ for every $x\in \mathbb R;$ by
H\"older's inequality the integral
\eqref{3.19a} converges for every $x\in \mathbb R,$\ $F(x)$ is absolutely continuous,
$F\in AC(R)^{(+)}$ (see \eqref{3.1})
and, in addition,
\begin{equation}\label{3.20a}
-\frac{1}{\theta(x)}F'(x)=f(x),\quad f\in L_p(\mathbb R).
\end{equation}
Thus \thmref{theorem3.1} together with \eqref{3.20a} reduce to the inequalities
\begin{align*}\|Kf\|_p&=\|\mu F\|_p \\
&\le (p)^{\frac{1}{p}}(p')^{\frac{1}{p'}}H_p
\Big\|\frac{1}{\theta}F' \Big\|_p\\
&=(p)^{\frac{1}{p}}(p')^{\frac{1}{p'}}H_p\|f\|_p.  \tag*{\QED}
\end{align*}
\end{proof}

\begin{remark}\label{rem3.1}
We have formulated the proof of \thmref{theorem3.2} especially for this
paper.
\end{remark}

\begin{proof}[Proof of Lemma 3.4]
The first assertion of the lemma is obvious.
Let us check \eqref{3.15} for $h\in[0,\hat d(x)]$ (the case $h\in[-\hat
d(x),0]$ can be treated
in a similar way).
From condition 4) it follows that
\begin{align}\label{3.25}
\begin{split}
\Phi(x+h,h+\hat d(x))&\ge\Phi(x,\hat d(x))=1 \\
\Phi(x+h,\hat d(x)-h)&\le \Phi(x,\hat d(x))=1.
\end{split}
\end{align}
From \eqref{3.25} and the definition of $\hat d(x)$, we obtain, respectively,
$\hat
d(x+h)  \le$ 
$\hat d(x)+h,$ \ $\hat d(x+h) \ge \hat d(x)-h,$ which is equivalent to
\eqref{3.15}.
Let us check \eqref{3.16}.
Denote $t-x=h.$
Then $t=x+h,$\ $|h|\le \ve\hat d(x)\le \hat d(x)$ since $\ve\in [0,1].$
Therefore \eqref{3.15} implies
$|\hat d(t)-\hat d(x)|=\le|h|\le\ve\hat d(x)$, and hence
\begin{align*}
\left|\fc{\hat d(t)}{\hat d(x)}-1\right|\le\ve \quad\text{ if }\;
|t-x|\le\ve\hat d(x),\; \ve\in[0,1].
\end{align*}
The last inequality is equivalent to \eqref{3.16}.
\QED \end{proof}

\begin{proof}[Proof of Lemma 3.7]
Let us verify that an $R(x,\varkappa(\cdot))$-covering exists for $[x,\iy).$
(The case $(-\iy,x]$ can be treated in a similar way.)
{}From \eqref{3.17} it follows that
\begin{equation}\label{3.26}
\lim_{t\to\iy} (t-\varkappa(t))=\iy.
\end{equation}
Set $\vp(t)=t-\varkappa(t)-x.$
Then $\vp(x)=-\varkappa(x)<0,$ and by \eqref{3.26} there is an $a>x$ such that
$\vp(a)>0.$
Since $\varkappa(t)$ is continuous, so is $\vp(t),$ and $\vp(x)<0,$\
$\vp(a)>0.$
Hence there is an $x_1\in (x,a)$ such that $\vp(x_1)=0,$ i.e.,
$x=x_1-\varkappa(x_1).$
Set $\Dl_1^\pm=x_1\pm \varkappa(x_1),$ and the segment
$\Dl_1=[\Dl_1^-,\Dl^+_1]$ is constructed.
The segments $\Dl_n,$\ $n\ge 2$ with the property $\Dl_n^+=\Dl_{n+1}^-$ are
constructed in a
similar way.
Let us verify that $\cupl_{n\ge 1}\Dl_n=[x,\iy).$
If this is not the case, then there is a $z\in(x,\iy)$ such that $\Dl_n^+<z$ for
all $n\ge 1.$
Since the sequence $\{x_n\}_{n=1}^\iy$ is increasing (by construction) and
bounded
$(x_n<\Dl_n^+<z,\ n\ge 1),$ it has a limit $x_0\le z.$
Moreover,
$\iy>z-x\ge 2\sum_{n=1}^\iy \varkappa(x_n)$, and hence
$$
\varkappa(x_n)\to 0\quad\text{for}\quad n\to\iy.
$$
Then $\varkappa(x_0)=0,$ a contradiction.
The lemma is proved.
\QED \end{proof}

\begin{remark}\label{rem3.4}
\lemref{lemma3.2} is proved by M.~O. Otelbaev's method (see \cite[Ch.1, \S4]{O},
\cite[Ch.III,~\S1]{MO}).
\end{remark}

\begin{proof}[Proof of Lemma 3.8]
According to \lemref{lemma3.2}, it is enough to verify that
\begin{equation}\label{3.27}
\underset{x\to\iy}\varliminf(x-d(x))=\iy,\quad
\left(\underset{x\to-\iy}\varlimsup (x+d(x))=-\iy\right).
\end{equation}
Equalities \eqref{3.27} are checked in the same way.
Let us obtain, say, the first one.
Assume that 
$\varliminf_{x\to\iy}(x-d(x))=c<\iy.$
Then there exists a sequence $\{x_n\}_{n=1}^\iy$ such that $x_n\to\iy$ as
$n\to\iy$ and
$x_n-d(x_n)\le c+1$ for $n\ge 1.$
Together with \eqref{2.2} this implies
$$2=\int_{x_n-d(x_n)}^{x_n+d(x_n)}\fc{q(t)}{r(t)}dt\ge
\int_{c+1}^{x_n+d(x_n)}\fc{q(t)}{r(t)}dt\to\iy\quad\text{as}\quad n\to
\iy,$$
a contradiction.
\QED \end{proof}

Below we need the following lemma.
\begin{lemma}\label{lemma3.4}
Let $S_1=\iy$ (see \eqref{2.1}).
The inequality $\eta\ge d(x)$ (resp., $0\le \eta\le d(x))$ holds if and only if
\begin{equation}\label{3.28}
\int_{x-\eta}^{x+\eta}\fc{q(t)}{r(t)}dt\ge 2\quad
\left(\text{resp.,}\
\int_{x-\eta}^{x+\eta}\fc{q(t)}{r(t)}dt\le 2\right).\end{equation}
\end{lemma}

\begin{proof}[Proof of Lemma 3.4] \ {\tt Necessity}.

\s
Let $\eta \ge d(x).$
Then $[x-d(x),x+d(x)]\subseteq [x-\eta,x+\eta]$, and by \eqref{2.2} we have
$$
\int_{x-\eta}^{x+\eta}\fc{q(t)}{r(t)}dt
\ge\int_{x-d(x)}^{x+d(x)}\fc{q(t)}{r(t)}dt=2.
$$
\end{proof}

\begin{proof}[Proof of Lemma 3.4] \ {\tt Sufficiency}.

\s
Suppose that \eqref{3.28} holds.
Assume the contrary: $\eta<d(x).$
Then \linebreak $[x-\eta,x+\eta]\subset [x-d(x),x+d(x)]$,
and by \eqref{2.2} we get
$$
2\le\int_{x-\eta}^{x+\eta}\fc{q(t)}{r(t)}dt
<\int_{x-d(x)}^{x+d(x)}\fc{q(t)}{r(t)}dt=2,$$
a contradiction.
Hence $\eta\ge d(x).$
\QED \end{proof}

\begin{proof}[Proof of Lemma 2.3] \ {\tt Necessity}.

\s
Let $S_1=\iy,$\ $d_0<\iy$ (see \eqref{2.1}, \eqref{2.10}).
Set $a=d_0.$
Then \eqref{2.2} implies $B(a)\bigm|_{a=d_0}\ge 2$ (see \eqref{2.11}):
$$B(d_0)=\inf_{x\in \mathbb R}\int_{x-d_0}^{x+d_0}\frac{q(t)}{r(t)}
dt\ge\inf_{x\in
R}\int_{x-d(x)}^{x+d(x)}\fc{q(t)}{r(t)}dt=2.$$
  \end{proof}

\begin{proof}[Proof of Lemma 2.3] \ {\tt Sufficiency}.

\s
Let $B(a)>0$ for some $a\in (0,\iy),$ and let $k$ be the smallest natural
number such that $(2k+1)B(a)\ge 2.$
Let $x$ be an arbitrary point from $\mathbb R,$\ $x_n=x+2na,$\
$n=\pm1,\pm2,\dots\ \pm
k .$
Then
\begin{align*}
\int_{x-(2k+1)a}^{x+(2k+1)a}\fc{q(t)}{r(t)}dt
&=\int_{x-a}^{x+a}\fc{q(t)}{r(t)}dt+\sum_{n=1}^k
\int_{x_n-a}^{x_n+a}\fc{q(t)}{r(t)}dt
+\sum_{n=-k}^{-1}\int_{x_n-a}^{x_n+a}\fc{q(t)}{r(t)}dt \\[0.2cm]
&\ge (2k+1)B(a)\ge 2.
\end{align*}
\lemref{lemma3.4} and the above inequality imply $d(x)\le (2k+1)a,$
i.e., 
$d_0=$\linebreak $\sup_{x\in \mathbb R}d(x)<\iy.$
We conclude that $S_1=\iy.$
\QED \end{proof}

\begin{proof}[Proof of Theorem 2.11]
Let $\Phi(x,d)$ be the function from \exampref{examp3.1}.
Then for $x\in \mathbb R$ and $\eta\ge0$, we get
\begin{align}\label{3.29}
\begin{split}
2\Phi(x,\eta)&=2\fc{q_1(x)}{r(x)}\eta
+\int_{x-\eta}^{x+\eta}\fc{q_2(t)}{r(t)}dt \\[0.2cm]
&\phantom{= }+\int_0^\eta\left[\fc{q_1(x+t)}{r(x+t)}-2\fc{q_1(x)}{r(x)}
+\fc {q_1(x-t)}{r(x-t)}\right]dt.
\end{split}
\end{align}
Set
$$
\eta_1(x)=\big(1+\varkappa_1(x)+\varkappa_2(x)\big)\fc{r(x)}{q_1(x)}.
$$
Then $\varkappa_1(x)+\varkappa_2(x)\le 1$ for $|x|\gg 1,$ and, together with
\eqref{3.28}, the equation
\eqref{3.29} implies
\begin{align*}
2 \Phi(x,\eta_1(x)) &\ge
2(1+\varkappa_1(x)+\varkappa_2(x))
-\left|\int_{x-\eta_1(x)}^{x+\eta_1(x)}\fc{q_2(t)}{r(t)}dt\right| \\
&\phantom{\ge }
-\left|\int_0^{\eta_1(x)}\left[\fc{q_1(x+t)}{r(x+t)}-2\fc{q_1(x)}{r(x)}
+\fc{q_1(x-t)}{r(x-t)}\right]dt\right|\\[0.2cm]
& \ge 2(1+\varkappa_1(x)+\varkappa_2(x))
-\sup_{|z|\le 2\fc{r(x)}{q_1(x)}}
\left|\int_{x-z}^{x+z}\fc{q_2(t)}{r(t)}dt\right| \\
&-\sup_{|x|\le 2\fc{r(x)}{q_1(x)}}
\left|\int_0^z\left[\fc{q_1(x+t)}{r(x+t)}-2\fc{q_1(x)}{r(x)}
+\fc{q_1(x-t)}{r(x-t)}\right]dt\right|\\[0.2cm]
&=2+\varkappa_1(x)+\varkappa_2(x) \ge 2.
\end{align*}
Together with \lemref{lemma3.4}, this implies $d(x)\le\eta_1(x)$ for $|x|\gg1.$
Similarly, we set
$
\eta_2=\big(1-\varkappa_1(x)-\varkappa_2(x)\big)\fc{r(x)}{q_1(x)}.
$
Then $\eta_2(x)>0$ for $|x|\gg 1,$ and \eqref{3.29} implies
$2 \Phi(x,\eta_2(x))\le 2$
Together with \lemref{lemma3.4}, this implies $d(x)\ge\eta_2(x)$ for $|x|\gg1.$
We thus get \eqref{2.24}.
Estimates \eqref{2.25} are easily derived from \eqref{2.24} taking into account
that the functions
$r(\cdot),$ $q_1(\cdot)$ and $d(\cdot)$ are continuous and positive on
$\mathbb R.$
\QED \end{proof}


\section{Proof of the main results in the case $p\in (1,\iy)$}\label{PMR}

In this section we prove Theorems \ref{Theorem2.1} and \ref{Theorem2.2} and
their corollaries.
We need the following assertion.

\begin{lemma}\label{lemma4.1}
Let $f \in L_p(\mathbb R),$\ $p\in[1,\iy].$
If the problem \eqref{1.1} -- \eqref{1.2} is solvable (not necessarily
correctly solvable), then
its solution $y(\cdot)$ is unique and can be represented by formula \eqref{1.12}
on $\mathbb R.$
\end{lemma}

\begin{proof}
Let $y(\cdot)$ be a solution of \eqref{1.1} -- \eqref{1.2}.
Then
\begin{equation}\label{4.1}
-y'(\xi)+\fc{q(\xi)}{r(\xi)}y(\xi)=\fc{f(\xi)}{r(\xi)},\quad \xi\in \mathbb R.
\end{equation}
Fix $x\in \mathbb R$ and multiply \eqref{4.1} by
$\left[-\exp\left(-\int_x^\xi\fc{q(s)}{r(s)}ds\right)\right]$.
We get
\begin{equation}\label{4.2}
\fc{d}{d\xi}\left[y(\xi)\exp\left(-\int_x^\xi\fc{q(s)}{r(s)}ds\right)\right]\\
=-\fc{1}{r(\xi)}\exp\left(-\int_x^\xi
\fc{q(s)}{r(s)}ds\right)f(\xi).
\end{equation}
Let $a\in(0,\iy).$
We integrate \eqref{4.2} along the segment $[x,x+a]$ and get
\begin{align}\label{4.3}
\begin{split}
y(x+a)\exp&\left(-\int_x^{x+a}\fc{q(s)}{r(s)}ds\right)-y(x)\\
&=-\int_x^{x+a}\fc{1}{r(\xi)}
\exp\left(-\int_x^\xi\fc{q(s)}{r(s)}ds\right)f(\xi)d\xi.
\end{split}
\end{align}
In \eqref{4.3} we take the limit as $a\to\iy.$
From \eqref{1.2} it follows that the limit in the left-hand side of
\eqref{4.3}
exists. This proves \eqref{1.12}.
\QED \end{proof}

\begin{proof}[Proof of Theorem 2.1] \ {\tt Necessity}.

\s
Suppose that \eqref{1.1} -- \eqref{1.2} is correctly solvable in $L_p(\mathbb R).$
By \lemref{lemma4.1}, its solution $y(\cdot)$ can be written in the form
$y(x)=(Gf)(x), x\in \mathbb R$ (see \eqref{1.12}).
Then
\begin{align}\label{4.4}
\begin{split}
y(x)&=\exp\left(\int_0^x
\fc{q(s)}{r(s)}ds\right)\int_x^\iy
\fc{1}{r(t)}\exp\left(-\int_0^t\fc{q(s)}{r(s)}ds\right)f(t)dt\\
&=\mu(x)\int_x^\iy\theta(t)f(t)dt \\[0.3cm]
&=(Kf)(x)
\end{split}
\end{align}
for all $f\in L_p(\mathbb R)$ (see \eqref{1.12} and \eqref{3.5}).
In \eqref{4.4} we have
\begin{align}\label{4.5}
\begin{split}
\mu(x)&=\exp\left(\int_0^x\fc{q(s)}{r(s)}ds\right)\\
\theta(x)&=\fc{1}{r(x)}\exp\left(-\int_0^x\fc{q(s)}{r(s)}ds\right).
\end{split}
\end{align}

Since \eqref{1.1} -- \eqref{1.2}  is correctly solvable, the operator $G\equiv
K: L_p(\mathbb R)\to L_p(\mathbb R)$
in the case \eqref{4.5} is bounded, and therefore $H_p<\iy$ because of
\thmref{theorem3.2}.
According to \eqref{4.5}, we have the following convenient representation for
$H_p(x)$ (see (3.6)):
\begin{align}
\begin{split}
H_p(x)&=\left[\int_{-\iy}^x
\exp\left(p\int_0^t\fc{q(s)}{r(s)}ds\right)dt\right]^{\frac{1}{p}} \\
&\quad \cdot \left [\int_x^\iy\fc{1}{r(t)^{p'}}
\exp\left(-p'\int_0^t\fc{q(s)}{r(s)}ds\right)dt\right]^{\frac{1}{p'}}\\[0.1cm]
&=\exp\left(-\int_0^x\fc{q(s)}{r(s)}ds\right)
\left[\int_{-\iy}^x\exp\left(p\int_0^t\fc{q(s)}{r(s)}ds
\right)dt\right]^{\frac{1}{p}}\\
&\quad \cdot \exp\left(\int_0^x\fc{q(s)}{r(s)}ds\right)
\left[\int_x^\iy\fc{1}{r(t)^{p'}}
\exp\left(-p'\int_0^t\fc{q(s)}{r(s)}ds\right)dt\right]^{\frac{1}{p'}}\\[0.1cm]
&=\left[\int_{-\iy}^x
\exp\left(-p\int_t^x\fc{q(s)}{r(s)}ds\right)dt\right]^{\frac{1}{p}} \\
&\quad \cdot\left[\int_x^\iy\fc{1}{r(t)^{p'}}
\exp\left(-p'\int_x^t\fc{q(s)}{r(s)}ds\right)dt\right]^{\frac{1}{p'}}.\label{4.6}
\end{split}
\end{align}
From \eqref{4.6} it follows that $H_p(x)=M_p(x)$ (see \eqref{2.4}); therefore
$M_p~=~H_p~<~\iy.$
We thus obtained \eqref{2.3}.
Let us check \eqref{2.5}.
Assume the contrary: $S_1<\iy$ (see \eqref{2.1}).
Then \eqref{2.3} implies
\begin{align*}
\iy& >M_p\\
&=\sup_{x\in \mathbb R}\left[\int_{-\iy}^x
\exp\left(-p\int_t^x\fc{q(s)}{r(s)}ds\right)dt\right]^{\frac{1}{p}}
\left[\int_x^\iy \fc{1}{r(t)^{p'}}
\exp\left(p'\int_x^t\fc{q(s)}{r(s)}ds\right)dt\right]^{\frac{1}{p'}}\\
&\ge \exp(-S_1)\sup_{x\in \mathbb R}
\left[\int_{-\iy}^x 1\,dx\right]^{\frac{1}{p}}
\left[\int_x^\iy\fc{1}{r(t)^{p'}}
\exp\left(-p'\int_x^t\fc{q(s)}{r(s)}ds\right)dt\right]^{\frac{1}{p'}}
=\iy,
\end{align*}
a contradiction.
Hence $S_1=\iy.$
Let us now verify that \eqref{2.6} holds.
We prove \eqref{2.6} in two separate  cases: $T_{p'}=\iy$ and $T_{p'}<\iy,$
where
\begin{equation}\label{4.7}
T_{p'}=\int_0^\iy\fc{dt}{r(t)^{p'}}.
\end{equation}
Let $T_{p'}=\iy.$
Then we get analogously $S_2=\iy$ (see \eqref{3.18}).
Let us now turn to \eqref{2.6} and prove it ad absurdum.
Let $A_{p'}=\iy$ (see \eqref{2.6}).
Let us show that there exists $F(t)\in L_p(\mathbb R)$ such that
$(GF)(t)\not\rightarrow 0$ as $|t|\to \iy,$
i.e., \eqref{1.2} does not hold for all $f(t)\in L_p(\mathbb R).$
This contradicts the correct solvability of \eqref{1.1} -- \eqref{1.2} in
$L_p(\mathbb R).$
Consider the function $A_{p'}(x),$\ $x\in \mathbb R$ (see \eqref{2.7}).
This function is continuous for $x\in \mathbb R$ (see \lemref{lemma3.1} and
\exampref{examp3.1}).
Let $\beta $ be a positive number which will be chosen later, and let
$a_{p'}=\inf\limits_{x\in \mathbb R}
A_{p'}(x).$
Since $A_{p'}(x)$ is continuous for $x\in \mathbb R,$ we have $0\le a_{p'}<\iy,$ and
for every integer
 $k\ge
a_{p'}+1$ there is a point $x_k$ such that
\begin{equation}\label{4.8}
k^\beta\le A_{p'}(x_k)\le (k+1)^\beta, \quad k\ge a_{p'}+1.
\end{equation}

Since $A_{p'}(x)$ is continuous, $|x_k|\to\iy$ as $k\to\iy.$
Since $S_1=S_2=\iy$ (see \eqref{2.1}, \eqref{3.18}), from \eqref{3.27} it
follows that one can choose
a subsequence $\{x_{k_n}\}_{n=1}^\iy$ such that $|x_{k_n}|\to\iy$ as $n\to\iy,$
and the segments
$$
\Dl_{k_n}=[\Dl_{k_n}^-,\Dl_{k_n}^+]=[x_{k_n}-d(x_{k_n}),x_{k_n}+d(x_{k_n})],
\quad n=1,2,\ldots,
$$
are disjoint.
Thus for a given $\beta>0$ there is a sequence $\{x_{k_n}\}_{n=1}^\iy$ such
that
\begin{alignat}{2}
\begin{aligned}\label{4.9}
&k_n^\beta\le \int_{\Dl_{k_n}}\fc{dt}{r(t)^{p'}}\le (k_n+1)^\beta
&\quad &\text{for } 1+a_{p'}\le k_n<k_{n+1}\\[0.2cm]
&\D_{k_n}\cap \Dl_{k_m}=\emptyset &\quad &\text{for }  n\ne m.
\end{aligned}
\end{alignat}
Let $\a$ be another positive number (also to be chosen later) and
\begin{align}\label{4.10}
f_{k_n}(t)&=\begin{cases}
 \frac{1}{(1+k_n)^\a}\fc{1}{r(t)^{p'-1}},\ & t\in
\Dl_{k_n}\\[0.1cm]
0,\ & t\notin \Dl_{k_n}
\end{cases}, \qquad n=1,2,\dots \\[0.1cm]
\label{4.11}
F(t)&=\sum_{n=1}^\iy f_{k_n}(t),\hspace{1.47cm} t\in \mathbb R.
\end{align}
Let us verify that $F \in L_p(\mathbb R)$ for $\a>\fc{1+\beta}{p}.$
Indeed, by \eqref{4.9} -- \eqref{4.10} we have
\begin{align}\label{4.12}
\begin{split}
\|F\|_p^p&=\sum_{n=1}^\iy
\int_{\Dl_{k_n}}\fc{1}{(1+k_n)^{\a p}}\fc{dt}{r(t)^{p(p'-1)}}\\
&=\sum_{n=1}^\iy \fc{1}{(1+k_n)^{\a p}}
\int_{\Dl_{k_n}}\fc{dt}{r(t)^{p'}} \\
&\le \sum_{n=1}^\iy \fc{(1+k_n)^\beta}{(1+k_n)^{p\a}} \\
&\le\sum_{n=1}^\iy\fc{1}{n^{p\a-\beta}} <\iy.
\end{split}
\end{align}
Let us estimate $(GF)(t)$ from below for $t=\Dl_{k_n}^-=x_{k_n}-d(x_{k_n}).$
By  \eqref{1.12},  \eqref{4.11},  \eqref{2.2},  \eqref{4.10} and \eqref{4.9},
we get
\begin{align} \label{4.13}
\begin{split}
(GF)(\Dl_{k_n}^-) &\ge
\int_{\Dl_{k_n}}\fc{1}{r(t)}f_{k_n}(t)
\exp\left(-\int_{\Dl_{k_n}^-}^t\fc{q(s)}{r(s)} ds\right)dt \\
&\ge \exp\left(-\int_{\Dl_{k_n}}\fc{q(s)}{r(s)}ds\right)
\int_{\Dl_{k_n}}\fc{f_{k_n}(t)}{r(t)}dt \\
&=\exp(-2)\fc{1}{(1+k_n)^\a}\int_{\Dl_{k_n}}\fc{dt}{r(t)^{p'}} \\
&\ge \exp(-2)\fc{k_n^\beta}{(1+k_n)^\a} \\[0.2cm]
&\ge c^{-1}(k_{n}+1)^{\beta-\a}.
\end{split}
\end{align}
For $\beta\ge \a$ we obtain from \eqref{4.13} that $(GF)(t)\not\rightarrow 0$
as $|t|\to\iy.$
Thus we shall reach the goal as soon as \eqref{4.12} and \eqref{4.13} will hold
together.
This will be true since there are $\a,\beta\in(0,\iy)$ such that $\beta\ge
\a>\fc{1+\beta}{p}.$
But it is not hard to satisfy these inequalities: one can take, say,
$\a=\beta,$\
$\beta>\fc{1}{p-1}$, a contradiction.
Hence $A_{p'}<\iy.$
Let now $T_{p'}<\iy$ (see \eqref{4.7}).
First, let us show that (regardless of the value of $T_{p'})$ there exists
$c_0\in \mathbb R$ such that the
inequality
\begin{equation}\label{4.14}
x-d(x)\ge c_0
\end{equation}
holds for all $x\ge 0.$
Assume the contrary.
Then there is a sequence $\{x_n\}_{n=1}^\iy$ such that $x_n-d(x_n)\le -n,$\
$x_n\ge 0,$\ $n=1,2,\dots\ $. Together with \eqref{2.2} and \eqref{2.5},
this implies
$$
2=\int_{x_n-d(x_n)}^{x_n+d(x_n)}\fc{q(t)}{r(t)}dt\ge
\int_{-n}^0\fc{q(t)}{r(t)}dt\to\iy \quad \text{as}\; n\to\iy,
$$
a contradiction.
Then, since $T_{p'}<\iy,$ we obtain that the function $A_{p'}(x)$ is absolutely
bounded for $x\ge0$.
Indeed, from \eqref{4.14} we get for $x \ge 0$
\begin{equation}\label{4.15} A_{p'}(x)
=\int_{x-d(x)}^{x+d(x)}\fc{dt}{r(t)^{p'}}
\le \left|\int_{c_0}^0\fc{dt}{r(t)^{p'}}\right|+\int_0^\iy\fc{dt}{r(t)^{p'}}
 =c+T_{p'}<\iy .
 \end{equation}
Taking into account \eqref{4.15}, let us assume now that
\begin{equation}\label{4.16}
\tilde A_{p'}=\sup_{x\le
0}\int_{x-d(x)}^{x+d(x)}\fc{dt}{r(t)^{p'}}=\iy.
\end{equation}
We then can repeat all the arguments from the preceding case $(T_{p'}=\iy).$
The only difference is that the initial sequence $\{x_k\}_{k=1}^\iy$ (see
\eqref{4.8}) is known to
satisfy the property $x_k\to -\iy$ as $k\to \iy.$
Taking into account this remark, we reduce this case $(T_{p'}<\iy)$ to the
preceding one
$(T_{p'}=\iy).$
We thus proved the necessity of all the conditions of the theorem.
\end{proof}

\begin{proof}[Proof of Theorem 2.1] \ {\tt Sufficiency}.

\s
Consider the function $y(x)=(Gf)(x),$\ $f(x)\in L_p(\mathbb R)$
(see \eqref{1.12}).
Since\linebreak $M_p < \iy,$ we have $M_p(x)<\iy$ for every
$x\in \mathbb R$ (see \eqref{2.3} -- \eqref{2.4}), and
therefore
\begin{equation}\label{4.17}
\int_x^\iy\fc{1}{r(t)^{p'}}
\exp\left(-p'\int_x^t\fc{q(s)}{r(s)}ds\right)dt<\iy \quad\text{for
every}\; x\in \mathbb R.
\end{equation}
Since $f \in L_p(\mathbb R),$ from \eqref{4.17} and H\"older's inequality , it
follows that the
integral
$(Gf)(x)$ converges for $x\in \mathbb R.$
We now differentiate \eqref{1.12} and conclude that $y(x)=(Gf)(x)$ satisfies
\eqref{1.1} almost
everywhere on $\mathbb R.$
{}From \eqref{4.6}, the equality $M_p(x)\equiv H_p(x),$\ $x\in \mathbb R,$ and
\thmref{theorem3.2} we get
\eqref{1.5}.
It remains to prove that equalities \eqref{1.2} hold for $y=(Gf)(x),$\
$x\in \mathbb R$ regardless of $f\in L_p(\mathbb R).$
{}From \eqref{1.12} and H\"older's inequality we get
\begin{align} \label{4.18}
\begin{split}
0&\le |(Gf)(x)| \\
&\le \left[\int_x^\iy\fc{1}{r(t)^{p'}}
\exp\left(-\int_x^t\fc{q(s)}{r(s)}ds\right)dt\right]^{\frac{1}{p'}}\\
&\phantom{\le } \cdot\left[\int_x ^\iy |f(t)|^p
\exp\left(-\int_x^t\fc{q(s)}{r(s)}ds\right)dt\right]^{\frac{1}{p}},\quad
x\in \mathbb R.
\end{split}
\end{align}
 \end{proof}
Below we need the following lemma.

\begin{lemma}\label{lemma4.2}
Suppose that conditions \eqref{2.5} -- \eqref{2.6} hold.
Then $J_\nu<\iy$ for any $\nu>0$, where
\begin{align}\label{4.19}
\begin{split}
J_\nu &\doe\sup\limits_{x\in \mathbb R} J_\nu(x) \\
J_\nu(x)&=\int_x^\iy\fc{1}{r(t)^{p'}}
\exp\left(-\nu\int_x^t\fc{q(s)}{r(s)}ds\right)dt,\quad x\in \mathbb R.
\end{split}
\end{align}
\end{lemma}

\renewcommand{\qedsymbol}{\openbox}
\begin{proof}
Consider two cases: $T_{p'}=\iy$ and $T_{p'}<\iy$ (see \eqref{4.7}).
If $T_{p'}=\iy$, then $S_2=\iy$ (see above), and by \lemref{lemma3.3} there is
an
$R(x,d(\cdot))$-covering of $[x,\iy).$
Let $\Dl_n,$\ $n\ge 1$ be the segments of the $R(x,d(\cdot))$-covering of
$[x,\iy).$
Then, if $t\in \Dl_n,$\ $n\ge 1,$ we have
\begin{equation}\label{4.20}
\int_x^t\fc{q(s)}{r(s)} ds\ge 2(n-1).
\end{equation}
Indeed, for $n=1$ estimate \eqref{4.20} is obvious, and for $n\ge 2$ we have
\begin{equation}\label{4.21}
\int_x^t\fc{q(s)}{r(s)}ds=
\sum_{k=1}^{n-1}\int_{\Dl_k}\fc{q(s)}{r(s)}ds
+\int_{\Dl_n^-}^t\fc{q(s)}{r(s)}ds\ge
\sum_{k=1}^{n-1}2=2(n-1).
\end{equation}
By the properties of a $R(x,d(\cdot))$-covering of $[x,\iy)$ and by
\eqref{4.21} and \eqref{2.6}, we
get
\begin{align*}
J_\nu(x)&=\sum_{n=1}^\iy\int_{\Dl_n}\fc{1}{r(t)^{p'}}
\exp\left(-\nu\int_x^t\fc{q(s)}{r(s)}ds\right)dt \\
&\le \sum_{n=1}^\iy\exp(-2(n-1)\nu)\int_{\Dl_n}\fc{dt}{r(t)^{p'}}\\
&\le A_{p'}\sum_{n=1}^\iy\exp(-2(n-1)\nu)\\[0.1cm]
&=c_\nu A_{p'}
\end{align*}
which implies $J_\nu\le c_\nu A_{p'}<\iy$.
Let now $T_{p'}<\iy.$
Since $S_1=\iy$ (see \eqref{2.5}), by \lemref{lemma3.3} there is an
$R({0},d(\cdot))$-covering of
$(-\iy,0].$
Let $\{\Dl_n\}_{n=-\iy}^{-1}$ be the segments forming this covering.
When estimating $J_\nu(x),$ first consider the case $x<0.$
Then $x\in \Dl_{n_0},$\ $n_0\le -1$ and \eqref{2.2} imply
\begin{align} \label{4.22}
\begin{split}
\int_x^0 \fc{1}{r(t)^{p'}}&
\exp\left(-\nu\int_x^t\fc{q(s)}{r(s)}ds\right)dt\\
&=\exp\left(\nu\int_{\Dl_{n_0}^-}^x\fc{q(s)}{r(s)}ds\right)
\int_x^0\fc{1}{r(t)^{p'}}
\exp\left(-\nu\int_{\Dl_{n_0}^-}^t\fc{q(s)}{r(s)}ds\right)dt\\
&\le \exp(2\nu)\int_{\Dl_{n_0}^-}^0\fc{1}{r(t)^{p'}}
\exp\left(-\nu\int_{\D_{n_0}^-}^t\fc{q(s)}{r(s)}ds\right)dt \\
&=\exp(2\nu)\sum_{k=n_0}^{-1}\int_{\Dl_k}\fc{1}{r(t)^{p'}}
\exp\left(-\nu\int_{\Dl_{n_0}^-}^t\fc{q(s)}{r(s)}ds\right)dt\\
&\le \exp(2\nu)\sum_{k=n_0}^{-1}\int_{\Dl_k}\fc{1}{r(t)^{p'}}
\exp\left(-\nu\int_{\Dl_{n_0}^-}^{\Dl_k^-}\fc{q(s)}{r(s)}ds\right)dt\\
&\le \exp(2\nu)A_{p'}\sum_{k=n_0}^{-1}\exp(-2\nu|n_0-k|)\\
&\le\exp(2\nu)A_{p'}\sum_{m=0}^\iy\exp(-2\nu m)\\[0.1cm]
&=c(\nu)A_{p'}.
\end{split}
\end{align}
Using \eqref{4.22}, we now get for $x\le 0$
\begin{align*}
J_\nu(x)&=\int_x^\iy\fc{1}{r(t)^{p'}}
\exp\left(-\nu\int_x^t\fc{q(s)}{r(s)}ds\right)dt \\
&=\int_x^0\fc{1}{r(t)^{p'}}\exp\left(-\nu\int_x^t\fc{q(s)}{r(s)}ds\right)dt
+\int_0^\iy\fc{1}{r(t)^{p'}}\exp\left(-\nu\int_x^t\fc{q(s)}{r(s)}ds\right)dt \\
&\le c(\nu)A_{p'}+\int_0^\iy\fc{dt} {r(t)^{p'}} \\[0.2cm]
&=c(\nu)A_{p'}+T_{p'}.
\end{align*}
For $x\ge0,$ $J_v(x)\le T_{p'}$ is obvious.
Thus, $ J_\nu\le c(\nu)A_{p'}+T_{p'}<\iy.$ \QED
\end{proof}
Now it is easy to see that $y(x)=(Gf)(x)\to0$ as $x\to\iy.$
Indeed from \eqref{4.18} and \lemref{4.2} it follows that
$$
0\le |(Gf)(x)|\le J_1^{\frac{1}{p'}}
\left[\int_x^\iy|f(t)|^pdt\right]^{\frac{1}{p}}\to 0
\quad\text{as}\; x\to\iy.
$$
To check \eqref{1.2} for $x\to-\iy,$ we use the following two lemmas.
\begin{lemma}\label{lemma4.3}
Let $S_1=\iy$ (see \eqref{2.1}).
Then for every $\eta\in(0,\iy)$ there is an $x_0(\eta)\le0$ such that for every
$x\le x_0(\eta)$ the
equation in $d\ge0$
\begin{equation}\label{4.23}
\Phi(d)\doe\int_x^{x+d}\fc{q(s)}{r(s)}ds=\eta
\end{equation}
has at least one solution $\hat d(x,\eta)$, and $x+\hat d(x,\eta)\le 0.$
\end{lemma}

\begin{proof}
Clearly, for $d\in [0,\iy)$ the function $\Phi(d)$ is continuous, non-negative,
and $\Phi(0)=0.$
Since $S_1=\iy,$ there is an $x_0(\eta)$ such that
\begin{equation}\label{4.24}
\int_{x_0(\eta)}^0\fc{q(s)}{r(s)}ds\ge 2\eta.
\end{equation}
Set $\mu(x)=-x$ for every $x\le x_0(\eta).$
Then
$$
\int_x^{x+\mu(x)}\fc{q(s)}{r(s)}ds\ge\int_{x_0(\eta)}^0\fc{q(s)}{r(s)}ds
\ge 2\eta>\eta.
$$
Thus $\Phi(0)=0,$\ $\Phi(\mu(x))>\eta$ and therefore, since $\Phi(d)$ is
continuous, in the segment
$[0,\mu(x)]$ there is at least one root $d=\hat d(x,\eta)$ of equation
\eqref{4.23}, and since $\hat d(x,\eta)\le \mu(x),$
we have $x+\hat d(x,\eta)\le 0.$
\QED
\end{proof}

Let $d(x,\eta),$\ $x\le x_0(\eta)$ be the smallest root of equation
\eqref{4.23}:
\begin{equation}\label{4.25}
d(x,\eta)=\inf_{d>0}\bigg\{d: \int_x^{x+d}\fc{q(s)}{r(s)}ds=\eta,\quad
x\le x_0(\eta)\bigg\}.
\end{equation}

\begin{lemma}\label{lemma4.4}
Let $S_1=\iy$ (see \eqref{2.1}).
Then
\begin{equation}\label{4.26}
\lim_{x\to-\iy}(x+d(x,\eta))=-\iy,\quad \eta\in(0,\iy).
\end{equation}
\end{lemma}
\begin{proof}
Assume that for some $\eta\in(0,\iy)$ we have
\begin{equation}\label{4.27}
\underset{x\to-\iy}\varlimsup(x+d(x,\eta))=c, \quad c\in \mathbb R.
\end{equation}
{}From \lemref{lemma4.3} and \eqref{4.25}, it follows that $c\le 0.$
{}From \eqref{4.27}, it follows that there is a sequence $\{x_n\}_{n=1}^\iy$
such that $x_n\to-\iy$ as
$n\to\iy$ and $x_n+d(x_n,\eta)\ge c-1$ for all $n\ge1.$
Hence
$$
\eta=\int_{x_n}^{x_n+d(x_n,\eta)}\fc{q(s)}{r(s)}ds\ge
\int_{x_n}^{c-1}\fc{q(s)}{r(s)}ds\to
\iy \quad \text{for}\; n\to \iy,
$$
a contradiction.
Thus $\underset{x\to-\iy}\varlimsup(x+d(x,\eta)=-\iy$ and therefore
$$
-\iy\le\underset{x\to-\iy}\varliminf(x+d(x,\eta))
\le\underset{x\to\iy}\varlimsup(x+d(x,\eta))=-\iy.
$$
The last relations immediately imply \eqref{4.26}.
\QED \end{proof}

Fix now $\eta\in(0,\iy).$
Let $x_0(\eta)$ be the number defined in \lemref{lemma4.3} and $x\le
x_0(\eta).$
{}From \eqref{4.18}, \lemref{lemma4.2} and \eqref{4.25}, it follows that
\begin{align} \label{4.28}
\begin{split}
0\le|(Gf)(x)|^p&\le
J_1^{\frac{p}{p'}}\cdot\int_x^\iy|f(t)|^p
\exp\left(-\int_{x}^t\fc{q(s)}{r(s)}ds\right)dt \\
&\le c\int_x^{x+d(x,\eta)}|f(t)|^pdt
+c\exp\bigg(-\int_x^{x+d(x,\eta)}\fc{q(s)}{r(s)}ds\bigg) \\
&\phantom{\le }\  \cdot \int_{x+d(x,\eta)}^\iy|f(t)|^{p}
\exp\left(-\int_{x+d(x,\eta)}^t\fc{q(s)}{r(s)}ds\right)dt \\
&\le c\int_x^{x+d(x,\eta)}|f(t)|^pdt+c\exp(-\eta)\|f\|_p^p.
\end{split}
\end{align}
Since $f \in L_p(\mathbb R),$ from \eqref{4.26} it follows that
$$
0 \le \int_x^{x+d(x,\eta)}|f(t)|^pdt\le
\int_{-\iy}^{x+d(x,\eta)}
|f(t)|^pdt\to0 \quad \text{as}\; x\to-\iy,
$$
 i.e., the first summand in \eqref{4.28} tends to zero as $x\to-\iy.$
Hence
\begin{equation}\label{4.29}
0\le\underset{x\to-\iy}\varlimsup|(Gf)(x)|^p\le c\exp(-\eta)\|f\|_p^p.
\end{equation}
In \eqref{4.29} the number $\eta\in (0,\iy)$ is arbitrary.
Hence in \eqref{4.29} we can take the limit as $\eta\to\iy.$
We get
$$
\underset{x\to-\iy}\varlimsup |(Gf)(x)|=0\ \Rightarrow\
0\le\underset{x\to-\iy}\varliminf|(Gf)(x)|
\le\underset{x\to-\iy}\varlimsup|(Gf)(x)|=0.
$$
The last relations finish the proof of \thmref{Theorem2.1}. \QED

\begin{proof}[Proof of \corref{Corollary2.1.1}]

\s
Since problem \eqref{1.1}--\eqref{1.2} is correctly solvable in
$L_p(\mathbb R),$ by \thmref{Theorem2.1} we
have $S_1=\iy$ and $A_{p'}<\iy$ (see \eqref{2.5}--\eqref{2.6}).
According to \eqref{1.12}, we have to prove that the integral operator
\eqref{3.5} with
\begin{align}\label{4.30}
\begin{split}
\mu(t)&=\left(\fc{q(t)}{r(t)}\right)^{\frac{1}{p}}
\exp\left(\int_0^t\fc{q(s)}{r(s)}ds\right) \\
\theta(t)&=\fc{1}{r(t)}\exp\left(-\int_0^t\fc{q(s)}{r(s)}ds\right)
\end{split}
\end{align}
is bounded in $L_p(\mathbb R).$
Below, when estimating $H_p(x)$ (see \eqref{3.6}), we use \lemref{lemma4.2}
and get
\begin{align*}
H_p(x) &=\left[\int_{-\iy}^x\fc{q(t)}{r(t)}
\exp\left(p\int_0^t\fc{q(s)}{r(s)}ds\right)dt\right]^{\frac{1}{p}} \\
&\phantom{= }\ \cdot \left[\int_x^\iy\fc{1}{r(t)^{p'}}
\exp\left(-p'\int_0^t\fc{q(s)}{r(s)}ds\right)dt\right]^{\frac{1}{p'}}\\
&=\left[\int_{-\iy}^x\fc{q(t)}{r(t)}
\exp\left(-p\int_t^x\fc{q(s)}{r(s)}ds\right)dt\right]^{\frac{1}{p}} \\
&\phantom{= }\ \cdot \left[\int_x^\iy\fc{1}{r(t)^{p'}}
\exp\left(-p'\int_x^t\fc{q(s)}{r(s)}ds\right)dt\right]^{\frac{1}{p'}}\\
&\le J_{p'}^{\frac{1}{p'}}\left[\fc{1}{p}
\exp\left(\left.-p\int_t^x\fc{q(s)}{r(s)}ds\right)
\right|_{-\iy}^x\right]^{\frac{1}{p}}\\
&=\fc{1}{p^{\frac{1}{p}}}
(J_{p'})^{\frac{1}{p'}},
\end{align*}
hence $H_p<\iy$.
It remains to apply \thmref{theorem3.2}.
\QED \end{proof}

\begin{proof}[Proof of \thmref{Theorem2.2}] \ {\tt Necessity}.

\s
The necessity of the condition $A_{p'}<\iy$ (see \eqref{2.6}) for correct
solvability of problem
\eqref{1.1}--\eqref{1.2} in $L_p(\mathbb R)$ follows from \thmref{Theorem2.1}.
\end{proof}

\begin{proof}[Proof of \thmref{Theorem2.2}] \ {\tt Sufficiency}.

\s
Suppose that \eqref{2.11} holds.
Then $S_1=\iy,$\ $d_0<\iy$, in view of \lemref{lemmama2.1}.
Since $B(a)>0,$ we have $S_2=\iy$ (see \eqref{3.18}).
Hence by \lemref{lemma3.3}, we conclude that for every $x\in \mathbb R$
there exist $R(x,d(\cdot))$-coverings for
$(-\iy,x]$ and $[x,\iy).$
Below we need the following lemma.

\begin{lemma}\label{lemma4.5}
Let $S_1=\iy$ and $d_0<\iy$ (see \eqref{2.1} and \eqref{2.10}).
Then $I_\nu<\iy$ for every $\nu>0$, where
\begin{align}\label{4.31}
\begin{split}
I_\nu &\doe \sup\limits_{x\in \mathbb R}I_\nu(x) \\
I_\nu(x)&=\int_{-\iy}^x\exp\left(-\nu\int_t^x\fc{q(s)}{r(s)}ds\right)dt.
\end{split}
\end{align}
\end{lemma}

\begin{proof}
Let $\{\Dl_n\}_{n=-\iy}^{-1} $ be the set of segments forming an
$R(x,d(\cdot))$-covering of $(-\iy,x].$
So, if $t\in \Dl_n,$\ $n\le -1,$ then
\begin{equation}\label{4.32}
\int_t^x\fc{q(s)}{r(s)}ds\ge 2(|n|-1),\quad n\le -1.
\end{equation}
Indeed, for $n=-1 $ estimate \eqref{4.32} is obvious, and for $n\le -2$ we have
$$
\int_t^x\fc{q(s)}{r(s)}ds
=\int_t^{\Dl_n^+}\fc{q(s)}{r(s)}ds
+\sum_{k=n+1}^{-1}\int_{\Dl_k}\fc{q(s)}{r(s)}ds
\ge \sum_{k=n+1}^{-1}2=2(|n|-1).
$$
The next chain of calculations
\begin{align*}
I_\nu(x)&=\sum^{-1}_{n=-\iy}\int_{\Dl_n}
\exp\left(-\nu\int_t^x\fc{q(s)}{r(s)}ds\right)dt \\
&\le 2\sum_{n=-\iy}^{-1} d(x_n)\exp(-2(|n|-1)\nu)\\
&\le 2d_0\sum_{n=1}^\iy\exp(-2(n-1)\nu) \\[0.2cm]
&=c(\nu)d_0<\iy.
\end{align*}
based on \eqref{4.32} and the properties of the
$R(x,d(\cdot))$-covering of $(-\iy,x]$ finishes the proof.
\QED
\end{proof}
Now, using \eqref{2.4}, \eqref{4.19}, \eqref{4.31} and Lemmas \ref{lemma4.2}
and \ref{lemma4.5}, it is easy to
see that
$$
M_p(x)=(I_p(x))^{\frac{1}{p}}(J_{p'}(x))^{\frac{1}{p'}}
\le(I_p)^{\frac{1}{p}}(J_{p'})^{\frac{1}{p'}}<\iy.
$$
Hence $M_p<\iy$ and therefore $\|G\|_{p\to p}<\iy$ (see \eqref{2.8}).
Consequently, all conditions of \thmref{Theorem2.1} are
satisfied. This means that problem \eqref{1.1}--\eqref{1.2} is correctly
solvable in $L_p(\mathbb R).$
\hfill \QED

\begin{proof}[Proof of \corref{Corollary2.2.1}] \ {\tt Necessity}.

\s
Suppose that \eqref{2.12} holds and problem \eqref{1.1}--\eqref{1.2} is
correctly solvable in $L_p(\mathbb R).$
Then $A_{p'}<\iy$ (see \eqref{2.6}), in view of \thmref{Theorem2.1}.
Together with \eqref{2.12} this implies $d_0<\iy$ (see \eqref{2.10}):
$$
\iy> A_{p'}=\sup_{x\in \mathbb R}\int_{x-d(x)}^{x+d(x)}\fc{dt}{r(t)^{p'}}
\ge \fc{2}{c^{p'}}\sup_{x\in \mathbb R}d(x)=cd_0.
$$
By \thmref{Theorem2.1}, we have $S_1=\iy.$
Then by \lemref{lemmama2.1}, there is $a_0\in(0,\iy)$ such that $B(a_0)>0$ (see
\eqref{2.11}).
Hence for every $x\in \mathbb R,$ from \eqref{2.12} we get
$$
c\int_{x-a_0}^{x+a_0}q(t)dt
\ge\int_{x-a_0}^{x+a_0}\fc{q(t)}{r(t)}dt
\ge B(a_0)>0,
$$
which implies (1.6).
\end{proof}

\begin{proof}[Proof of \corref{Corollary2.2.1}] \ {\tt Sufficiency}.

\s
Suppose that for some $a=a_0$ we have \eqref{1.6}.
Then, from \eqref{2.12} it follows that $B(a_0)>0$ (see \eqref{2.11}):
$$
0<\inf_{x\in \mathbb R}\int_{x-a_0}^{x+a_0} q(t)dt
\le c\inf_{x\in \mathbb R}\int_{x-a_0}^{x+a_0}\fc{q(t)}{r(t)}dt.
$$
Hence by \lemref{lemmama2.1}, we get $d_0<\iy.$
Therefore, by \eqref{2.12} we have $A_{p'}<\iy$ (see \eqref{2.6}):
$$
A_{p'}=\sup_{x\in \mathbb R}\int_{x-d(x)}^{x+d(x)}\fc{dt}{r(t)^{p'}}
\le 2c^{p'}\sup_{x\in \mathbb R}d(x)=cd_0<\iy.
$$
It remains to refer to \thmref{Theorem2.2}. \corref{Corollary2.2.1}
is proved.
\QED
\end{proof}
\end{proof}


\section{Proof of the main results in the case $p=1$}\label{ProofMRp1}

In this section we prove Theorems \ref{Theorem2.3} and \ref{Theorem2.4} and
their corollaries.

\begin{proof}[Proof of \thmref{Theorem2.3}] \ {\tt Necessity}.

\s
Suppose that problem \eqref{1.1}--\eqref{1.2} is correctly solvable in
$L_1(\mathbb R).$
By \linebreak \lemref{lemma4.1}, its solution $y(\cdot)$ is of the form
\eqref{1.12}, and the
equalities
\eqref{4.4}--\eqref{4.5} hold. From \eqref{1.5} it follows that the operator
$G:L_1(\mathbb R)\to L_1(\mathbb R)$ (see \eqref{1.12}) is bounded, and
therefore by \thmref{theorem3.3} we have the relations
\eqref{2.15}--\eqref{2.16}. From
\eqref{2.15}--\eqref{2.16}, it follows that $S_1=\iy$ (see \eqref{2.13}).
Indeed,
assuming $S_1<\iy$ we arrive at a contradiction:
\begin{align*}
\iy>M_1&=\sup_{x\in \mathbb R}\fc{1}{r(x)}\int_{-\iy}^x
\exp\left(-\int_t^x\fc{q(s)}{r(s)}ds\right)dt \\
&\ge\exp(-S_1)\sup_{x\in \mathbb R}\fc{1}{r(x)}\int_{-\iy}^x1 dt=\iy.
\end{align*}

Let us now verify that $r_0>0$ (see \eqref{2.14}).
Assume the contrary: $r_0=0.$
Let $\a$ be a positive number which will be chosen later.
Since $r_0=0$, from \eqref{1.3} it follows that there is a sequence
$\{x_n\}_{n=1}^\iy$ such that
$r(x_n)=n^{-\a},$\ $n=1,2,\dots\ .$
Clearly, $|x_n|\to\iy$ as $n\to\iy$ (otherwise, there is a point $x_0$ \
$(|x_0|<\iy)$ such that
$r(x_0)=0,$ a contradiction (see \eqref{1.3}).
From \eqref{1.3} it also follows that there are numbers $\dl_n>0,$\
$n=1,2,\dots$ such that
\begin{equation}\label{5.1}
\fc{1}{2n^\a}\le r(t)\le \fc{2}{n^\a},\quad t\in[x_n-\dl_n,x_n+\dl_n],\;
n=1,2,\dots\ .
\end{equation}
Let $\om_n=\min\{\dl_n,d(x_n),1\},$\ $n=1,2,\dots\ .$
Making $\dl_n$ smaller (if necessary), we can choose $\omega_n,$\ $n=1,2,\dots$
so that the segments
$\Dl_n=]x_n-\om_n,x_n+\om_n],$\ $n=1,2,\dots$, are disjoint.
Let us introduce the functions
\begin{align}\label{5.2}
f_n(t)&=\begin{cases}\fc{1}{2\om_n}\ \fc{1}{n^\a},\quad &t\in\Dl_n\\[0.1cm]
0,\quad &t\notin\Dl_n
\end{cases} \qquad n=1,2,\dots \\
\label{5.3}
F(t)&=\sum_{n=1}^\iy f_n(t),\hspace{0.7cm} t\in \mathbb R.
\end{align}
Let us verify that $F \in L_1(\mathbb R)$ for $\a>1.$
Indeed,
$$\|F\|_1=\sum_{n=1}^\iy\int_{\Dl_n}|f_n(t)|dt=\sum_{n=1}^\iy
\fc{1}{n^\a}\
\fc{2\om_n}{2\om_n}=\sum_{n=1}^\iy\fc{1}{n^\a}<\iy.$$
We now estimate $(GF)(t)$ from below for $t_n=x_n-\om_n,$\ $n=1,2,\ldots:$
\begin{align*}
(GF)(t_n)&\ge\int_{x_n-\om_n}^{x_n+\om_n}\fc{1}{r(t)}
\exp\left(-\int_{x_n-\om_n}^t\fc{q(s)}{r(s)}ds\right)f_n(t)dt\\
&\ge\fc{1}{2\om_n\cdot n^\a}\int_{x_n-\om_n}^{x_n+\om_n}\fc{1}{r(t)}
\exp\left(-\int_{x_n-\om_n}^{x_n+\om_n}\fc{q(s)}{r(s)}ds\right)dt\\[0.1cm]
&\ge \fc{\exp(-2)}{2}.
\end{align*}
Thus, $(GF)(t)\not\rightarrow 0$ as $|t|\to\iy$,
a contradiction.
Hence $r_0>0,$ and the necessity of the conditions of the theorem is proved.
\end{proof}

\begin{proof}[Proof of \thmref{Theorem2.3}] \ {\tt Sufficiency}.

\s
Consider the function $y(x)=(Gf)(x),$ $x\in \mathbb R,$
$f \in L_1(\mathbb R)$ (see \eqref{1.12}).
Since $r_0>0$ (see \eqref{2.14}), the integral $(Gf)(x)$ converges for all
$x\in \mathbb R,$ and we have
\begin{align}\label{5.4}
\begin{split}
|y(x)|=|(Gf)(x)|&\le\fc{1}{r_0}\int_x^\iy
\exp\left(-\int_x^t\fc{q(s)}{r(s)}ds\right)|f(t)|dt \\
&\le\fc{1}{r_0}\int_x^\iy|f(t)|dt.
\end{split}
\end{align}
When differentiating \eqref{1.12}, we immediately conclude that the
function\linebreak $y(\cdot)=(Gf)(\cdot)$ satisfies
\eqref{1.1} almost everywhere on $R.$
Furthermore, since $M_1<\iy$ (see \eqref{2.15}), by \thmref{theorem3.3} the
operator $G: L_1(\mathbb R)\to
L_1(\mathbb R)$ is bounded, and therefore \eqref{1.5} hods.
Clearly, $y(x)=(Gf)(x)\to0$ as $x\to\iy$ because of \eqref{5.4}.
The proof of \eqref{1.2} for $x\to-\iy$ goes exactly as in \thmref{Theorem2.1};
the only difference is as
follows: instead of the   equality \eqref{4.18} one has to use the first
inequality of~\eqref{5.4}.
\QED \end{proof}

\begin{proof}[Proof of \corref{Corollary2.3.1}]
In the following relations, we use formula~\eqref{1.12}, \thmref{Theorem2.3}
(see \eqref{2.13}--\eqref{2.14}) and Fubini's Theorem:
\begin{align} \label{5.5}
\begin{split}
\left\|\fc{q(x)}{r(x)}y(x)\right\|_1
&=\int_{-\iy}^\iy\fc{q(x)}{r(x)}\left|\int_x^\iy\fc{1}{r(t)}
\exp\left(-\int_x^t\fc{q(s)}{r(s)}ds\right)f(t)dt\right|dx\\
&\le\int_{-\iy}^\iy\fc{q(x)}{r(x)}\left[\int_x^\iy\fc{1}{r(t)}
\exp\left(-\int_x^t\fc{q(s)}{r(s)}ds \right)|f(t)|dt\right]dx\\
&=\int_{-\iy}^\iy\fc{|f(t)|}{r(t)}\left[\int_{-\iy}^t\fc{q(x)}{r(x)}
\exp\left(-\int_x^t\fc{q(s)}{r(s)}ds\right)dx\right]dt\\
&\le\fc{1}{r_0}\int_{-\iy}^\iy|f(t)|
\bigg[\exp\left.\left(-\int_x^t\fc{q(s)}{r(s)}ds\right)
\right|_{x=-\iy}^{x=t}\bigg]dt \\
&=\fc{1}{r_0}\|f(t)\|_1.
\end{split}
\end{align}
Inequality \eqref{2.17} now follows from \eqref{5.5}, \eqref{1.1} and the
triangle inequality for norms.
\QED
\end{proof}

\begin{proof}[Proof of \corref{Corollary2.3.2}]\ {\tt Necessity}.

\s
Suppose that problem \eqref{1.1}--\eqref{1.2} is correctly solvable in
$L_1(\mathbb R).$
Then\linebreak $S_1=\iy$ and $M_1~<~\iy$ by \thmref{Theorem2.3} (see
\eqref{2.13}--\eqref{2.14}).
Hence the function $d(\cdot)$ (see \eqref{2.2}) is uniquely determined.
Therefore from \eqref{2.12} we get
\begin{align*}
\iy> M_1&=\sup_{x\in \mathbb R}\fc{1}{r(x)}\int_{-\iy}^x
\exp\left(-\int_t^x\fc{q(s)}{r(s)}ds\right)dt \\
&\ge\sup_{x\in \mathbb R}\fc{c}{r(x)}\ \fc{1}{c}\int_{x-d(x)}^x
\exp\left(-\int_t^x\fc{q(s)}{r(s)}ds\right)dt\\
&\ge \fc{1}{c}\sup_{x\in \mathbb R}\int_{x-d(x)}^x
\exp\left(-\int_{x-d(x)}^{x+d(x)}\fc{q(s)}{r(s)}ds\right)dt \\[0.1cm]
&\ge \fc{\exp(-2)}{c}\sup_{x\in \mathbb R}d(x)=c_1d_0.
\end{align*}
Since $S_1=\iy$ and $d_0<\iy,$ \lemref{lemmama2.1} implies $B(a)>0$ for some
$a\in (0,\iy)$ (see \eqref{2.11}).
This gives \eqref{1.6} (see the proof of \corref{Corollary2.2.1}, necessity).
\end{proof}

\begin{proof}[Proof of \corref{Corollary2.3.2}]\ {\tt Sufficiency}.

\s
Suppose that \eqref{1.6} holds.
Then by \eqref{2.12} we have \eqref{2.11} (see the proof of
\corref{Corollary2.2.1}, sufficiency).
Hence $S_1=\iy$ and $d_0<\iy$ because of \lemref{lemmama2.1}.
Let us verify that $M_1<\iy$ (see \eqref{2.15}).
In the following relations we use \lemref{lemma3.3}, \eqref{2.12} and
\eqref{4.32}:
\begin{align*}
M_1&=\sup_{x\in \mathbb R}\fc{1}{r(x)}\int_{-\iy}^x
\exp\left(-\int_t^x\fc{q(s)}{r(s)}ds\right)dt \\
&\le c\sup_{x\in\mathbb R}\sum_{n=-\iy}^{-1}\int_{\Dl_n}
\exp\left(-\int_t^x\fc{q(s)}{r(s)}ds\right)dt\\
&\le 2c \sup_{x\in \mathbb R}\sum_{n=-\iy}^{-1}d(x_n)\exp (-2(|n|-1))\\
&\le c_1d_0\sum_{n=1}^\iy\exp (-2(n-1))
<\iy.
\end{align*}
The corollary now follows from \thmref{Theorem2.3}.
\QED \end{proof}

\begin{proof}[Proof of \thmref{Theorem2.4}]\ {\tt Necessity}.

\s
If problem \eqref{1.1}--\eqref{1.2} is correctly solvable in $L_1(\mathbb R)$,
then $r_0>0$ because of \thmref{Theorem2.3}.
\end{proof}

\begin{proof}[Proof of \thmref{Theorem2.4}]\ {\tt Sufficiency}.

\s
Suppose that condition \eqref{2.11} holds.
Then from \lemref{lemmama2.1} we get $S_1=\iy$ and $d_0<\iy$ (see \eqref{2.1}
and \eqref{2.10}).
Then for $\varkappa(t)\equiv d(t)$,\ $t\in \mathbb R$, we have equalities
\eqref{3.17}, and by \lemref{lemma3.2} for
every $x\in \mathbb R$  there exist $R(x,d(\cdot))$-coverings of $(-\iy,x]$ and
$[x,\iy).$ By \thmref{Theorem2.3}, it remains to prove that $M_1<\iy$
(see \eqref{2.15}).
Below we use the properties of the $R(x,d(\cdot))$-covering of $(-\iy,x],$
\eqref{4.32} and \eqref{2.14}:
\begin{align*}
M_1&=\sup_{x\in \mathbb R}\fc{1}{r(x)}\int_{-\iy}^x
\exp\left(-\int_t^x\fc{q(s)}{r(s)}ds\right)dt \\
&\le\fc{1}{r_0}\sup_{x\in\mathbb R}\sum_{n=-\iy}^{-1}\int_{\Dl_n}
\exp\left(-\int_t^x\fc{q(s)}{r(s)}ds\right)dt\\
&\le \fc{2}{r_0}\sup_{x\in \mathbb R}\sum_{n=-\iy}^{-1}d(x_n)\exp(-2(|n|-1))\\
&\le\fc{2d_0}{r_0}\sum_{n=1}^\iy\exp(-2(n-1))
<\iy.
\end{align*}
Thus, \thmref{Theorem2.4} is proved. \QED
\end{proof}


\section{Proofs of the Main Results in the Case $p=\iy$}\label{PMRpiy}

In this section we prove Theorems \ref{Theorem2.5} and \ref{Theorem2.6}
and their corollaries.

\renewcommand{\qedsymbol}{}
\begin{proof}[Proof of \thmref{Theorem2.5}]\ {\tt Necessity}.

\s
Suppose that problem \eqref{1.1}--\eqref{1.2} is correctly solvable in $C(\mathbb R).$
Then its solution $y(\cdot)$ can be written in the form \eqref{1.12}, in view of
\lemref{lemma4.1}.
Hence for $f(x)\equiv1,$\ $x\in \mathbb R$, we get $y(x)=(G1)(x)\to0$ as $|x|\to\iy$
which coincides with \eqref{2.18}.
\end{proof}

\begin{proof}[Proof of \thmref{Theorem2.5}]\ {\tt Sufficiency}.

\s
From \eqref{1.3}--\eqref{1.4} and \eqref{2.18} it follows that $A_0<\iy$,
where
\begin{equation}\label{6.1}
A_0=\sup_{x\in \mathbb R} A(x),\quad A(x)=\int_x^\iy
\fc{1}{r(t)}\exp\left(-\int_x^t\fc{q(s)}{r(s)}ds\right)dt,\quad x\in \mathbb R.
\end{equation}
Let $\mu(\cdot)$ and $\theta(\cdot)$ be defined by \eqref{4.5}.
Then, since $A_0<\iy,$ from \eqref{3.12} it follows that the operator
$G:C(\mathbb R)\to C(\mathbb R)$ (see \eqref{1.12})
is bounded, and therefore \eqref{1.5} holds.
The equalities \eqref{1.2} follow from \eqref{2.18} and the following standard
estimate:
$$
|y(x)|=|(Gf)(x)|\le\left[\int_x^\iy
\fc{1}{r(t)}\exp\left(-\int_x^t\fc{q(s)}{r(s)}ds\right)dt\right]\cdot
\|f\|_{C(\mathbb R)},\quad x\in \mathbb R.
$$

Let us now prove that $S_1=\iy$ (see \eqref{2.1}).
Assume the contrary: $S_1<\iy.$
Then for every $x\le 0$ we have
\begin{align*}
\int_x^\iy\fc{1}{r(t)}\exp\left(-\int_x^t\fc{q(s)}{r(s)}ds\right)dt
\ge \exp(-S_1)\int_x^0\fc{dt}{r(t)}.
\end{align*}
From \eqref{2.18} and the above inequalities, it follows that
$$
0=\lim_{x\to-\iy}\int_x^\iy\fc{1}{r(t)}
\exp\left(-\int_x^t\fc{q(s)}{r(s)}ds\right)dt\ge\exp(-S_1)
\lim_{x\to -\iy}\int_x^0\fc{dt}{r(t)}>0,
$$
a contradiction. Hence $S_1=\iy.$ \thmref{Theorem2.5} is proved.
\QED
\end{proof}

\begin{proof}[Proof of \corref{Corollary2.5.1}]\ {\tt Necessity}.

Suppose that problem \eqref{1.1}--\eqref{1.2} is correctly solvable in
$C(\mathbb R)$ and \eqref{2.12} holds.
Then, by \thmref{Theorem2.5}, we have \eqref{2.18} and therefore for every
$a\in(0,\iy)$, we derive from
\eqref{2.12}
\begin{align*}
0&=\lim_{|x|\to\iy}\int_{x-a}^\iy\fc{1}{r(t)}
\exp\left(-\int_{x-a}^t\fc{q(s)}{r(s)}ds\right)dt \\
&\ge\fc{1}{c}\lim_{|x|\to\iy}\int_{x-a}^{x+a}
\exp\left(-\int_{x-a}^{x+a}\fc{q(s)}{r(s)}ds\right)dt\\
&\ge\fc{2a}{c}\lim_{|x|\to\iy}\exp\left(-c\int_{x-a}^{x+a}q(s)ds\right)\ge0
\end{align*}
which implies (1.7).
\end{proof}

\begin{proof}[Proof of \corref{Corollary2.5.1}]\ {\tt Sufficiency}.

\s
Suppose that \eqref{1.7} holds.
Then, obviously, for every $a\in(0,\iy)$ we have
\begin{align}\label{6.2}
\lim_{|x|\to\iy}\int_x^{x+a}q(t)dt=\iy.
\end{align}
It is also clear that there is $a_0\gg1$ such that
\begin{equation}\label{6.3}
\a:=\inf_{x\in \mathbb R}\int_{x-a_0}^{x+a_0}q(t)dt>0.
\end{equation}
Let $a_1=\max\{a_0,\a^{-1}\}.$
Then for all $x\in \mathbb R$ we have
\begin{equation}\label{6.4}
\int_{x-a_1}^{x+a_1}q(t)dt\ge\int_{x-a_0}^{x+a_0}q(t)dt\ge\a\ge\fc{1}{a_1}
\end{equation}
Let $\{\sigma_n\}_{n=1}^\iy $ be the system of segments constructed according
to the following rule:
\begin{tightitemize}
\item[1)] $\sigma_n=[\sigma_n^-,\sigma_n^+]=[x_n-a_1,x_n+a_1],\ n=1,2,\dots\ .$
\item[2)] $\sigma_{n+1}^-=\sigma_n^+,\ n=1,2,\dots;\ \sigma_1^-=x.$
\end{tightitemize}
Then, from \eqref{2.12} and \eqref{6.4} for $n\ge2$, it follows that
\begin{align*}
\int_{\sigma_1^-}^{\sigma_n^-}\fc{q(s)}{r(s)}ds
=\sum_{k=1}^{n-1}\int_{\sigma_k}\fc{q(s)}{r(s)}ds
\ge\fc{1}{c}\sum_{k=1}^{n-1}\int_{x_k-a_1}^{x_k+a_1}q(s)ds
\ge\fc{n-1}{ca_1}
\end{align*}
Thus for $n\ge 2$ and $\dl=(ca_1)^{-1}$, we have
\begin{equation}\label{6.5}
\int_{\sigma_1^-}^{\sigma_n^-}\fc{q(s)}{r(s)}ds\ge\dl(n-1),\quad
\delta>0.
\end{equation}
Clearly, \eqref{6.5} is also valid for $n=1.$
Let us verify that under conditions \eqref{2.12} and \eqref{1.7} we have
$A_0<\iy$ (see \eqref{6.1}).
Indeed,
\begin{align*}
A(x)&=\int_x^\iy\fc{1}{r(t)}\exp\left(-\int_x^t\fc{q(s)}{r(s)}ds\right)dt\\
&=\sum_{k=1}^\iy\int_{\sigma_k}\fc{1}{r(t)}
\exp\left(-\int_{\sigma_1^-}^t\fc{q(s)}{r(s)}ds\right)dt\\
&\le  2a_1 c\sum_{k=1}^\iy
\exp\left(-\int_{\sigma_1^-}^{\sigma_k^-}\fc{q(s)}{r(s)}ds\right) \\
&\le 2a_1 c \sum_{k=1}^\iy\exp(-\dl(k-1))=c_1<\iy.
\end{align*}

Let us now prove \eqref{2.18}.
Let $\ve>0.$
For $x\in \mathbb R$ from $A_0<\iy$ and \eqref{2.12}, it follows that
\begin{align} \label{6.6}
\begin{split}
A(x)&=\int_x^{x+\ve}\fc{1}{r(t)}\exp\left(-\int_x^t\fc{q(s)}{r(s)}ds\right)dt \\
&\phantom{= } +\exp\left(-\int_x^{x+\ve}\fc{q(s)}{r(s)}ds\right)A(x+\ve) \\
&\le c\ve+\exp\left(-\fc{1}{c}\int_x^{x+\ve}q(s)ds\right)A_0.
\end{split}
\end{align}
Now, from \eqref{6.6} and \eqref{6.2}, we get
$$
0\le \underset{|x|\to\iy}\varliminf A(x)\le\underset{|x|\to\iy}\varlimsup
A(x)\le  c\ve
$$
which implies $\lim_{|x|\to\iy}A(x)=0$.
Thus \eqref{2.18} holds.
It remains to refer to \thmref{Theorem2.5} to finish the proof of
\corref{Corollary2.5.1}.
\QED
\end{proof}

\m
\begin{proof}[Proof of \thmref{Theorem2.6}]\ {\tt Necessity}.

\s
Let $B(a)>0$ for some $a\in (0,\iy).$
Then by \lemref{lemmama2.1} we get $S_1=\iy$ and $d_0<\iy$ (see \eqref{2.1} and
\eqref{2.10}).
In the following relations, we use \thmref{Theorem2.5} (see \eqref{2.18}),
\eqref{2.2} and the
inequality
$d_0<\iy:$
\begin{align*}
0&=\lim_{|x|\to\iy}\int_{x-d(x)}^\iy\fc{1}{r(t)}
\exp\left(-\int_{x-d(x)}^t\fc{q(s)}{r(s)}ds\right)dt\\
&\ge\lim_{|x|\to\iy}\int_{x-d(x)}^{x+d(x)}\fc{1}{r(t)}
\exp\left(-\int_{x-d(x)}^t\fc{q(s)}{r(s)}ds\right)dt\\
&\ge\lim_{|x|\to\iy}
\exp\left(-\int_{x-d(x)}^{x+d(x)}\fc{q(s)}{r(s)}ds\right)
\cdot\int_{x-(d(x)}^{x+d(x)}\fc{dt}{r(t)} \\
&=\exp(-2)\lim_{|x|\to\iy}\int_{x-d(x)}^{x+d(x)}\fc{dt}{r(t)}\ge0.
\end{align*}
Thus \eqref{2.20} holds.
\end{proof}

\begin{proof}[Proof of \thmref{Theorem2.6}]\ {\tt Sufficiency}.

\s
Let $B(a)>0$ for some $a\in(0,\iy)$ and $\tilde A_0=0.$ Then $S_1=\iy$,\
$d_0<\iy $ by Lemma 2.3, $S_2=\iy$ (because $d_0<\iy$),  there is a
$R(x,d(\cdot))$-covering $\{\Delta_n\}_{n=1}^\iy$ of $[x,\iy),$ \
$x\in\mathbb
R$ (see Lemma 3.8) and, finally, $\tilde A=\|\tilde A(x)\|_{C(\mathbb
R)}<\iy.$
Furthermore, for a given $\varepsilon>0$, there is $x_0(\varepsilon)\gg1$
such
that $\tilde A(x)\le\varepsilon/4$ for $|x|\ge x_0(\varepsilon).$ Let $n_0$
be
a natural number such that $4\tilde A\exp(-2n_0)\le\varepsilon.$ Then the
following relations are fulfilled for $|x|\ge\tilde
x_0(\varepsilon)=x_0(\varepsilon)+2n_0d_0:$
$$\cupl_{k=1}^{n_0}\Delta_k\subseteq[x,x+2n_0d_0],\quad
[x,x+2n_0d_0]\cap[-x_0(\varepsilon),x_0(\varepsilon)]
=\emptyset.$$ Consequently, for $|x|\ge\tilde x_0(\varepsilon),$ we obtain:
\begin{equation}
\begin{aligned}\label{6.3}
A(x)&=\sum_{n=1}^\iy\int_{\Delta_n}\frac{1}{r(t)}\exp\left(-\int_x^t\frac{q(s)}{r(s)}ds\right)
dt\\ &\le\sum_{n=1}^\iy\exp(-2(n-1))\int_{\Delta_n}\frac{dt}{r(t)}\\
&=\sum_{n=1}^{n_0}\exp(-2(n-1))\int_{\Delta_n}\frac{dt}{r(t)}+
\sum_{n=n_0+1}^\iy
\exp(-2(n-1))\int_{\Delta_n}\frac{dt}{r(t)}\\
&\le \frac{\varepsilon}{4}\sum_{n=1}^\iy\exp(-2(n-1))+\tilde
A\exp(-2n_0)\sum_{n=1}^\iy\exp(-2(n-1))\le\varepsilon. \end{aligned}
\end{equation}
{}From \eqref{6.3} it follows that $A_0=0.$
 It remains to apply \thmref{Theorem2.5}.
\QED \end{proof}


\section{Correct non-solvability of the boundary value problem}
\label{Ptbvp}
In this section we prove \thmref{Theorem2.7}.

\begin{proof}[Proof of \thmref{Theorem2.7}]
Suppose that \eqref{2.11} holds.
Then $S_1=\iy$ and $d_0<\iy$ because of \lemref{lemmama2.1} (see \eqref{2.1}
and \eqref{2.10}).
Since we have \eqref{1.9}, let, for example, $q(x)\to0$ as $x\to\iy.$
Then for a given $\ve>0$, there is $x_0(\ve)$ such that $q(x)\le\ve$ for $x\ge
x_0(\ve).$
Below we assume $t\ge x_0(\ve)+d_0.$
We have $q(t)\le\ve$ for $t\ge x-d(x)$ since $x-d(x)\ge x_0(\ve)+d_0-d(x)\ge
x_0(\ve).$
To prove the theorem, we consider separately three cases:\ I)~$p\in(1,\iy);$\
II)~$p=1;$\ III)~$p=\iy.$

\s
{\tt Case I)}. Suppose that problem \eqref{1.1}--\eqref{1.2} is correctly solvable in
$L_p(\mathbb R).$
Then by \thmref{Theorem2.1} we have $A_{p'}<\iy.$
{}From \eqref{2.2} it follows from the definition of $A_{p'}$ that
\begin{align*}
2
\le\ve\int_{x-d(x)}^{x+d(x)}\fc{dt}{r(t)}
\le c\ve d(x)^{\frac{1}{p}}
\bigg[\int_{x-d(x)}^{x+d(x)}\fc{dt}{r(t)^{p'}}\bigg]^{\frac{1}{p'}}
\le c\ve d(x)^{\frac{1}{p}}A_{p'}^{\frac{1}{p'}}.
\end{align*}
Hence $d(x)\ge c^{-1}\ve^{-p}$ for $x\ge x_0(\ve)+d_0.$
This means that $d(x)\to\iy$ as $x\to\iy$ which is a
contradiction, since $d_0<\iy.$
Hence problem \eqref{1.1}--\eqref{1.2} is not correctly solvable in
$L_p(\mathbb R).$

\s
{\tt Case II)}. Suppose that problem \eqref{1.1}--\eqref{1.2} is correctly
solvable in $L_1(\mathbb R).$
Then $r_0>0$ by \thmref{Theorem2.3} (see \eqref{2.14}).
From \eqref{2.2} it follows that
$$
2=\int_{x-d(x)}^{x+d(x)}\fc{q(t)}{r(t)}dt\le\ve\int_{x-d(x)}^{x+d(x)}\fc{dt}{r(t)}\le \fc{2\ve}{r_0}
d(x).
$$
Hence $d(x)\ge r_0\ve^{-1}$ for $x\ge x_0(\ve)+d_0.$
Therefore, $d(x)\to\iy$ as $x\to\iy$ which is a
contradiction, since $d_0<\iy.$
Hence problem \eqref{1.1}--\eqref{1.2} is not correctly solvable
in $L_1(\mathbb R).$

\s
{\tt Case III)}. Suppose that \eqref{1.1}--\eqref{1.2} is correctly solvable
in $C(\mathbb R).$  Then
$$
2=\int_{x-d(x)}^{x+d(x)}\fc{q(t)}{r(t)}
\le\ve\int_{x-d(x)}^{x+d(x)}\fc{dt}{r(t)},
$$
for each $\ve>0$. This
contradicts   condition \eqref{2.20} of \thmref{Theorem2.6}.
Hence problem \eqref{1.1}--\eqref{1.2} is not correctly solvable in
$C(\mathbb R).$
\QED
\end{proof}


\section{Example}\label{Example}

Let $\alpha,\beta\in \mathbb R,$\ $\gamma>0,$ and
\begin{equation}\label{8.1}
r(x)=e^{\a|x|},\quad q(x)=e^{\beta|x|}+e^{\beta|x|}\cos e^{\g|x|}, \quad
x\in
\mathbb R.
\end{equation}
In this section, we study a condition of correct solvability in
$L_p(\mathbb
R)$,\ $p\in[1,\iy]$, for the problem \eqref{1.1} -- \eqref{1.2} in the case
\eqref{8.1}. (Below for brevity, we write ``problem \eqref{8.1}").

\subsection{Necessary conditions for correct solvability}

\begin{lemma}\label{lemma8.1}
The problem \eqref{8.1}, is correctly solvable in $L_p(\mathbb R)$,\
$p\in[1,\iy]$, only if $\beta\ge\alpha$ and $ \beta\ge0.$ In the latter
case,
$S_1=\iy$,\ $d_0<\iy$ (see \eqref{2.1}, \eqref{2.10}).
\end{lemma}

\begin{proof}
By Theorems 2.1, 2.5 and 2.8, the problem \eqref{1.1} -- \eqref{1.2} is
correctly solvable in $L_p(\mathbb R)$,\ $p\in[1,\iy]$ only if $S_1=\iy.$
It is
clear that in our case $S_1=\iy$ only if $\beta\ge\alpha.$ We shall show
that
in this case \eqref{2.11} is fulfilled.  Let $a=3/2\gamma,$\ $x\ge a$ and
$\xi$
is some
point on the segment $[x-a,x+a].$ We denote \begin{equation}\label{8.2}
f(t)=e^{\theta t}+e^{\theta t}\cos e^{\g t},\quad t\ge 0,\quad
\theta=\beta-\alpha\ge0. \end{equation} Below we use the mean value theorem
\cite[\S12.3]{Tit}:
\begin{align*}
F(x)&=\int_{x-a}^{x+a}\frac{q(t)}{r(t)}dt =\int_{x-a}^{x+a}f(t)dt\ge
\int_{x-a}^{x+a}(1+\cos e^{\gamma
t})dt\\
&=2a+\frac{1}{\gamma}\int_{x-a}^{x+a}e^{-\gamma t}[\gamma e^{\gamma t}\cos
e^{\gamma t}]dt=2a+\frac{1}{\gamma e^{x-a}}\left(\left.(\sin e^{\gamma
t}\right|_{x-a}^\xi\right)\\ &\ge 2a-\frac{2}{\gamma}=\frac{1}{\gamma}>0.
\end{align*}

For $z\in[0,a]$, the function $F(z)$ is continuous and positive.  Therefore
$F_0(a)=\min\limits_{z\in [0,a]}F(z)>0.$
 Let $\a=\min\left\{\gamma^{-1}, F_0(a)\right\}.$
Then $F(z)\ge\a$ for all $z\ge0.$ The case $z\le 0$ can be considered in a
similar way.
 Thus, \eqref{2.11} holds.  By Lemma 2.3, we have $d_0<\iy.$
Moreover, $\beta\ge0$ by Theorem 2.10, and therefore the lemma is proved.
\QED
\end{proof}

\subsection{Estimates of the auxiliary function of growth}

\quad{}

 Below we assume that conditions $\beta>\alpha,$\ $\beta\ge0$ hold. We
shall
establish inequalities for $d(\cdot),$ on $\mathbb R $ in three separate
cases:
$\theta<\g,$\ $\theta=\g,$\ $\theta>\g$. We need the following notation.
Let
the functions $\vp(\cdot)$ and $\psi(\cdot)$ be positive and continuous on
$\mathbb R.$   We write $\vp(x)\asymp \psi(x)$ if there is $c\in [1,\iy)$
such
that $c^{-1}\vp(x)\le\psi(x)\le c\vp(x),$\ $x\in\mathbb R.$ Moreover, in
the
proofs, we consider only the case $x\ge0$ where the functions $r(\cdot)$
and
$q(\cdot)$ in \eqref{8.1} are even.

\begin{lemma}\label{lemma8.4}
Let $0\le\theta\le \g,$\ $\theta=\beta-\a.$ Then
\begin{equation}\label{8.3}
d(x)\asymp e^{-\theta|x|},\quad x\in\mathbb  R.
\end{equation}
\end{lemma}

\begin{proof}
To apply \thmref{Theorem2.8}, set
$$
q_1(x)=e^{\beta |x|},\quad q_2(x)=e^{\beta|x|}\cos(e^{\g|x|}),\quad
x\in\mathbb
R.
$$
Let $\theta<\g.$ Let us estimate $\varkappa_1(x),$(see \eqref{2.22} for
$x\to\iy:$
 \begin{align*}
 \varkappa_1(x)&=\sup_{|z|\le2e^{-\theta
 x}}\left|\int_0^z[e^{\theta(x+t)}-2e^{\theta
 x}+e^{\theta(x-t)}]dt\right|\\
 &=e^{\theta x}\sup_{|z|\le 2e^{-\theta x}} \sup_{|z|\le 2e^{-\theta
 x}}\left|\int_0^z[\theta^2t^2+\dots]dt\right|
  \le ce^{-2\theta x}.
 \end{align*}
Below, when estimating $\varkappa_2(x),$\ $x\gg1$ we use the    mean value
theorem \cite[\S12.3]{Tit}
\begin{equation}
 \begin{aligned}\label{8.4}
  \varkappa_2(x)&=\sup_{|z|\le 2e^{-\theta
x}}\left|\int_{x-z}^{x+z}e^{\theta
   t}\cos(e^{\g t})dt\right|\\
 &=\sup_{|z|\le 2e^{-\theta
 x}}\left|\int_{x-z}^{x+z}e^{(\theta-\g)t}[\gamma e^{\g t}\cos(e^{\g
t}]dt\right|
   \le ce ^{(\theta-\g)x}.
 \end{aligned}
 \end{equation}
Since $\theta<\g$, we have $\varkappa_1(x)\to0,$\ $\varkappa_2(x)\to0$ as
$x\to\iy.$
 By \thmref{Theorem2.8} this implies \eqref{8.3}.
Consider now the cases $\theta=0$ and $\theta=\g.$   If $\theta=0$, then
$\varkappa_1(x)=0$ and \eqref{8.4} holds. Therefore, we obtain \eqref{8.3}
as
above. Let $\theta=\gamma.$ We set $\eta(x)=(1+\gamma^{-1})e^{-\gamma x},$\
$x\gg1.$ Then
\begin{align*}
 \int_{x-\eta(x)}^{x+\eta(x)}f(t)dt& =e^{\g x}
\fc{e^{\g\eta(x)}-e^{-\g\eta(x)}}{\g}
+\left.\fc{\sin(e^{\g x})}{\g}\right|_{x-\eta(x)}^{x+\eta(x)}\\
&\ \ge \fc{e^{\g
x}}{\g}\left[2\gamma\eta(x)+2\fc{(\g\eta(x))^3}{3!}+\dots\right]-\fc{2}{\g}\ge
  2.
\end{align*}
{}From here $d(x)\le\eta(x)$ for $x\gg1$ (see \eqref{3.28}). Let
$\ve=\min\{4^{-1},(2\g)^{-1}\},$\ $\eta(x)= \ve e^{-\g x}.$ Then, for
$x\gg1,$
we obtain:
\begin{align*}
  \int_{x-\eta(x)}^{x+\eta(x)}f(t)dt&\le 2\int_{x-\eta(x)}^{x+\eta(x)}e^{\g
t}dt
 =  \frac{4}{\g}e^{\g
x} \left[\g\eta(x)+\frac{(\g\eta(x))^3}{3!}+\cdots\right]\\
&\le \frac{4}{\g}e^{\g x}\frac{\g\eta(x)}{1-(\g\eta(x))^2}\le
\frac{16}{3}e^{\g\eta(x)}\le \frac4{3}<2.
\end{align*}
Hence $d(x)\ge\eta(x)$ for $x\gg1$ (see \eqref{3.28}).  This implies
\eqref{8.3}.
\QED \end{proof}

\begin{lemma}\label{lemma8.3}
Let $\theta> \g.$ Denote by $\{x_k\}_{k=-\iy}^\iy$ points such that
$|x_k|=$\linebreak $ \g^{-1}ln[(2|k|+1)\pi],$\ $k=0,\pm1,\pm2,\dots$ Then
\begin{equation}\label{8.5}
d_k\doe d(x_k)\asymp(2|k|+1)^{-\fc{\theta+2\g}{3\g}},\quad
k=0,\pm1,\pm2,\dots
\end{equation}
\end{lemma}

\begin{proof}
Let   $k\gg1,$\ $d\in[0,1],$\ $ t\in\omega_k=[x_k-d,x_k+d].$ Since
$f(x_k)=f'(x_k)=0$ (see \eqref{8.2}), by Taylor's formula we obtain
\eqref{8.6}
and as a consequence \eqref{8.7}:
\begin{equation}\label{8.6}
f(t)=
\fc{f''(x_k)}{2!}(t-x_k)^2+\fc{f'''(\xi(t))}{3!}(t-x_k)^3,\quad\xi(t)\in\omega_k,
\end{equation}
\begin{equation}\label{8.7}
\int_{\omega_k}f(t)dt=
\fc{f''(x_k)}{3}d^3+\frac{1}{3!}\int_{\omega_k}f'''(\xi(t)) (t-x_k)^3dt.
\end{equation}
The following relations are obvious: \begin{equation}\label{8.8}
f''(x_k)=\g^2[(2k+1)\pi]^{\frac{\theta+2\g}{\g}}\asymp
k^{\frac{\theta+2\g}{\g}},\quad k\gg1, \end{equation}
\begin{equation}\label{8.9}
|f'''(\xi)|\le ce^{(\theta+3\g)\xi}\le ce^{(\theta+3\g)x_k}\le
ck^{\fc{\theta+3\g}{\g}} ,\quad k\gg1,\quad \xi\in\omega_k.
\end{equation}
\begin{equation}\label{8.10}
\left|\int_{\omega_k}f'''(\xi(t))(t-x_k)^3dt\right|\le
ck^{\frac{\theta+3\g}{\g}}d^4.
\end{equation}

Let $\mu\in(0,\iy)$ and $\eta_k=\mu[f''(x_k)]^{-1/3}.$ Then due to
\eqref{8.7},
\eqref{8.8} and \eqref{8.10}, we get for $\mu=(9/2)^{1/3}$ and $k\gg1:$
$$\int_{x-\eta_k}^{x+\eta_k}f(t)dt\le\frac{3}{2}+ck^{-\frac{\theta-\g}{3\g}}<2.$$
Hence, $d_k\ge\eta_k$ by \lemref{lemma3.4}. Similarly, we have for
$\mu=9^{1/3};$
$$
\int_{x_k-d}^{x_k+d}f(t)dt\ge 3-ck^{-\frac{\theta-\gamma}{3\gamma}}\ge
2,\quad
k\gg1,$$ and hence $d_k\le \eta_k$ by \lemref{lemma3.4}. The relation
\eqref{8.5} then follows.
\QED \end{proof}

\begin{lemma}\label{lem8.4}
The following inequalities hold:
\begin{equation}\label{8.11}
d(x)\le cd_k\quad\text{for}\quad x\in[x_k,x_{k+1}],\quad
k=0,\pm1,\pm2,\dots
\end{equation}
\end{lemma}
To prove \eqref{8.11} (for $x\ge0)$, we need Lemmas \ref{lem8.5} --
\ref{lem8.11} below.  The case $x\le 0$ can be studied in a similar way.

\begin{lemma}\label{lem8.5}
For all $k=0,1,2,\dots,$ the function $f(t)$ (see \eqref{8.2}) has a unique
extremum (maximum) on the interval $(x_k,x_{k+1}).$ If $\tilde z_k$ is a
corresponding extreme point, then
\begin{equation}\label{8.12}
z_k<\tilde z_k<z_k+\mu k^{-2}\quad\text{for}\quad k\gg1.
\end{equation}
Here $\mu$ is some absolutely positive constant,
$z_k=\g^{-1}ln[(2k+2)\pi],$\
$k=0,1,2,\dots $
\end{lemma}
\begin{proof}

The equality $f'(t)=0$ can be easily brought  to the form
\begin{equation}\label{8.13}
\g\left(\cos\frac{1}{2}e^{\g t}\right)^2e^{-\g t}\varphi(t)=0,\quad
\varphi(t)\doe \theta\g^{-1}e^{-\g t}-tg\left(\frac{1}{2}e^{\g t}\right).
\end{equation}
Since $f(t)>0,$\ $t\in(x_k,x_{k+1})$ and $f(x_k)=f(x_{k+1})=0,$ then $f(t)$
has
maximum on the interval $(x_k,x_{k+1}).$ Furthermore, the first two factors
in
\eqref{8.13} are positive on $(x_k,x_{k+1}),$ and $\varphi'(t)<0$ for
$t\in(x_k,x_{k+1}),$\ $\varphi(x_k)\varphi(x_{k+1})<0.$ This means that
$\varphi(t)$ has a unique root on $(x_k,x_{k+1})$, as desired. The lower
estimate \eqref{8.12} follows from an inequality $f'(z_k)>0.$ Let $k\gg1,$
\
$\hat z_k=z_k+\mu k^{-2}.$ By Taylor's formula, we obtain \eqref{8.14} with
$\xi\in(z_k,\hat z_k):$
\begin{equation}\label{8.14}
f'(\hat z_k)=f'(z_k)\left[1+\frac{f''(z_k)}{f'(z_k)}\
\frac{\mu}{k^2}+\frac{1}{2!}\ \frac{f'''(\xi)}{f'(z_k)}\
\frac{\mu^2}{k^4}\right].
\end{equation}
Together with the obvious relation \begin{equation}\label{8.15}
f'(z_k)\asymp
k^{\frac\theta\g},\quad |f''(z_k)|\asymp
k^{\frac{\theta+2\gamma}{\gamma}},\quad f''(z_k)<0,\quad k\gg1,
\end{equation}
we use \eqref{8.14} and \eqref{8.9} and obtain
$$f'(\hat z_k)\le f(z_k)(1-c^{-1}\mu+c\mu^2k^{-1}).$$
Clearly, $\hat f'(\hat z_k)<0$ for $\mu=2c$ and $k\gg1$.  Therefore, the
upper
estimate \eqref{8.12} is true.
\QED \end{proof}

\begin{lemma}\label{lem8.6} For all $k\gg1,$ the following inequality
holds:
\begin{equation}\label{8.16}
f(x_k+d_k)\ge f(x_k-d_k).
\end{equation}
\end{lemma}
\begin{proof}
The inequality \eqref{8.16} is equivalent to (see \eqref{8.2}):
 \begin {equation}\label{8.17}
 e^{\theta d_k}[1+\cos e^{\g(x_k+d_k)}]\ge e^{-\theta d_k}[1+\cos
  e^{\g(x_k-d_k)}],\quad k\gg1.
  \end{equation}
{}From \eqref{8.17} and from the following obvious equalities
$$\cos[e^{\g(x_k+d_k)}]=-\cos[(2{k+1})\pi(e^{\g d_k}-1)],\quad k\gg1,$$
$$\cos[e^{\g(x_k-d_k)}]=-\cos[(2{k+1})\pi(e^{-\g d_k}-1)],\quad k\gg1,$$
it follows that \eqref{8.17} is equivalent to
\begin{equation}\label{8.18}
e^{\theta d_k}\ge\left|\frac{\sin\left[\left(k+\frac{1}{2}]\right)\pi(e^{\g
d_k}-1)e^{-\g d_k}]\right]}{\sin\left[\left(k+\frac{1}{2}\right)\pi(e^{\g
d_k}-1)\right]}\right|,\quad k\gg1.
\end{equation}
In this connection, we note that $$0<\left(k+\frac{1}{2}\right)\pi(e^{\g
dk}-1)\le ckd_k\le ck^{-\frac{\theta-\g}{3\g}}\to0\quad\text{as}\quad
k\to\iy.$$ This means that in \eqref{8.18}, the arguments of both sines
tend to
zero (as $k\to\iy)$ and are positive. Therefore, \eqref{8.18} follows from
the
monotonicity of the function $\sin x$ in the neighborhood of the point
$x=0.$
\QED \end{proof}

\begin{lemma}\label{lem8.7}
The following equalities hold:
\begin{equation}\label{8.19}
\lim_{k\to\iy}\frac{\tilde z_k-x_k}{d_k}=\iy,\qquad
\lim_{k\to\iy}\frac{x_{k+1}-\tilde z_k}{d_k}=\iy.
\end{equation}
\end{lemma}
\begin{proof}
Below we use \eqref{8.5} and \eqref{8.12}
$$
\frac{\tilde z_k-x_k}{d_k} >\frac{z_k-x_k}{d_k}=\frac{1}{\g}\
\frac{ln\left(1+\frac{1}{2k+1}\right)}{d_k}\ge
ck^{\frac{\theta-\gamma}{3\g}}\to\iy\quad\text{as}\quad k\to\iy,$$
 \begin{align*}\frac{x_{k+1}-\tilde z_k}{d_{k+1}}&\ge
 \frac{x_{k+1}-z_k-\mu k^{-2}}{d_{k+1}}\ge \frac{1}{\g}
  \frac{ln\left(1+\frac{1}{2k+2}\right)}{d_k}-\frac{\mu}{k^2d_k}
 \\
&\ge c^{-1}k^{\frac{\theta-\g}{3\g}}\to\iy\quad\text{as}\quad k\to \iy.
\qquad \QED
\end{align*}
 \end{proof}

\begin{corl}\label{cor8.8}
 For all $k\gg1,$ the function $f(t)$ decreases monotonically on the
segment
 $[x_k-d_k,x_k]$ and increases monotonically on the segment
$[x_k,x_k+d_k].$
\end{corl}
 \begin{lemma}\label{lem8.9}
 The following equality holds:
 \begin{equation}\label{8.20}
 \lim_{k\to\iy}\frac{f(\tilde z_k)}{f(z_k)}=1.
 \end{equation}
 \end{lemma}
 \begin{proof}
The following relations are obvious:
$$f(\tilde z_k)=f(z_k)+f'(\xi)(\tilde z_k-z_k),\qquad \xi\in(z_k,\tilde
z_k)$$
$$|f'(\xi)|\le ce^{(\theta+\g)\xi}\le ce^{(\theta+\g)z_k}\le
ck^{\frac{\theta+\g}{\g}},\quad k\gg 1.$$ {}From here and taking into
account
\eqref{8.12}, we obtain:
$$\left|\frac{f(\tilde z_k)}{f(z_k)}-1\right|\le c\frac{|f'(\xi)|}{f(z_k)}\
\frac{1}{k^2}\le
ck^{\frac{\theta+\g}{\g}-\frac{\theta}{\g}-2}=\frac{c}{k}\to0\quad\text{as}\quad
k\to\iy. \qquad \QED $$
\end{proof}
\begin{lemma}\label{lem8.10}
For $k\gg1,$ the following inequalities hold:
\begin{equation}\label{8.21}
f(\tilde z_k+d_k)\ge f(x_k+d_k),\qquad f(\tilde z_k-d_{k+1})\ge
f(x_{k+1}+d_{k+1}).
\end{equation}
\end{lemma}
\begin{proof}
Both inequalities can be verified in the same way. We prove, for example,
the
first one.  Below we use Taylor's formula
$$f(\tilde z_k+d_k)=f(\tilde z_k)+\frac{f''(\xi)}{2!}d_k^2,\qquad
\xi\in(\tilde
z_k,\tilde z_k+d_k),$$  and the estimates of type \eqref{8.15} and
equality
\eqref{8.20}:
\begin{align*}
f(\tilde z_k+d_k)&\ge f(\tilde z_k)\left(1-\frac{|f''(\xi)|}{2f(\tilde
z_k)}d_k^2\right)\ge
c^{-1}f(z_k)\left(1-ck^{\frac{\theta+\g}{\g}-\frac{\theta}{\g}-\frac{2}{3}\
\frac{\theta+2\g}{\g}}\right)\\
&\ge c^{-1}k^{\frac{\theta}{\g}}\left(1-ck^{ -\frac{2}{3}\
\frac{\theta-\g}{\g}}\right)\ge c^{-1}k^{\frac{\theta}{\g}},\quad k\gg1.
\end{align*}
Below we again use Taylor's formula:
$$f(x_k+d_k)=f(x_k)+f'(x_k)d_k+\frac{f''(\xi)}{2!}d_k^2=\frac{f''(\xi)}{2!}d_k^2,\quad
\xi\in(x_k, x_k+d_k)$$ and estimates of type \eqref{8.8}:
 $$f(x_k+d_k)=\frac{f''(\xi)}{2!}d_k^2\le
ck^{\frac{\theta+2\g}{\g}-\frac{2}{3}\
  \frac{\theta+2\g}{\g}} =c k^{\frac{\theta+2\g}{3\g}},\quad k\gg1.
$$
Based on what we have found, we obtain for $k\gg1:$
$$f(\tilde z_k+d_k)\ge c^{-1}k^{\frac{\theta}{\g}}\ge
ck^{\frac{\theta+2\g}{3\g}}\ge f(x_k+d_k). \qquad \QED$$
\end{proof}

\begin{lemma}\label{lem8.11} The following inequality is true:
\begin{equation}\label{8.22}
f(x_k-2d_k)\ge f(x_k+d_k),\qquad k\gg1.
\end{equation}
\end{lemma}

\begin{proof}
The following relations are obvious:
$$\frac{2f(x_k-2d_k)}{f''(x_k)d_k^2}=4+\frac{8}{3}\
\frac{f'''(\xi_1)}{f''(x_k)}d_k\ge 4-\frac{8}{3}\
\frac{|f'''(\xi_1)|}{f''(x_k)}d_k,$$
$$\frac{2f(x_k+d_k)}{f''(x_k)d_k^2}=1+\frac{f''(\xi_2)}{3f''(x_k)}d_k\le
1+\frac{|f'''(\xi_2)|}{3f''(x_k)}d_k.$$ Here $\xi_1\in(x_k-2d_k,x_k),$\
$\xi_2\in(x_k,x_k+d_k).$ Clearly, \eqref{8.22} holds if
\begin{equation}\label{8.23}
4-\frac{8}{3}\ \frac{|f'''(\xi_1)|}{f''(x_k)}d_k\ge
1+\frac{|f'''(\xi_2)|}{f''(x_k)}d_k,\quad k\gg1.
\end{equation}
 {}From
\eqref{8.5}, \eqref{8.8} and \eqref{8.9}, it follows that \eqref{8.23} is
indeed true for all $k\gg1.$
\QED \end{proof}

 \begin{proof}[Proof of \lemref{lem8.4}]
 Below we consider \eqref{8.11} separately in the cases\linebreak
I)~$x\in[x_k,\tilde
 z_k]$ and
 II)~$x\in[\tilde z_k,x_{k+1}].$
 Furthermore, we assume that $k\ge
 k_0\gg1$.  Here $k_0$ is chosen so that for $k\ge k_0$ it is possible to
use
 Lemmas \ref{lem8.5} -- \ref{8.11}. It should be noted that inequality
 \eqref{8.11} is obvious for $k\le k_0.$

\textbf{Case  I)}\ If $x\in[x_k,x_k+2d_k],$ then
$[x_k-d_k,x_k+d_k]\subseteq[x-3d_k,x+3d_k]$; and therefore
$$\int_{x-3d_k}^{x+3d_k}f(t)dt\ge\int_{x_k-d_k}^{x_k+d_k}f(t)dt= 2.$$
{}From here $d(x)\le 3d_k$ by Lemma 3.11.  If $x\in[x_k+2d_k,\tilde z_k],$
then
$[x-d_k,x+d_k]\subseteq [x_k+d_k,\tilde z_k+d_k].$ Hence if
$t\in[x-d_k,x+d_k]$
and $\xi\in[x_k-d_k,x_k+d_k],$ then $f(t)\ge f(\xi)$ by Lemmas
\ref{lem8.6}, \
    \ref{lem8.10} and \corref{cor8.8}.  Thus, we get
$$\int_{x-d_k}^{x+d_k}f(t)dt\ge\int_{x_k-d_k}^{x_k+d_k}f(\xi)d\xi=2,$$ and
$d(x)\le d_k$ by Lemma 3.11.

\textbf{Case  II)}\ If $x\in[x_{k+1}-3d_{k+1},x_{k+1}],$ then
$[x_{k+1}-d_{k+1},x_{k+1}+d_{k+1}]\subseteq [x-4d_{k+1},x+4d_{k+1}];$ and
we
have
$$\int_{x-4d_{k+1}}^{x+4d_{k+1}}f(t)dt\ge
\int_{x_{k+1}-d_{k+1}}^{x_{k+1}+d_{k+1}}f(t)dt=2.$$ {}From here it follows
that
$d(x)\le 4d_{x+1}$ by Lemma 3.11. If $x\in[\tilde z_k,x_{k+1}-3d_{x+1}],$
then
$[x-d_{k+1},x+d_{k+1}]\subseteq [\tilde z_k-d_{k+1},x_{k+1}-2d_{k+1}].$
Hence
if $t\in[x-d_{x+1},x+d_{k+1}]$ and
$\xi\in[x_{k+1}-d_{k+1},x_{k+1}+d_{x+1}],$
then $f(t)\ge f(\xi)$ by Lemmas 8.6, 8.9, and 8.10 and by Corollary 8.7.1.
{}From here we obtain
$$\int_{x-d_{k+1}}^{x+d_{k+1}}f(t)dt\ge\int_{x_{k+1}-d_{k+1}}^{x_{k+1}+d_{k+1}}f(\xi)d\xi=2$$
and $d(x)\le d_{k+1}$ by Lemma 3.11.
\\
It remains to note that $d_k\asymp d_{k+1} $ (see \eqref{8.5}). \QED
\end{proof}

\subsection{Precise conditions for correct solvability of problem
\eqref{8.1}}

{}\quad

  Below we study problem (8.1) with the help  of
Theorems 2.4, 2.7 and 2.9 and Lemmas 2.3 and 8.1. Since the requirement
$d_0<\iy$ in the case \eqref{8.1} is fulfilled ``by necessity" (see Lemma
8.1),
then $r(t)\asymp r(x)$ for $|t-x|\le d(x),$\ $x\in\mathbb R$ (see
\eqref{2.10}
and \eqref{8.1}). We use these relations   together with conditions
$\beta\ge\alpha$,\ $\beta\ge0$ (see Lemma 8.1), without additional
stipulation.

\begin{theorem}\label{thm8.12} Let $p\in(1,\iy).$ Then problem \eqref{8.1}
is
correctly solvable in $L_p(\mathbb R)$ if and only if   one of the
following
conditions hold:
\begin{alignat}{2}
&1)\ && \beta=\alpha=0;\label{8.24}\\
&2)\ && \beta>0,\quad \gamma\ge\beta-\alpha>0,\quad p\ge
1-\frac{\alpha}{\beta};\label{8.25}\\
&3)\ && \beta>0,\quad
\beta-\alpha>\gamma,\quad\beta+2\alpha+2\gamma>0,\quad
p>1-\frac{3\alpha}{\beta+2\alpha+2\gamma}.\label{8.26} \end{alignat}
\end{theorem}

\begin{proof} Let $\gamma\ge\beta-\alpha.$ Then \eqref{2.7} and \eqref{8.3}
imply:
$$A_{p}'(x)=\int_{x-d(x)}^{x+d(x)}\frac{dt}{r(t)}\asymp\frac{d(x)}{e^{\alpha
p'|x|}}\asymp \frac{1}{e^{(\beta-\alpha+\alpha p')|x|}}.$$ Therefore,
$A_p'<\iy$ if and only if $\beta-\alpha+\alpha p'\ge0.$ Then, by Theorem
2.4,
we have the following relations:
 \begin{equation*}
\left.\begin{aligned}
  \beta\ge\a,\ \beta\ge0\quad  \\
 \beta-\a-\g\le0\ \  \\
 \beta-\a+\a p'\ge0\end{aligned}\right\}\Rightarrow \begin{cases}\text{a)}\
\text{if}\
 \beta=0\ \Rightarrow\ \a\le 0,
\a(p'-1)\ge0\ \Rightarrow \beta= \a \\
\quad \\
\text{b)}\ \text{if}\ \beta>0\ \Rightarrow\ 0\le\beta-a\le \g,\  p\ge
1-\fc{\a}{\beta}.\end{cases}
\end{equation*}
Thus conditions \eqref{8.24} and \eqref{8.25} are obtained.  Let now
$\beta-\alpha>\gamma.$ We use \eqref{2.6}, \eqref{2.7} and \eqref{8.5} to
obtain
\begin{align*}
A_p'\ge\sup_{|k|\ge0}A_p'(x_k)&=\sup_{|k|\ge0}\int_{x_k-d_k}^{x_k+d_k}\frac{dt}{e^{\alpha
p'|t|}}\ge c^{-1}\sup_{|k|\ge0}\frac{d_k}{e^{\alpha p'|x_k|}}\\
&\ge
c^{-1}\sup_{|k|\ge0}(2|k|+1)^{-\frac{\theta+2\gamma}{3\gamma}-\frac{\alpha
p'}{\gamma}}.
\end{align*}
On the other hand, if $x\in[x_k,x_{k+1}],$ then from \eqref{8.11} and
\eqref{8.5} it follows that
\begin{align*}&A_p(x)=\int_{x-d(x)}^{x+d(x)}\frac{dt}{e^{\alpha p'|t|}}\le
c\frac{d(x)}{e^{\alpha p'|x|}}\le c\frac{d_k}{e^{\alpha p'|x_k|}}\le
  c(2|k|+1)^{-\frac{\theta+2\gamma}{3\gamma}-\frac{\alpha p'}{\gamma}}
  \\
 &\Rightarrow A_p'\le
   c\sup_{|k|\ge0}(2|k|+1)^{-\frac{\theta+2\gamma}{3\gamma}-\frac{\alpha
  p'}{\gamma}}.
\end{align*}
Thus, $A_p'<\iy$ if and only if the following relations are fulfilled:
\begin{align*}
\left.\begin{aligned}
&\beta\ge0,\   \beta-\a>\g  \\
 &\fc{\beta-\alpha }{3\g}+\fc{\a p'}{\g}\ge0\end{aligned}\right\}
 \quad\Rightarrow\quad
\begin{cases}
  \beta>0, \beta-\a>\g, \\
 \beta+2\a+2\g>0\\
p\ge 1-\fc{3\a}{\beta+2\a+2\g}\end{cases}\quad\Rightarrow\quad\eqref{8.26}.
\end{align*}
It remains to quote Theorem 2.4.
\QED \end{proof}

\begin{theorem}\label{thm8.13} The problem \eqref{8.1} is correctly
solvable in
$L_1(\mathbb R)$ if and only if $\beta\ge\alpha\ge0$. \end{theorem}

\begin{proof} This statement follows from Lemma 8.1 and Theorem 2.7.
\QED \end{proof}

\begin{theorem}\label{thm8.14}
The problem \eqref{8.1} is correctly solvable in $C(\mathbb R)$ if and only
if
either one of the following conditions hold:
\begin{alignat}{2}
&1)\ && \beta> 0,\quad
\g\ge\beta-\alpha\ge0;\qquad\qquad\qquad\qquad\qquad\qquad\qquad\qquad\label{8.27}\\
&2)\ && \beta>0,\quad \beta-\alpha>\g,\quad
\beta+2\alpha+2\gamma>0.\label{8.28}  \end{alignat}
\end{theorem}

\begin{proof}
Let $\g\ge\beta-\alpha.$ Then as above, we obtain
$$\int_{x-d(x)}^{x+d(x)}\frac{dt}{e^{\alpha|t|}}\asymp
\frac{d(x)}{e^{\alpha|x|}}\asymp \frac{1}{e^{\beta|x|}},\quad x\in\mathbb
R.$$
Therefore, $\tilde A_0=0$ (see Theorem 2.9) if and only if $\beta>0$ and
\eqref{8.27} is fulfilled. Let $\beta-\alpha>\gamma.$ If $\tilde A_0=0,$
then
(see \eqref{8.5})
\begin{align*}
0=\lim_{|x\|to\iy}\int_{x-d(x)}^{x+d(x)}\frac{dt}{e^{\alpha|t|}}&=\lim_{|k|\to\iy}\int_{x_k-d_k}^{x_k+d_k}
\frac{dt}{e^{\alpha|t|}}\ge
c^{-1}\lim_{|k|\to\iy}\frac{d_k}{e^{\alpha|x_k|}}\\
&\ge
c^{-1}\lim_{|k|\to\iy}(2|k|+1)^{-\frac{\beta-\alpha+2\g}{3\g}-\frac{\alpha}{\g}}\ge0.\end{align*}
This implies $\beta+2\alpha+2\g>0.$ On the other hand, for
$x\in[x_k,x_{k+1}]$
by \eqref{8.5} and \eqref{8.11}, we obtain
$$\int_{x-d(x)}^{x_d(x)}\frac{dt}{e^{\alpha|t|}}\le
c\frac{d(x)}{c^{\alpha|x|}}\le c\frac{d_k}{e^{\alpha|x_k|}}\le
c(2|k|+1)^{-\frac{\beta-\alpha+2\g}{3\g}-\frac{\alpha}{\g}},$$ and
therefore
the condition $\beta+2\alpha+2\gamma>0$ must fulfill the equality $\tilde
A_0=0.$ Thus, $\tilde A_0=0$ if and only if
\begin{align*}
\left.\begin{aligned}
&\beta\ge0,\   \beta-\a>\g  \\
 & \beta+2\alpha +2\g \end{aligned}\right\}
 \quad\Rightarrow\quad
\begin{cases}
  \beta>0, \beta-\a>\g, \\
 \beta+2\a+2\g>0\end{cases}\quad\Rightarrow\quad \eqref{8.28}.
\end{align*}
It remains to use Theorem 2.9.
\QED \end{proof}


\section{Appendix}

This section contains the proof of \thmref{theorem1.2} (see \S$1$).
In addition, we present here Theorems \ref{theorem9.1} and \ref{theorem9.3}
which arise in connection with
the analysis of \thmref{theorem1.2} and can be viewed as a complement to
Theorems \ref{Theorem2.6} and
\ref{Theorem2.7}.

\begin{theorem}\label{theorem9.1}
Let $S_1=S_2=\infty$ (see \eqref{2.1} and \eqref{3.18}, and suppose that there
is $\delta>0$ such that
for $x\in \mathbb R$ the following inequality holds (see \eqref{2.2}):
\begin{equation}\label{9.1}
d(t)\ge \delta d(x)\qquad\text{for}\quad |t-x|\le d(x).
\end{equation}
Then, if condition  \eqref{1.9}  holds, problem \eqref{1.1}--\eqref{1.2} cannot
be correctly solvable
in $L_p(\mathbb R)$ for any $p\in[1,\infty].$
\end{theorem}

\begin{remark}\label{rem9.2}
According  to Example \ref{examp3.1}, the function $d(x)$ always satisfies
inequalities  \eqref{3.16}.
Inequality \eqref{9.1} slightly strengthens the a priori property \eqref{3.16},
and therefore
\thmref{theorem9.1} can be applied to a broad class of equations \eqref{1.1}.
\end{remark}

\begin{proof}[Proof of \thmref{theorem9.1}]
We consider the cases 1)\ $p=1;$\quad 2) \ $p\in(1,\iy);$\quad 3)\ $p=\infty$\
separately.
Assume that \begin{equation}\label{9.2}q(x)\to0 \quad \text{as}\quad
x\to\infty\end{equation}
 (the case $x\to-\infty$
can be treated in a similar way).

1)\ Let $p=1.$ Assume the contrary:  problem \eqref{1.1}--\eqref{1.2}  is
correctly solvable in
$L_1(\mathbb R).$
Then $M_1<\infty$ (see \eqref{2.15}), and for any $t\in \mathbb R$ we get
\begin{align*}
M_1&\ge\frac{1}{r(t)}\int_{-\infty}^t
\exp\left(-\int_\xi^t\frac{q(s)}{r(s)}ds\right)d\xi \\
&\ge  \frac{1}{r(t)}\int_{t-d(t)}^t\exp\left(-\int_\xi^t\frac{q(s)}{r(s)}ds\right)d\xi\nonumber\\
&\ge \frac{d(t)}{r(t)}
\exp\left(-\int_{t-d(t)}^{t+d(t)}\frac{q(s)}{r(s)}ds\right)
\end{align*}
which implies
\begin{align}
\frac{d(t)}{r(t)}\le  e^2M_1,\quad t\in \mathbb R.\label{9.3}
\end{align}
Let us integrate inequality \eqref{9.3} along the interval $[x-d(x),x+d(x)],$\
$x\in \mathbb R$ and use
\eqref{9.1} to get
\begin{align*}
2e^2M_1d(x)&\ge \int_{x-d(x)}^{x+d(x)}\frac{d(t)}{r(t)}dt\ge \delta
d(x)\int_{x-d(x)}^{x+d(x)}\frac{dt}{r(t)},\quad x\in \mathbb R
\end{align*}
and hence
\begin{align}
\sup_{x\in\mathbb R}\int_{x-d(x)}^{x+d(x)}\frac{dt}{r(t)}
\le2e^2\delta^{-1}M_1<\infty.\label{9.4}
\end{align}

Let $\varepsilon$ be a given positive number.
Then $q(x)\le\varepsilon$ for $x\ge c(\varepsilon)\gg1.$
Since $S_2=\infty,$ there is $x_0(\varepsilon)\ge c(\varepsilon)\gg1$ such that
$x-d(x)\ge
c(\varepsilon)$ for $x\ge x_0(\varepsilon)$ (see \eqref{3.27} and
\eqref{3.26}).
For $x\ge x_0(\varepsilon)$, we get
$$
  2=\int_{x-d(x)}^{x+d(x)}
 \frac{q(t)}{r(t)}dt\le\varepsilon
 \int_{x-d(x)}^{x+d(x)}\frac{dt}{r(t)}
$$
and hence
\begin{equation}\label{9.5}
\lim_{x\to\infty}\int_{x-d(x)}^{x+d(x)}\frac{dt}{r(t)}=\infty.
\end{equation}
This leads to a contradiction between \eqref{9.4} and \eqref{9.5}.

\s
2)\ Let $p\in(1,\infty).$
Assume the contrary: problem \eqref{1.1}--\eqref{1.2} is correctly solvable in
$L_p(\mathbb R)$. Then
$M_p<\infty$ (see \eqref{2.3}).
Denote $z_1(x)=x-d(x),$\ $z_2(x)=x+d(x)$.
Below we use \eqref{2.3}, \eqref{2.4}, \eqref{2.2} and \eqref{9.1}:
\begin{align*}
M_p&\ge\left[\int_{-\infty}^{z_1(x)}
\exp\left(-p\int_t^{z_1(x)}\frac{q(\xi)}{r(\xi)}d\xi\right)dt
\right]^{\frac{1}{p}}
\\ &\phantom{\ge }\cdot
\left[\int_{z_1(x)}^{\infty}\frac{1}{r(t)^{p'}}
\exp\left(-p'\int_{z_1(x)}^t\frac{q(\xi)}{r(\xi)}d\xi\right)dt
\right]^{\frac{1}{p'}} \\
&\ge\left[\int_{z_1(x)-d(z_1(x))}^{z_1(x)}
\exp\left(-p\int_t^{z_1(x)}\frac{q(\xi)}{r(\xi)}d\xi\right)dt\right]
^{\frac{1}{p}}
\\&\phantom{\ge }\cdot
\left[\int_{z_1(x)}^{z_2(x)}\frac{1}{r(t)^{p'}}
\exp\left(-p'\int_{z_1(x)}^t\frac{q(\xi)}{r(\xi)}d\xi
\right)dt\right]^{\frac{1}{p'}} \\
&\ge\exp\left(-\int_{z_1(x)-d(z_1(x))}^{z_1(x)+d(z_1(x))}
\frac{q(\xi)}{r(\xi)}d\xi\right)
\cdot\exp\left(- \int_{z_1(x)}^{z_2(x)}\frac{q(\xi)}{r(\xi)}d\xi\right) \\
&\phantom{\ge }\cdot d(z_1(x))^{\frac{1}{p}}
\left[\int_{z_1(x)}^{z_2(x)}\frac{dt}{rt)^{p'}}\right]^{\frac{1}{p}'} \\
&\ge e^{-4}\delta^{\frac{1}{p}}d(x)^{\frac{1}{p}}
\left[\int_{x-d(x)}^{x+d(x)}\frac{dt}{r(t)^{p'}}\right]^{\frac{1}{p}'}
\end{align*}
which implies
\begin{align}
\sup_{x\in \mathbb R}d(x)^{\frac{1}{p}}
\left[\int_{x-d(x)}^{x+d(x)}\frac{dt}{r(t)^{p'}}\right]^{\frac{1}{p'}}
\le\delta^{-\frac{1}{p}}e^4M_p<\infty,\quad x\in \mathbb R.\label{9.6}
\end{align}

On the other hand, let $\varepsilon>0$ be given.
Below, for $x\ge x_0(\varepsilon)$ (see the proof of~\eqref{9.5}) it holds
\begin{align*}
2=\int_{x-d(x)}^{x+d(x)}\frac{q(t)}{r(t)}dt
\le\varepsilon\int_{x-d(x)}^{x+d(x)}\frac{dt}{r(t)}
\le 2^{\frac{1}{p}}\varepsilon d(x)^{\frac{1}{p}}
\left[\int_{x-d(x)}^{x+d(x)}\frac{dt}{r(t)^{p'}}\right]^{\frac{1}{p'}}\\
\end{align*}
which implies
\begin{align}
\lim\limits_{x\to\infty}d(x)^{\frac{1}{p}}
\left[\int_{x-d(x)}^{x+d(x)}\frac{dt}{r(t)^{p'}}\right]^{\frac{1}{p'}}
=\infty.\label{9.7}
\end{align}
This leads to a contradiction between \eqref{9.6} and \eqref{9.7}.

\s
3)\ Let $p=\infty.$
Assume the contrary: problem \eqref{1.1}--\eqref{1.2} is correctly solvable in
$C(\mathbb R).$
We have $(x-d(x))\to\infty$ as $x\to\infty$ (see case 1) above).
Since
\begin{equation}
\begin{aligned}\label{9.8}
\int_{x-d(x)}^{\infty}\frac{1}{r(t)}
\exp&\left(-\int_{x-d(x)}^t\frac{q(\xi)}{r(\xi)}d\xi\right)dt  \\
&\ge\int_{x-d(x)}^{x+d(x)}\frac{1}{r(t)}
\exp\left(-\int_{x-d(x)}^t\frac{q(\xi)}{r(\xi)}d\xi\right)dt\\
&\ge \exp\left(-\int_{x-d(x)}^{x+d(x)}\frac{q(\xi)}{r(\xi)}d\xi\right)
\int_{x-d(x)}^{x+d(x)}\frac{dt}{r(t)} \\
&=e^{-2}\int_{x-d(x)}^{x+d(x)}\frac{dt}{r(t)}>0,
\end{aligned}
\end{equation}
From \eqref{9.8} and \eqref{2.18} it follows that
\begin{equation}\label{9.9}
\lim_{x\to\infty}\int_{x-d(x)}^{x+d(x)}\frac{dt}{r(t)}=0.
\end{equation}

On the other hand, let $\varepsilon>0$ be given.
Below, for $x\ge x_0(\varepsilon)$ (see the proof of \eqref{9.5}) it holds
$$
2=\int_{x-d(x)}^{x+d(x)}\frac{q(t)}{r(t)}dt
\le\varepsilon\int_{x-d(x)}^{x+d(x)}\frac{dt}{r(t)},
$$
which implies
\begin{equation}\label{9.10}
\lim_{x\to\infty}\int_{x-d(x)}^{x+d(x)}\frac{dt}{r(t)}=\infty.
\end{equation}
This leads to a contradiction between \eqref{9.9} and \eqref{9.10}.
\QED \end{proof}

\begin{theorem}\label{theorem9.3}
Suppose that the functions $r(x)$ and $q(x)$ satisfy the conditions
\begin{align}
\text{\bf 1)}\ S_1=\infty \ \text{ (see \eqref{2.1}) } \hspace{8cm} \nonumber\\
\text{\bf 2)}\ q_0>0 \text{ (see \eqref{1.16} and } \quad
\label{9.11} q(x)\to\infty\ \text{ as} \quad |x|\to\infty.
\hspace{2.57cm}
\end{align}
Then problem \eqref{1.1}--\eqref{1.2} is correctly solvable in $C(\mathbb R).$
\end{theorem}

\begin{proof}
Let us check that $A(x)\to0$ as $|x|\to\infty$ (see \eqref{6.1} and
\eqref{2.18}).
Fix $\varepsilon>0.$
Then there is an interval $(x_1,x_2)$ such that
\begin{equation}\label{9.12}
q(x)\ge 3\varepsilon^{-1}\qquad\text{for}\qquad x\notin (x_1,x_2).
\end{equation}
To estimate $A(x)$ for $x\notin (x_1,x_2)$, we consider the cases\ a)~$x\ge
x_2$ and \ b)~$x\le x_1$ separately.
In case~a) we have
\begin{align*}A(x)&=\int_x^\infty\frac{1}{r(t)}
\exp\left(-\int_x^t\frac{q(\xi)}{r(\xi)}d\xi\right)dt\le
\frac{\varepsilon}{3}\int_x^\infty\frac{q(t)}{r(t)}
\exp\left(-\int_x^t\frac{q(\xi)}{r(\xi)}d\xi\right)dt\\
&\le
\frac{\varepsilon}{3}<
\varepsilon.
\end{align*}
To estimate $A(x)$ in case  b), we write $A(x)$
for $x\le x_1$ in the form
\begin{equation}\begin{aligned}\label{9.13}
A(x)&=\int_{x}^{x_1}\frac{1}{r(t)}
\exp\left(-\int_x^t\frac{q(\xi)}{r(\xi)}d\xi\right)dt \\
&\phantom{= }+\int_{x_1}^{x_2}\frac{1}{r(t)}
\exp\left(-\int_x^t\frac{q(\xi)}{r(\xi)}d\xi\right)dt\\
&\phantom{= }+\int_{x_2}^\infty\frac{1}{r(t)}
\exp\left(-\int_x^t\frac{q(\xi)}{r(\xi)}d\xi\right)dt \\[0.3cm]
&:= A_1(x)+A_2(x)+A_3(x).
\end{aligned}
\end{equation}
We estimate each summand in \eqref{9.13} separately for $x\le x_1$:
\begin{align} \label{9.14}
\begin{split}
A_1(x)&=\int_x^{x_1}\frac{1}{r(t)}
\exp\left(-\int_x^t\frac{q(\xi)}{r(\xi)}d\xi\right)dt \\
&\le \frac{\varepsilon}{3}\int_x^{x_1}\frac{q(t)}{r(t)}
\exp\left(-\int_x^t\frac{q(\xi)}{r(\xi)}d\xi\right)dt  \\[0.2cm]
&\le\frac{\ve}{3}\\
A_2(x)&=\int_{x_1}^{x_2}\frac{1}{r(t)}
\exp\left(-\int_x^t\frac{q(\xi)}{r(\xi)}d\xi\right)dt \\
&\le \frac{1}{q_0}\int_{x_1}^{x_2}\frac{q(t)}{r(t)}
\exp\left(-\int_x^t\frac{q(\xi)}{r(\xi)}d\xi\right)dt \\
&\le \frac{1}{q_0}
\exp\left(-\int_x^{x_1}\frac{q(\xi)}{r(\xi)}d\xi\right).
\end{split}
\end{align}
Since $S_1=\infty,$ there is $x_0=x_0(\varepsilon)\ll x_1$ such that
$$
\frac{1}{q_0} \exp\left(-\int_x^{x_1}\frac{q(\xi)}{r(\xi)}d\xi\right)
\le\frac{\varepsilon}{3}
$$
for $x\le x_0$, we get $A_2(x)\le \frac{\varepsilon}{3}$ for $x\le  x_0.$
Finally, from
$$
A_3(x)=\int_{x_2}^\infty\frac{1}{r(t)}
\exp\left(-\int_x^t\frac{q(\xi)}{r(\xi)}d\xi\right)dt\le\int_{x_2}^\infty
\frac{q(t)}{r(t)}
\exp\left(-\int_{x_2}^t\frac{q(\xi)}{r(\xi)}d\xi\right)\frac{dt}{q(t)}
$$
it follows $A_3(x)\le\frac{\varepsilon}{3}$. Hence for $x\notin(x_0,x_2)$,
we have the estimate
$$
A(x)=A_1(x)+A_2(x)+A_3(x)
\le\frac{\varepsilon}{3}+\frac{\varepsilon}{3}+\frac{\varepsilon}{3}
=\varepsilon
$$
which implies $\lim_{ |x|\to\infty}A(x)=0.$
 It remains to refer to \thmref{Theorem2.5}.
\QED \end{proof}

\begin{proof}[Proof of \thmref{theorem1.2}]
To prove \thmref{theorem1.2}, we need some lemmas.
When stating them, we assume that the hypotheses of \thmref{theorem1.2} are
satisfied.
Below we often  use an obvious statement which, for convenience, is formulated
as a separate assertion.
  \end{proof}
\begin{lemma}\label{lemma9.4}
Let $\varphi(x)$ and $\psi(x)$ be positive and continuous functions for $x\in
R.$
If there exist a constant $c\in[1,\infty)$ and an interval $(x_1,x_2)$ such
that
\begin{equation}\label{9.15}
c^{-1}\psi(x)\le \varphi(x)\le c\psi(x)\qquad \text{for}\quad x\notin
(x_1,x_2),
\end{equation}
then equalities \eqref{9.15} remain true for all $x\in \mathbb R,$ possibly
after the replacement of $c$ with a bigger constant.
\end{lemma}

\begin{proof} The function $f(x)=\frac{\varphi(x)}{\psi(x)}$ is continuous and
positive for $x\in [x_1,x_2].$
Hence its minimum  $m$ and maximum $M$ on the segment $[x_1,x_2]$ are finite
positive numbers.
Let $c_1=\max\{c, m^{-1},M\}.$
Then $c_1^{-1}\psi(x)\le\varphi(x)\le c_1\psi(x)$ for $x\in \mathbb R.$
\QED \end{proof}

\begin{lemma}\label{lemma9.6}
Let $x\in \mathbb R$ be  given.
Let us construct a sequence $\{x_k\}_{k=-\infty}^\infty$ as follows:
\begin{equation}\label{9.16}
x_0=x,\quad x_{k+1}=x_k+b\frac{r(x_k)}{q(x_k)}\qquad \text{for}\quad
k=0,1,2,\dots
\end{equation}
\begin{equation}\label{9.17}
x_0=x,\quad x_{k-1}=x_k-b\frac{r(x_k)}{q(x_k)}\qquad \text{for}\quad
k=0,-1,-2,\dots
\end{equation}
Here $b$ is taken from \eqref{1.13}.
Then we have
\begin{equation}\label{9.18}
\lim_{k\to-\infty}x_k=-\infty,\qquad \lim_{k\to\infty}x_k=\infty.
\end{equation}
\end{lemma}

\begin{proof}
Both equalities from \eqref{9.18} are checked in a similar way.
Let us prove, for example, the second one.
Assume the contrary.
The sequence \eqref{9.16} is, by construction, monotone increasing.
If \eqref{9.18} does not hold, then there is $z<\infty$ such that $x_k<z$ for
$k\ge0.$
Then the sequence \eqref{9.16} has a limit $z_0\le z.$
Moreover,
\begin{align*}\infty&>z-x
\ge\sum_{k=0}^\infty(x_{k+1}-x_k)=b\sum_{k=0}^\infty\frac{r(x_k)}{q(x_k)},
\end{align*}
which implies $\lim_{k\to\infty}\frac{r(x_k)}{q(x_k)}=0$, in
contradiction to
$\lim_{k\to\infty}\frac{r(x_k)}{q(x_k)}=\frac{r(z_0)}{q(z_0)}\ne0.$
\QED
\end{proof}

\begin{lemma}\label{lemma9.7}
Let $\theta\in[0,b]$ (see \eqref{1.13}).
Denote
\begin{equation}\label{9.19}
\omega^{(+)}(x)=\left[x,x+\theta\frac{r(x)}{q(x)}\right],\quad
\omega^{(-)}(x)=\left[x-\theta\frac{r(x)}{q(x)},x\right],\quad x\in \mathbb R.
\end{equation}
Then for $x\notin (\alpha,\beta)$ (see \eqref{1.13}), the following
inequalities hold:
\begin{equation}\label{9.20}
\frac{\theta}{a^2}\le\int_{\omega^{(+)}(x)}\frac{q(t)}{r(t)}dt,
\quad\int_{\omega^{(-)}(x)}\frac{q(t)}{r(t)}dt\le
\theta a^2.
\end{equation}
\end{lemma}

\begin{proof}
Inequalities \eqref{9.20} follow from \eqref{1.13}.
For example,
\begin{align*}
\int_{\omega^{(+)}(x)}\frac{q(t)}{r(t)}dt
&=\int_{\omega^{(+)}(x)}\frac{q(t)}{q(x)}\cdot\frac{q(x)}{r(x)}
\cdot \frac{r(x)}{r(t)}dt
\ge\frac{1}{a^2}\cdot\frac{q(x)}{r(x)}\cdot\theta\frac{r(x)}{q(x)}
=\frac{\theta}{a^2} \\
\int_{\omega^{(+)}(x)}\frac{q(t)}{r(t)}dt
&=\int_{\omega^{(+)}(x)}\frac{q(t)}{q(x)}\cdot\frac{q(x)}{r(x)}
\cdot \frac{r(x)}{r(t)}dt
\le {a^2}\cdot\frac{q(x)}{r(x)}\cdot\theta\frac{r(x)}{q(x)}
= {\theta}{a^2}.  \tag*{\QED }
\end{align*}
\end{proof}

\begin{lemma}\label{lemma9.8}
We have (see \eqref{2.1} and \eqref{3.18})
\begin{equation}\label{9.21}
S_1=\int_{-\infty}^0\frac{q(t)}{r(t)}dt=\infty,\qquad
S_2=\int_0^\infty\frac{q(t)}{r(t)}dt=\infty.
\end{equation}
\end{lemma}

\begin{proof}
In \eqref{9.16} set $x_0=0.$
By \lemref{lemma9.6}, there is $k_0\gg 1$ such that the points $x_k$ for $k\ge
k_0$ are outside the interval
$(\alpha,\beta)$ from condition \eqref{1.13}.
Then by \lemref{lemma9.7} we have
$$
\infty\ge S_2\ge\int_{x_{k_0}}^\infty\frac{q(t)}{r(t)}dt
=\sum_{k=k_0}^\infty\int_{x_k}^{x_{k+1}}\frac{q(t)}{r(t)}dt
\ge\sum_{k=k_0}^\infty \frac{b}{a^2}=\infty ,
$$
which implies $S_2=\infty.$ The equality  $S_1=\infty$ can be checked in a
similar way.
\QED \end{proof}

\begin{lemma}\label{lemma9.9}
Let $a\ge1$,\ $b>0,$ and $3\gamma\le1$ (see \eqref{1.14}).
Then $b\ge a^2.$
\end{lemma}

\begin{proof}
Assume the contrary: $b<a^2$.
Then $3\le 3a^2\le e^{b/a^2}\le e$, a contradiction.
\QED \end{proof}

\begin{lemma}\label{lemma9.10}
For a given $x\in \mathbb R$, the equation in $d\ge 0$
\begin{equation}\label{9.23}
\int_{x-d}^{x+d}\frac{q(\xi)}{r(\xi)}d\xi=2
\end{equation}
has a unique positive solution $d=d(x).$
Moreover,
\begin{equation}\label{9.24}
\frac{1}{a^2}\frac{r(x)}{q(x)}\le d(x)\le a^2\frac{r(x)}{q(x)}\qquad
\text{for}\quad x\notin (\alpha,\beta),
\end{equation}
\begin{equation}\label{9.25}
c^{-1}\frac{r(x)}{q(x)}\le d(x)\le c\frac{r(x)}{q(x)}\qquad x\in \mathbb R.
\end{equation}
\end{lemma}

\begin{proof}
According to \eqref{9.21} and Example \ref{examp3.1}, we only have to prove
estimates \eqref{9.24} and \eqref{9.25}.
Let
$$
\eta_1(x)=\frac{1}{a^2}\frac{r(x)}{q(x)},\qquad x\notin (\alpha,\beta).
$$
Since $b\ge a^2\ge a^{-2}$ because of \lemref{lemma9.9}, from \eqref{9.20} it
follows that
$$\int_{x-\eta_1(x)}^{x+\eta_1(x)}\frac{q(t)}{r(t)}dt=\int_{x-\eta_1(x)}^x\frac{q(t)}{r(t)}dt+\int_x^{x+\eta_1(x)}
\frac{q(t)}{r(t)}dt\le a^2\frac{1}{a^2}+a^2\frac{1}{a^2}=2.$$
Hence $d(x)\ge\eta_1(x)$ by \lemref{lemma3.4}.
Let now
$$\eta_2(x)=a^2\frac{r(x)}{q(x)},\qquad x\notin (\alpha,\beta ).$$
Since $b\ge a^2 $ by \lemref{lemma9.9}, from \eqref{9.20} it follows that
$$
\int_{x-\eta_2(x)}^{x+\eta_2(x)}\frac{q(t)}{r(t)}dt=\int_{x-\eta_2(x)}^x
\frac{q(t)}{r(t)}dt+\int_x^{
x+\eta_2(x)}\frac{q(t)}{r(t)}dt\ge a^2\frac{1}{a^2}+a^2\frac{1}{a^2}=2.
$$
Hence $d(x)\le \eta_2(x)$ by \lemref{lemma3.4}, which implies \eqref{9.24}.
Since the function $d(x)$ is continuous and positive (see Example
\ref{examp3.1}), inequalities \eqref{9.25} follows
From \lemref{lemma9.4}
\QED \end{proof}

\begin{lemma}\label{em9.11}
Let $x\notin (\alpha,\beta)$ (see \eqref{1.13}), and let
$\{x_k\}_{k=-\infty}^\infty$ be the sequence from
\lemref{lemma9.6}.
Then the following inequalities hold:
\begin{alignat}{2}\label{9.26}
\int_{x_0}^{x_k}\frac{q(t)}{r(t)}dt&\ge\frac{bk}{a^2}
&\quad &\text{\rm for}\; x\ge\beta,\ k=0,1,2,\dots  \\
\label{9.27}
\int_{x_k}^{x_0}\frac{q(t)}{r(t)}dt&\ge\frac{b|k|}{a^2}
&\quad &\text{\rm for}\; x\le\alpha,\ k=0,-1,-2,\dots
\end{alignat}
\end{lemma}

\begin{proof}
If $k=0,$ relation \eqref{9.26} is obvious.
For $k\ge1$ it follows from Lemmas \ref{lemma9.6} and \ref{lemma9.7}:
$$
\int_{x_0}^{x_k}\frac{q(t)}{r(t)}dt
=\sum_{\ell=0}^{k-1}\int_{x_\ell}^{x_{\ell+1}}\frac{q(t)}{r(t)}dt
\ge\sum\limits_{\ell=0}^{k-1}\frac{b}{a^2}
=\frac{b}{a^2}k.
$$
Inequality \eqref{9.27} can be checked in a similar way.
\QED \end{proof}

\begin{lemma}\label{lemma9.12}
Let $x\notin (\alpha,\beta)$ (see \eqref{1.13}), and let
$\{x_k\}_{k=-\infty}^\infty$ be the sequence from
\lemref{lemma9.6}.
Then the following inequalities hold:
\begin{alignat}{2}\label{9.28}
a^{-k}&\le\frac{r(x_k)}{r(x_0)},\quad \frac{q(x_k)}{q(x_0)}\le a^k
&\quad &\text{\rm for}\; x\ge \beta,\ k=1,2,\dots \\
\label{9.29}
a^{-|k|}&\le\frac{r(x_k)}{r(x_0)},\quad \frac{q(x_k)}{q(x_0)}\le a^{|k|}
&\quad &\text{\rm for}\; x\le\alpha,\ k=-1,-2,\dots
\end{alignat}
\end{lemma}

\begin{proof}
Let $x\ge\beta.$
Then from \eqref{1.13} and \eqref{9.16}, it follows that
$$\frac{1}{a}\le \frac{r(x_\ell)}{r(x_{\ell-1})},\quad
\frac{q(x_{\ell})}{q(x_{\ell-1})}\le a\quad\text{for}\quad
\ell=1,2,\dots,k;\ k\ge 1.$$
After multiplying these inequalities, we obtain \eqref{9.28}.
Estimates \eqref{9.29} can be checked in a similar way.
\QED \end{proof}

\begin{lemma}\label{lemma9.13}
For $p\in[1,\infty)$ and $x\in \mathbb R,$ the following inequalities hold
with $c=c(p)$ (see \eqref{4.31}):
\begin{equation}\label{9.30}
c^{-1}\frac{r(x)}{q(x)}
\le I_p(x)=\int_{-\infty}^x
\exp\left(-p\int_t^x\frac{q(\xi)}{r(\xi)}d\xi\right)dt
\le c\frac{r(x)}{q(x)}.
\end{equation}
\end{lemma}
\begin{proof}
The proof of the lower bound in \eqref{9.30} is based on \lemref{lemma9.10}:
\begin{align*}
I_p(x)&=\int_{-\infty}^x
\exp\left(-p\int_t^x\frac{q(\xi)}{r(\xi)}d\xi\right)dt\\
&\ge\int_{x-d(x)}^x\exp\left(-p\int_t^x\frac{q(\xi)}{r(\xi)}d\xi\right)dt\\
&\ge d(x)\exp\left(-p\int_{x-d(x)}^{x+d(x)}
\frac{q(\xi)}{r(\xi)}d\xi\right)\\
&=d(x)e^{-2p}\\
&\ge c^{-1}\frac{r(x)}{q(x)}.
\end{align*}
To prove the upper bound in \eqref{9.30}, consider two separate cases:\
1)~$x\le \alpha;$\ 2)~$x\ge\beta.$
In case 1), we use below the sequence \eqref{9.17} and relations \eqref{9.27},
\eqref{9.29} and \eqref{1.14}:
\begin{align*}
I_p(x)&=\int_{-\infty}^x\exp\left(-p\int_t^x\frac{q(\xi)}{r(\xi)}d\xi\right)dt\\
&=\sum_{k=-\infty}^0\int\limits_{x_{k-1}}^{x_k}
\exp\left(-p\int_t^{x_0}\frac{q(\xi)}{r(\xi)}d\xi\right)dt  \\
&\le\sum_{k=-\infty}^0(x_k-x_{k-1})
\exp\left(-p\int_{x_k}^{x_0}\frac{q(\xi)}{r(\xi)}d\xi\right) \\
&\le  b\sum_{k=-\infty} ^0\frac{r(x_k)}{q(x_k)}
\exp\left(-p\frac{b}{a^2}|k|\right)  \\
&=b\frac{r(x_0)}{q(x_0)}
\sum_{k=-\infty}^0\frac{r(x_k)}{r(x_0)}\frac{q(x_0)}{q(x_k)}
\exp\left(-p\frac{b}{a^2}|k|\right) \\
&\le b\frac{r(x)}{q(x)}\sum_{k=-\infty}^0a^{2|k|}
\exp\left(-p\frac{b}{a^2}|k|\right)\\
&\le b\frac{r(x)}{q(x)}\sum_{k=-\infty}^0
\left(a^2\exp\left(-\frac{b}{a^2}\right)\right)^{|k|} \\
&\le b\frac{r(x)}{q(x)}\sum_{k=0}^\infty \frac{1}{3^k} \\
&=c\frac{r(x)}{q(x)}.
\end{align*}

Consider now case 2).
Let us write it down in the following form:
\begin{align} \label{9.31}
\begin{split}
I_p(x)&=\int_{-\infty}^x
\exp\left(-p\int_t^x\frac{q(\xi)}{r(\xi)}d\xi\right)dt \\
&=\exp\left(-p\int_0^x\frac{q(\xi)}{r(\xi)}d\xi\right)
\int_{-\infty}^x\exp\left(p\int_0^t\frac{q(\xi)}{r(\xi)}d\xi\right)dt\\
&=\exp\left(-p\int_0^x\frac{q(\xi)}{r(\xi)}d\xi\right)f(x)
\end{split}
\end{align}
where
$$
f(x)\doe\int_{-\infty}^x\exp\left(p\int_0^t\frac{q(\xi)}{r(\xi)}d\xi\right)dt,\quad
x\ge\beta.
$$
From \eqref{9.31} and case~1) above, it follows that the integral $I_p(x)$
exists for $x\ge \beta.$
Moreover, according to \eqref{9.21}, we have the equality
\begin{equation}\label{9.32}
\lim_{x\to\infty}f(x)
=\lim_{x\to\infty}\int_{-\infty}^x
\exp\left(p\int_0^t\frac{q(\xi)}{r(\xi)}d\xi\right)dt
=\infty.
\end{equation}

Let us define the integral
\begin{equation}\label{9.33}
I_p(x,\beta)=\int_\beta^x
\exp\left(-p\int_t^x\frac{q(\xi)}{r(\xi)}d\xi\right)dt,\quad x\ge\beta.
\end{equation}
Let us write down for $I_p(x,\beta)$ an analogue of the representation
\eqref{9.31}:
\begin{align}\label{9.34}
\begin{split}
I_p(x,\beta)&=\int_{\beta}^x
\exp\left(-p\int_t^x\frac{q(\xi)}{r(\xi)}d\xi\right)dt \\
& =\exp\left(-p\int_0^x\frac{q(\xi)}{r(\xi)}d\xi\right)
\int_\beta^x\exp\left(p\int_0^t\frac{q(\xi)}{r(\xi)}d\xi\right)dt \\
&=\exp\left(-p\int_0^x\frac{q(\xi)}{r(\xi)}d\xi\right)f_\beta(x),
\end{split}
\end{align}
where
$$
f_\beta(x)\doe\int_\beta^x
\exp\left(p\int_0^t\frac{q(\xi)}{r(\xi)}d\xi\right)dt,\quad x\ge \beta.
$$
Here, according to \eqref{9.21}, we have
\begin{equation}\label{9.35}
\lim_{x\to\infty}f_\beta(x)
=\lim_{x\to\infty}\int_\beta^x
\exp\left(p\int_0^t\frac{q(\xi)}{r(\xi)}d\xi\right)
dt=\infty.
\end{equation}
{}From \eqref{9.31}, \eqref{9.32}, \eqref{9.34}, \eqref{9.35} and L'H\^opital's
rule, it follows that
\begin{equation}\label{9.36}
\lim_{x\to\infty}\frac{I_p(x)}{I_p(x,\beta)}
=\lim_{x\to\infty}\frac{f(x)}{f_\beta(x)}=1.
\end{equation}
Let $m\ge\beta$ be such that for $x\ge m$, the following inequality holds:
\begin{equation}\label{9.37}
I_p(x)\le 2I_p(x,\beta),\qquad x\ge m\ge\beta.
\end{equation}

Consider the sequence \eqref{9.17}.
By \lemref{lemma9.6}, for $x\ge m\ge \beta$ there is $\ell\le0$  such that
\begin{equation}\label{9.38}
x_\ell\ge\beta,\qquad x_{\ell-1}\le\beta.
\end{equation}
Let us show that here one can choose $m$ so that for all $x\ge m$ the number
$\ell$ in inequalities \eqref{9.38}
satisfies the inequality $\ell\le -1.$
Assume the contrary.
Let $\{m_s\}_{s=1}^\infty$ be any monotone sequence increasing to infinity with
$m_1>\beta.$
By the assumption, for every $m_s,$ \ $s\ge 1,$ there is $x^{(s)}\ge m_s$ such
that
$$x^{(s)}-b\frac{r(x^{(s)})}{q(x^{(s)})}\le\beta.$$
This means that inequalities \eqref{1.13} can be extended to the interval
$[\beta,x^{(s)}]$ because
$$
[\beta,x^{(s)}]\subseteq
\left[x^{(s)}-b\frac{r(x^{(s)})}{q(x^{(s)})},x^{(s)}\right].
$$
Then
$$
\int_\beta^{x^{(s)}}\frac{q(t)}{r(t)}dt
=\int_\beta^{x^{(s)}}\frac{q(t)}{q(x^{(s)})}\cdot
\frac{q(x^{(s)})}{r(x^{(s)})}\cdot \frac{r(x^{(s)})}{r(t)}dt\le a^2b.
$$
Since here $x^{(s)}\ge m_s\to\infty$ as $s\to\infty,$ the integral $S_2$
converges (see \eqref{9.21}), a contradiction.
Therefore, in the sequel we choose $m$ big enough so that $m>\beta,$
\eqref{9.37} holds, and for all $x\ge m$ we always
have $\ell\le -1$ in \eqref{9.38}.

\s
To estimate $I_p(x,\beta)$, we use the sequence \eqref{9.17} and relations
\eqref{9.38}, \eqref{9.27},
\eqref{9.29} and~\eqref{1.14}:
\begin{align}\label{9.39}
\begin{split}
I_p(x,\beta)&=\int_\beta^x
\exp\left(-p\int_t^x\frac{q(\xi)}{r(\xi)}d\xi\right)dt \\
&=\sum_{k=\ell+1}^{0}\int_{x_{k-1}}^{x_k}
\exp\left( \int_t^{x_0}\frac{q(\xi)}{r(\xi)}d\xi\right)dt\\
&\phantom{= }+\int_\beta^{x_\ell}
\exp\left(-p\int_t^{x_0}\frac{q(\xi)}{r(\xi)}d\xi\right)dt \\
&\le\sum_{k=\ell+1}^0(x_k-x_{k-1})
\exp\left(-p\int_{x_k}^{x_0}\frac{q(\xi)}{r(\xi)}d\xi\right)\\
&\phantom{= }+(x_\ell-\beta)
\exp\left(-p\int_{x_\ell}^{x_0}\frac{q(\xi)}{r(\xi)}d\xi\right)\\
&\le b\sum_{k=\ell}^0\frac{r(x_k)}{q(x_k)}
\exp\left(-p\frac{b}{a^2}|k|\right)\\
&=b\frac{r(x_0)}{q(x_0)}
\sum_{k=\ell}^0\frac{r(x_k)}{r(x_0)}\cdot
\frac{q(x_0)}{q(x_k)}\exp\left(-p\frac{b}{a^2}|k|\right) \\
&\le b\frac{r(x)}{q(x)}\sum_{k=0}^{|\ell|}a^{2k}
\exp\left(-p\frac{b}{a^2}k\right)\\
&\le b\frac{r(x)}{q(x)}
\sum_{k=0}^\infty\left[a^2 \exp\left(-\frac{b}{a^2}\right)\right]^k \\
&=\frac{br(x)}{q(x)}\sum_{k=0}^\infty 3^{-k}=c\frac{r(x)}{q(x)}.
\end{split}
\end{align}

From \eqref{9.39} and \eqref{9.37}, we obtain estimates \eqref{9.30} for
$x\ge m.$
Thus inequalities \eqref{9.30} are proved for $x\notin (\alpha,m).$
To complete the proof of \eqref{9.30}, it remains to apply \lemref{lemma9.4}.
\QED \end{proof}

\begin{proof}[Proof of \thmref{theorem1.2} for $p=1$] {\tt Necessity}.

\s
Suppose that problem \eqref{1.1}--\eqref{1.2} is correctly solvable in
$L_1(\mathbb R).$
Then $r_0>0$ and $M_1<\infty$ because of \thmref{Theorem2.3} (see \eqref{2.14}
and \eqref{2.15}).
{}From \eqref{2.15} and \eqref{9.30}, it follows that
\begin{align*}
M_1=\sup_{x\in \mathbb R}\frac{1}{r(x)}
\int_{-\infty}^x\exp\left(-\int_t^x\frac{q(\xi)}{r(\xi)}d\xi\right)dt
=\sup_{x\in \mathbb R}\frac{1}{r(x)}I_1(x)
\ge
c^{-1}\sup_{x\in \mathbb R}\frac{1}{q(x)}
\end{align*}
which implies
$$
q_0=\inf_{x\in \mathbb R}q(x)\ge
\left(\sup_{x\in \mathbb R}\frac{1}{q(x)}\right)^{-1}
 \ge \frac{c^{-1}}{M_1}>0.
$$
\end{proof}

\begin{proof}[Proof of \thmref{theorem1.2} for  $p=1$] \ {\tt Sufficiency}.

\s
Since $S_1=\infty$ (see \eqref{9.21}) and $r_0>0$ (see \eqref{1.16}), in the
space $L_1(\mathbb R)$ correct solvability of
problem \eqref{1.1}--\eqref{1.2} is guaranteed by the inequality $M_1<\infty$
(see \thmref{Theorem2.3}).
Below we use \lemref{lemma9.13} and condition $q_0>0$ (see \eqref{1.16}) to
check this requirement:
\begin{align*}
M_1&=\sup_{x\in \mathbb R}\frac{1}{r(x)}\int_{-\infty}^x
\exp\left(-\int_t^x\frac{q(\xi)}{r(\xi)}d\xi\right)dt \\
&=\sup_{x\in \mathbb R}\frac{1}{r(x)}I_1(x)\\
&\le c\sup_{x\in \mathbb R}\frac{1}{r(x)}\cdot\frac{r(x)}{q(x)}\\
&=c\sup_{x\in \mathbb R}\frac{1}{q(x)}\le\frac{c}{q_0}<\infty. \tag*{\QED}
\end{align*}
\end{proof}

\begin{proof}[Proof of \thmref{theorem1.2} for $p\in(1,\infty)$] {\tt
Necessity}.

\s
Suppose that for some $p\in(1,\infty),$ problem \eqref{1.1}--\eqref{1.2} is
correctly solvable in $L_p(\mathbb R).$
Then $M_p<\infty$ by \thmref{Theorem2.1} (see \eqref{2.3}--\eqref{2.4}).
Let $x\in \mathbb R$ be arbitrary.
In the following relations, we use Lemmas \ref{lemma9.8} and \ref{lemma9.10}:
\begin{align}\label{9.40}
\begin{split}
\infty >M_p &\ge M_p(x) \\
&=\left[\int_{-\infty}^x
\exp\left(-p\int_t^x\frac{q(\xi)}{r(\xi)}d\xi\right)dt\right]^{\frac{1}{p}} \\
&\phantom{= } \cdot \left[\int_x^\infty\frac{1}{r(t)^{p'}}
\exp\left(-p'\int_x^t\frac{q(\xi)}{r(\xi)}d\xi\right)dt\right]^{\frac{1}{p'}}\\
&\ge\left[\int_{x-d(x)}^x
\exp\left(-p\int_t^x\frac{q(\xi)}{r(\xi)}d\xi\right)dt\right]^{\frac{1}{p}} \\
&\phantom{= } \cdot \left[\int_x^{x+d(x)}\frac{1}{r(t)^{p'}}
\exp\left(-p'\int_x^t\frac{q(\xi)}{r(\xi)}d\xi\right)dt\right]^{\frac{1}{p}'}\\
&\ge\left[\int_{x-d(x)}^x
\exp\left(-p\int_{x-d(x)}^x\frac{q(\xi)}{r(\xi)}d\xi\right)
dt\right]^{\frac{1}{p}}\\
&\phantom{= } \cdot\left[\int_{x}^{x+d(x)}\frac{1}{r(t)^{p'}}
\exp\left(-p'\int_x^{x+d(x)}\frac{q(\xi)}{r(\xi)}d\xi\right)dt\right]
^{\frac{1}{p'}}\\
&=\exp\left(-\int_{x-d(x)}^{x+d(x)}
\frac{q(\xi)}{r(\xi)}d\xi\right)d(x)^{\frac{1}{p}}\left[\int_x^{x+d(x)}
\frac{dt}{r(t)^{p'}}\right]^{\frac{1}{p'}}\\
&=e^{-2}d(x)^{\frac{1}{p}}
\left[\int_x^{x+d(x)}\frac{dt}{r(t)^{p'}}\right]^{\frac{1}{p'}}.
\end{split}
\end{align}
Below we assume that $x\notin (\alpha,\beta)$ and continue estimate
\eqref{9.40} using \eqref{1.13} and \eqref{9.24}:
\begin{align*}
\infty >e^2M_p&\ge d(x)^{\frac{1}{p}}
\left[\int_{x}^{x+d(x)}\frac{dt}{r(t)^{p'}}\right]^{\frac{1}{p'}}\\
&\ge\left(\frac{1}{a^2}\cdot
\frac{r(x)}{q(x)}\right)^{\frac{1}{p}}
\left(\int_{x}^{x+a^{-2}\frac{r(x)}{q(x)}}
\left(\frac{r(x)}{r(t)}\frac{1}{r(x)}\right)^{p'}dt\right]^{\frac{1}{p'}}\\
&\ge \frac{1}{a^3}\left(\frac{r(x)}{q(x)}\right)^{\frac{1}{p}}
\frac{1}{r(x)}\left(\frac{r(x)}{q(x)}\right)^{\frac{1}{p'}}\\
&=\frac{1}{a^3q(x)}
\end{align*}
which implies $q(x)\ge(e^2a^3M_p)^{-1}$ and hence $q_0>0$.

\s
Furthermore, by \thmref{Theorem2.1}, $A_{p'}$ is also finite (see
\eqref{2.6}--\eqref{2.7}).
In the following relations we use \eqref{1.13} and \lemref{lemma9.10}  for
$x\notin (\alpha,\beta):$
\begin{align*}
\infty>A_{p'}&\ge\int_{x-d(x)}^{x+d(x)}\frac{dt}{r(t)^{p'}}
\ge\int_{x-\frac{r(x)}{a^2q(x)}}^{x+\frac{r(x)}{a^2q(x)}}
\left(\frac{r(x)}{r(t)}\frac{1}{r(x)}\right)^{p'}dt
\ge
\frac{2}{a^{p'+2}}\frac{1}{r(x)^{p'-1}q(x)},
\end{align*}
which implies
$$
r(x)^{\frac{1}{p}}q(x)^{\frac{1}{p}'}
\ge\frac{1}{A_{p'}^{\frac{1}{p}'}}\frac{1}{a^3}
$$
and hence $\sigma_{p'}>0$.
\end{proof}

\begin{proof}[Proof of \thmref{theorem1.2} for  $p\in(1,\infty)$] \ {\tt
Sufficiency}.

\s
Below we need the following assertion.

\begin{lemma}\label{lemma9.14}
Let $\sigma_{p'}>0$ (see \eqref{1.16}).
Then the following inequalities hold:
\begin{equation}\label{9.41}
\frac{c^{-1}}{r(x)^{p'-1}q(x)}\le K_p(x)\le\frac{c}{r(x)^{p'-1}q(x)},
\quad x\in \mathbb R.
\end{equation}
Here $p'\ge1$, and
\begin{equation}\label{9.42}
 K_p(x)\doe \int_x^\infty\frac{1}{r(x)^{p'}}\exp\left(-p'\int_x^t
 \frac{q(\xi)}{r(\xi)}d\xi
\right)dt,\quad
x\in \mathbb R.
\end{equation}
\end{lemma}

\begin{proof} Let us first verify that the integral \eqref{9.42} converges for
every $x\in \mathbb R.$
Indeed, from the condition $\sigma_{p'}>0$ follows
$r(t)^{\frac{1}{p}}q(t)^{\frac{1}{p}'}\ge\sigma_{p'}>0$
and hence $r(t)^{p'-1}q(t)\ge\sigma_{p'}^{p'}$ for $t\in \mathbb R$,
which implies
\begin{equation}\label{9.43}
\frac{1}{r(t)^{p'}}\le\frac{1}{\sigma_p^{p'}}\cdot \frac{q(t)}{r(t)},\quad
t\in \mathbb R.
\end{equation}
From \eqref{9.43} and \eqref{9.21} for $x\in \mathbb R$ we now obtain
\begin{align*}
K_p(x)&=\int_x^\infty\frac{1}{r(t)^{p'}}
\exp\left(-p'\int_x^t\frac{q(\xi)}{r(\xi)}d\xi\right)dt \\
&\le\frac{1}{\sigma_{p'}^{p'}}\int_x^\infty\frac{q(t)}{r(t)}
\exp\left(-p'\int_x^t\frac{q(\xi)}{r(\xi)}d\xi\right)dt\\
&=\frac{1}{p'\sigma_{p'}^{p'}}
\left.\left[-\exp\left(-p'\int_x^t
\frac{q(\xi)}{r(\xi)}d\xi\right)\right|_x^\infty\right] \\
&=\frac{1}{p'\sigma_{p'}^{p'}}<\infty.
\end{align*}
To check the lower bound from \eqref{9.41}, we assume $x\notin (\alpha,\beta).$
Then according to \eqref{9.24} and \eqref{1.13}, we have
\begin{align*}
K_p(x)&\ge
\int_{x}^{x+d(x)}\frac{1}{r(t)^{p'}}
\exp\left(-p'\int_x^t\frac{q(\xi)}{r(\xi)}d\xi\right)dt\\
&\ge\exp\left(-p'\int\limits_{x-d(x)}^{x+d(x)}
\frac{q(\xi)}{r(\xi)}d\xi\right)\int_x^{x+d(x)}\frac{dt}{r(t)^{p'}}\\
&\ge e^{-2p'}\int_x^{x+r(x)/a^2q(x)}
\left(\frac{r(x)}{r(t)}\frac{1}{r(x)}\right)^{p'}dt \\
&\ge\frac{e^{-2p'}}{a^{p'+2}}\frac{1}{r(x)^{p'-1}q(x)}.
\end{align*}
Taking into account \lemref{lemma9.4}, the latter inequality gives the lower
bound from \eqref{9.41} for all $x\in \mathbb R.$
To prove the upper bound from \eqref{9.41}, we consider separate cases: \
1)~$x\ge\beta$ and \ 2)~$x\le\alpha.$
In case 1) we use below the sequence $\{x_k\}_{k=0}^\infty$ (see \eqref{9.16})
and inequalities \eqref{9.26},
\eqref{9.28}, \eqref{1.13} and \eqref{1.14}:
\begin{align*}
K_p(x)&=\int_x^\infty\frac{1}{r(t)^{p'}}
\exp\left(-p'\int_x^t\frac{q(\xi)}{r(\xi)}d\xi\right)dt\\
&=\sum_{k=0}^\infty\int_{x_k}^{x_{k+1}}\frac{1}{r(t)^{p'}}
\exp\left(-p'\int_{x_0}^t\frac{q(\xi)}{r(\xi)}d\xi\right)dt  \\
&\le \sum_{k=0}^\infty\int_{x_k}^{x_{k+1}}\frac{1}{r(t)^{p'}}
\exp\left(-p'\int_{x_0}^{x_k}\frac{q(\xi)}{r(\xi)}d\xi\right)dt \\
&\le\sum_{k=0}^\infty e^{-p'\frac{b}{a^2}k}\int_{x_k}^{x_{k+1}}
\left(\frac{r(x_k)}{r(t)}\frac{1}{r(x_k)}\right)^{p'}dt \\
&\le a^{p'}b\sum_{k=0}^\infty\frac{e^{-p'\frac{b}{a^2}k}}{r(x_k)^{p'-1}q(x_k)}\\
&=\frac{c_1}{r(x_0)^{p'-1}q(x_0)}\sum_{k=0}^\infty
\left(\frac{r(x_0)}{r(x_k)}\right)^{p'-1}
\left(\frac{q(x_0)}{q(x_k)}\right)e^{-p'\frac{b}{a^2}k}  \\
&\le \frac{c_2}{r(x)^{p'-1}q(x)}
\sum_{k=0}^\infty\left(ae^{-\frac{b}{a^2}}\right)^{kp'} \\
&\le\frac{c_3}{r(x)^{p'-1}q(x)}\sum_{k=0}^\infty\frac{1}{3^k}\\
&= \frac{c}{r(x)^{p'-1}q(x)}.
\end{align*}
Thus estimate \eqref{9.41}  holds for $x\ge \beta.$

\s
Consider case 2).
Let us introduce the following function:
\begin{equation}\label{9.44}
K_p(x,\alpha)=\int_x^\alpha\frac{1}{r(t)^{p'}}
\exp\left(-p'\int_x^t\frac{q(\xi)}{r(\xi)}d\xi\right)dt,\quad
x\le \alpha.
\end{equation}
Let $x\le m<\alpha$ (we shall choose $m$ later).
Then
\begin{align} \label{9.45}
\begin{split}
K_p(x)&=K_p(x,\alpha)
+\exp\left(-p'\int_x^\alpha\frac{q(\xi)}{r(\xi)}d\xi\right)K_p(\alpha)\\
&\le K_p(x,\alpha)+\frac{1}{p'\sigma_{p'}^{p'}}
\exp\left(-p'\int_x^\alpha\frac{q(\xi)}{r(\xi)}d\xi\right) \\
&=K_p(x,\alpha)\left\{1+\frac{c}{K_p(x,\alpha)}
\exp\left(-p'\int_x^\alpha\frac{q(\xi)}{r(\xi)}d\xi\right)\right\}.
\end{split}
\end{align}
\QED \end{proof}
Here we have
\begin{align}\label{9.45a}
\begin{split}
\exp\left(p'\int_x^\alpha\frac{q(\xi)}{r(\xi)}d\xi\right)
K_p(x,\alpha)&=\int_x^\alpha\frac{1}{r(t)^{p'}}
\exp \left(p'\int_t^\alpha\frac{q(\xi)}{r(\xi)}d\xi\right)dt\\
&\ge\int_x^\alpha\frac{dt}{r(t)^{p'}}\\
&\ge\int_m^\alpha\frac{dt}{r(t)^{p'}}.
\end{split}
\end{align}
Denote
$$
c(\alpha)=\int_{-\infty}^\alpha\frac{dt}{r(t)^{p'}},
\quad\delta(m)=\int_m^\alpha\frac{dt}{r(t)^{p'}}
$$
and choose $m$ as follows:
$$
m=\begin{cases}
\theta_1\quad &\text{if}\quad
c(\alpha)=\infty\quad\text{and}\quad
\int_{\theta_1}^\alpha\frac{dt}{r(t)^{p'}}=1 \\
\theta_2\quad&\text{if}\quad c(\alpha)<\infty
\quad\text{and}\quad\int_{\theta_2}^\alpha
\frac{dt}{r(t)^{p'}}=\frac{c(\alpha)}{2}.
\end{cases}
$$
With such a choice of $m,$ from \eqref{9.45} and \eqref{9.45a} it follows that
\begin{equation}\label{9.47}
K_p(x)\le cK_p(x ,\alpha),\quad x\le m,\quad
c=1+\frac{\delta(m)^{-1}}{p'\sigma_{p'}^{p'}}.
\end{equation}
Let now $x\le m.$
For $x_0=x$ consider the sequence \eqref{9.16}.
By \lemref{lemma9.6}, for $x\le m\le\alpha$, there is $\ell\ge0$ such that
\begin{equation}\label{9.48}
x_\ell\le\alpha,\qquad x_{\ell+1}>\alpha.
\end{equation}

Let us show that $m$ can be chosen so small that for all $x\le m$ the number
$\ell$ in inequalities \eqref{9.48}
satisfies the inequality $\ell\ge1.$
Assume the contrary.
Let $\{m_s\}_{s=1}^\infty$ be any monotone sequence  decreasing to   $-\infty$
with $m_1<\alpha.$
By the assumption, for every $m_s,$\ $s\ge 1,$ there is $x_s\le m_s$ such that
$$x_s+b\frac{r(x_s)}{q(x_s)}\ge \alpha.$$
This means that inequalities \eqref{1.13} can be extended to the interval
$[x_s,\alpha]$ because
$$
[x_s,\alpha]\subseteq\left[x_s,x_s+b\frac{r(x_s)}{q(x_s)}\right].
$$
Then
$$
\int_{x_s}^\alpha\frac{q(t)}{r(t)}dt
=\int_{x_s}^\alpha\frac{q(t)}{q(x_s)}\cdot\frac{q(x_s)}{r(x_s)}\cdot
\frac{r(x_s)}{r(t)}dt\le a^2b<\infty;
$$
and since here $x_s\le m_s\to-\infty$ as $s\to\infty,$ the integral $S_1$
converges (see \eqref{9.21}).
Contradiction.
Therefore, below $m$ is chosen so small that $m<\alpha,$ \eqref{9.47} holds,
and in inequalities \eqref{9.48} we always
have $\ell\ge 1.$

\s
When estimating $K_p(x,\alpha),$ we use sequences \eqref{9.16} and relations
\eqref{9.48}, \eqref{9.27},
\eqref{1.3}, \eqref{9.29} and \eqref{1.14}:
\begin{align}\label{9.49}
\begin{split}
K_p(x,\alpha)&=\int_x^\alpha\frac{1}{r(t)^{p'}}
\exp\left(-p'\int_x^t\frac{q(\xi)}{r(\xi)}d\xi\right)dt\\
&=\sum\limits_{k=0}^{\ell-1}\int_{x_k}^{x_{k+1}}\frac{1}{r(t)^{p'}}
\exp\left(-p'\int_{x_0}^t\frac{q(\xi)}{r(\xi)}d\xi\right)dt\\
&\phantom{= } +\int_{x_\ell}^\alpha\frac{1}{r(t)^{p'}}
\exp\left(-p'\int_{x_0}^t\frac{q(\xi)}{r(\xi)}dt\right)dt \\
& \le \sum_{k=0}^{\ell-1}\int_{x_k}^{x_{k+1}}\frac{1}{r(t)^{p'}}
\exp\left(-p'\int_{x_0}^{x_k}\frac{q(\xi)}{r(\xi)}d\xi\right)dt\\
&\phantom{= } +\int_{x_\ell}^\alpha\frac{1}{r(t)^{p'}}
\exp\left(-p'\int_{x_0}^{x_\ell}\frac{q(\xi)}{r(\xi)}d\xi\right)dt \\
&\le\sum\limits_{k=0}^{\ell-1}e^{-p'\frac{b}{a^2}k}
\int_{x_k}^{x_{k+1}}\left(\frac{r(x_k)}{r(t)}\frac{1}{r(x_k)}\right)
^{p'}dt \\
&\phantom{= } +e^{-p'\frac{b}{a^2}\ell}\int_{x_\ell}^\alpha
\left(\frac{r(x_\ell)}{r(t)}\frac{1}{r(x_\ell)}\right)^{p'}dt \\
&\le a^{p'}b\left[\sum_{k=0}^{\ell-1}
\frac{1}{r(x_k)^{p-1}q(x_k)}e^{-p'\frac{b}{a^2}k}\right]\\
&\phantom{= } +a^{p'}\frac{(\alpha-x_\ell)}{r(x_\ell)^{p'}}
e^{-p'\ell\frac{b}{a^2}}\\
&\le a^{p'}b\left[\sum_{k=0}^{\ell-1}\frac{1}{r(x_k)^{p-1}q(x_k)}
e^{-p'\frac{b}{a^2}k}\right]+a^{p'}
\frac{x_{\ell+1}-x_\ell}{r(x_\ell)^{p'}}e^{-p'\frac{b}{a^2}\ell} \\
&=a^{p'}b\sum_{k=0}^\ell
\frac{1}{r(x_k)^{p-1}q(x_k)}e^{-p'\frac{b}{a^2}k} \\
&=\frac{c_1}{r(x_0)^{p'-1}q(x_0)}\sum_{k=0}^\ell
\left(\frac{r(x_0)}{r(x_k)}\right)^{p'-1}
\left(\frac{q(x_0)}{q(x_k)}\right)e^{-p'\frac{b}{a^2}k}  \\
&\le c_1\frac{1}{r(x)^{p'-1}q(x)}\sum_{k=0}^\ell
\frac{a^{kp'}}{e^{p'\frac{b}{a^2}k}}\\
&\le\frac{c_1}{r(x)^{p'-1}q(x)}\sum_{k=0}^\infty\left
 (\frac{1}{3}\right)^{kp'}\\
&=\frac{c_2}{r(x)^{p'-1}q(x)}.
\end{split}
\end{align}
Thus the upper bound from \eqref{9.41} holds for $x\notin (m,\beta).$
To finish the proof of \eqref{9.41}, it remains to apply \lemref{lemma9.4}.
\QED \end{proof}

Let us now go  to the proof of the theorem.
Here we use \thmref{Theorem2.1}.
Since $S_1=\infty$ (see \eqref{9.21} and \eqref{2.5}), to apply
\thmref{Theorem2.1} it is enough to prove that
$M_p<\infty$ and $A_{p'}<\infty$ (see \eqref{2.3}--\eqref{2.4} and
\eqref{2.6}--\eqref{2.7}).
In the next estimate for $M_p$ we use Lemmas \ref{lemma9.13} and
\ref{lemma9.14} and condition \eqref{1.16}:
\begin{align*}
M_p&=
=\sup_{x\in \mathbb R}(I_p(x))^{\frac{1}{p}}(K_p(x))^{\frac{1}{p'}}\\
&\le c\sup_{x\in \mathbb R}\left(\frac{r(x)}{q(x)}\right)^{\frac{1}{p}}
\frac{1}{r(x)^{\frac{1}{p}}q(x)^{\frac{1}{p'}}}\\
&=c\sup_{x\in \mathbb R}\frac{1}{q(x)}
<\infty .
\end{align*}

To check the inequality $A_{p'}<\infty,$ let us first estimate the function
$A_{p'}(x)$ for $x\notin(\alpha,\beta).$
Below we use \lemref{lemma9.10}, \eqref{1.13} and \eqref{1.16}:
\begin{align}
A_{p'}(x)
\le\int_{x-a^2r(x)/q(x)}^{x+a^2r(x)/q(x)}\left(\frac{
r(x)}{r(t)}\frac{1}{r(x)}\right)^{p'}dt
\le 2a^{p'+2}\frac{1}{r(x)^{p'-1}q(x)}
\le\frac{c}{\sigma_{p'}^{p'}}.
\label{9.50a}
\end{align}
Note that the function $A_{p'}(x)$ is continuous for $x\in \mathbb R$ because so is
$d(x)$ (see Example \ref{examp3.1}).
Therefore, $A_{p'}(x)$ is bounded on $[\alpha,\beta].$
Together with \eqref{9.50a}, this leads to the inequality   $A_{p'}~<~\infty.$
Thus problem \eqref{1.1}--\eqref{1.2} is correctly solvable in $L_p(\mathbb R)$ for
$p\in(1,\infty)$ by \thmref{Theorem2.1}.
\QED

\renewcommand{\qedsymbol}{}
\begin{proof}[Proof of \thmref{theorem1.2} for $p=\infty$]\ {\tt Necessity.}

\s
Suppose that problem \eqref{1.1}--\eqref{1.2} is correctly solvable in $C(\mathbb R).$
Then equality \eqref{2.18} holds, and $S_1=\infty$ (see \eqref{9.21} and
\eqref{2.1}).
For $x\in \mathbb R$, \lemref{lemma9.10} implies
\begin{equation}\label{9.50}
\int_{x}^\infty\frac{1}{r(t)}
\exp\left(-\int_x^t\frac{q(\xi)}{r(\xi)}d\xi\right)dt
\ge\int_x^{x+d(x)}
\frac{1}{r(t)}\exp\left(-\int_x^t\frac{q(\xi)}{r\xi)}d\xi\right)dt>0.
\end{equation}
From \eqref{2.18} and
\eqref{9.50}, it follows that
\begin{equation}\label{9.51}
\lim_{|x|\to\infty}\int_x^{x+d(x)}\frac{1}{r(t)}
\exp\left(-\int_x^t\frac{q(\xi)}{r(\xi)}d\xi\right)dt=0.
\end{equation}
Furthermore, for $x\notin(\alpha,\beta)$, using \eqref{1.13}  and \eqref{9.24}
we obtain
\begin{equation}
\begin{aligned}\label{9.52}
\int_x^{x+d(x)}\frac{1}{r(t)}&
\exp\left(-\int_x^t\frac{q(\xi)}{r(\xi)}d\xi\right)dt \\
&\ge\exp\left(-\int\limits_{x-d(x)}^{x+d(x)}
\frac{q(\xi)}{r(\xi)}d\xi\right)\int_x^{x+d(x)}\frac{dt}{r(t)}\\
&=e^{-2}\int_x^{x+d(x)}\frac{dt}{r(t)}\\
&\ge e^{-2}\int_x^{x+r(x)/a^2q(x)}\frac{r(x)}{r(t)}\frac{dt}{r(x)}\\
&\ge\frac{e^{-2}}{a^3}\frac{1}{q(x)}>0.
\end{aligned}
\end{equation}
From \eqref{9.51} and \eqref{9.52} we get
$\lim\limits_{|x|\to\infty}\frac{1}{q(x)}=0$.
Hence $q(x)\to\infty$ as $|x|\to\infty.$
\end{proof}

\begin{proof}[Proof of \thmref{theorem1.2} for $p=\infty$]\ {\tt Sufficiency.}

\s
In this case the statement of the theorem is an obvious consequence of
\lemref{lemma9.8} and \thmref{theorem9.3}.
\QED
\end{proof}

\end{document}